\documentclass[reqno, 11pt]{amsart}
\usepackage[left=2cm, right=2cm, top=2cm, bottom=2cm]{geometry}
\usepackage[T1]{fontenc}
\usepackage[utf8]{inputenc}
\DeclareUnicodeCharacter{000B}{ }	
\usepackage[english]{babel}
\usepackage[babel]{csquotes}
\usepackage{amssymb}
\usepackage{bbm}
\usepackage{amsmath}
\usepackage{amsthm}
\usepackage{latexsym}
\usepackage{xcolor}
\usepackage{mathrsfs}
\usepackage{bbm}
\usepackage{verbatim}
\usepackage[colorlinks=true, pdfstartview=FitV, linkcolor=blue, citecolor=blue, urlcolor=blue]{hyperref}
\usepackage[style=numeric, sorting=nty, sortcites=true, backend=bibtex, maxbibnames=20, giveninits=true]{biblatex}
\usepackage{faktor}
\usepackage[shortlabels]{enumitem}
\numberwithin{equation}{section}
\newtheorem{thm}{Theorem}[section]
\newtheorem{cor}[thm]{Corollary}
\newtheorem{lem}[thm]{Lemma}
\newtheorem{prop}[thm]{Proposition}
\theoremstyle{definition}
\newtheorem{defin}[thm]{Definition}
\newtheorem{remark}[thm]{Remark}

\renewcommand{\b}[1]{\boldsymbol #1}
\renewcommand{\d}{{\mathrm d}} 
\renewcommand{\div}{\operatorname{div}} 
\newcommand{\ii}[4]{\int_{#1}^{#2} #3 \: \d#4} 
\newcommand{\pd}[2]{\dfrac{\partial #1}{\partial #2}} 
\newcommand{\tpd}[2]{\frac{\partial #1}{\partial #2}} 
\newcommand{\td}[2]{\dfrac{\d #1}{\d #2}} 
\newcommand{\dq}[1]{\partial_t^h #1} 
\newcommand\restr[2]{{
		\left.\kern-\nulldelimiterspace 
		#1 
		\vphantom{\big|} 
		\right|_{#2} 
}} 
\newcommand{\review}[1]{{\color{black}#1}}
\newcommand{\reviewF}[1]{{\color{black}#1}}

\AtEveryBibitem{%
	\clearfield{doi}%
	\clearfield{url}%
	\clearfield{issn}%
	\clearfield{isbn}%
	\ifentrytype{article}
	{}
	{\clearfield{eprint}}
} 
\renewbibmacro*{publisher+location+date}{%
	\ifentrytype{book}
	{\printlist{publisher}%
		\iflistundef{location}
		{\setunit*{\addcomma\space}}
		{\setunit*{\addcomma\space}}%
		\printlist{location}}
	{\printlist{location}%
		\iflistundef{publisher}
		{\setunit*{\addcomma\space}}
		{\setunit*{\addcomma\space}}%
		\printlist{publisher}}
	\setunit*{\addcomma\space}%
	\usebibmacro{date}%
	\newunit} 

\begin{document}
	\title[On a diffuse interface model for diblock copolymers interacting with an electric field]{On a diffuse interface model for diblock copolymers interacting with an electric field}
	
	\author{Helmut Abels}
	\address{Fakult\"{a}t f\"{u}r	 Mathematik, Universit\"{a}t Regensburg, 93040 Regensburg, Germany}
	\email{helmut.abels@ur.de}
	\author{Andrea Di Primio}
	\address{Dipartimento di Matematica,
		Politecnico di Milano, 20133 Milano, Italy}
	\email{andrea.diprimio@polimi.it}
	\author{Harald Garcke}
	\address{Fakult\"{a}t f\"{u}r	 Mathematik, Universit\"{a}t Regensburg, 93040 Regensburg, Germany}
	\email{harald.garcke@ur.de}
	
	\begin{abstract}
		We consider a diffuse interface model describing a ternary system constituted by a conductive diblock copolymer and a homopolymer acting as solvent. The resulting dynamics is modeled by two Cahn--Hilliard--Oono equations for the copolymer blocks, accounting for long-range interactions; a classical Cahn--Hiliard equation for the homopolymer and the Maxwell equation for the electric displacement field. A multiphase singular potential is employed in order to ensure physical consistency. First, we show existence of global weak solutions in two and three dimensions. Uniqueness of weak solutions is established in the constant mobility case, and a conditional result is given in the general case. Instantaneous regularization and long-time behavior are also investigated, the latter in the case of affine-linear electric permittivity, showing in particular that solutions converge to a single stationary state.
	\end{abstract}
	\maketitle
	
	\section{Introduction} \label{sec:intro} 
	Diffuse interface models are widely used to effectively model phase separation phenomena happening in multicomponent systems. In this framework, the interface between two different interacting species of a system is modeled as a region with positive thickness, rather than as a lower-dimensional manifold. Starting from the pioneering work of J. W. Cahn and J. E. Hilliard (see \cite{CahnHilliard58}), diffuse interface models, often coupled with other physical laws, have been extensively used to describe a variety of natural phenomena. In this contribution, we are interested in a Cahn--Hilliard type system arising in the realm of materials science, specifically concerning the dynamics of certain polymeric solutions. Polymer blends are mixtures of some macromolecules and may possess peculiar mechanical properties. Diblock copolymers, in particular, are composed of two different subunits, which we shall refer to as blocks. When a diblock copolymer is mixed with a homopolymer acting as solvent (as, e.g., water does), two different kinds of phase separation phenomena may occur: at the macrophase level, separation occurs between the homopolymer and the copolymer; on the other hand, at the microphase level, the two blocks of the block copolymer also exhibit segregation dynamics. The result of the interplay between these two processes is the creation of strikingly regular patterns at the mesoscopic spatial scale, which have been experimentally investigated (see, for example, \cite{Avalos16, Avalos18}) as well as numerically simulated (see also \cite{Li20,SodiniMartini21}). Moreover, the model proposed in \cite{Avalos16} has also been mathematically analyzed in \cite{DiPG22}, where existence and uniqueness of global weak solutions is shown, as well as their instantaneous regularization properties and their convergence to a unique stationary state. Here, we shall address similar issues for a model proposed in \cite{Shen21}. This model is a refined version of the aforementioned copolymer model: indeed, it takes into account possible interactions with an electric field. Maxwellian interactions of this kind bear significant importance in many applications, as, for instance, electromagnetic fields are able to play a pivotal role in self-assembling procedures involved in the creation of nanotechnologies (we refer to \cite{Hamley1} for an extensive treatment of the subject, see also \cite{Boker2003}). In \cite{Shen21}, a ternary system formed by the two blocks of a diblock copolymer and a homopolymer acting as solvent is modeled through a set of partial differential equations, obtained by the variational derivatives of a suitable energy functional. Let us briefly introduce how to derive this model. \\ Let $d \in \{2,3\}$ and let $\Omega \subset \mathbb{R}^d$ be a bounded domain with a sufficiently smooth boundary $\partial \Omega$. In the following, we partially carry over the notation employed in \cite{Shen21}, namely, we let the indices $A$ and $B$ denote the two blocks of the copolymer, while $S$ denotes the homopolymer acting as solvent. For the sake of convenience, let us introduce the two sets of indexes
	\[
	\mathcal{I} := \{A,\, B,\, S\}, \qquad \mathcal{J} := \{1,\,2,\,3\}
	\]
	and the bijection $\ell : \mathcal{J} \to \mathcal{I}$, defined in such a way that
	\[
	\ell(1) = A,\qquad \ell(2) = B, \qquad \ell(3) = S.
	\]
	For any $i \in \mathcal{I}$, we denote with $c_i = c_i(\b{x}, t)$ the corresponding order parameter, representing the local concentration of the respective substance, and the symbol $\b{c}$ denotes the vector whose $j$-th component is equal to $c_{\ell(j)}$ for all $j \in \mathcal{J}$. On account of this choice, it holds
	\[
	c_A + c_B + c_S = 1.
	\]
	We now give the precise structure of the energy functional through a general overview of the various energetic contributions. For precise assumptions, we refer the reader to the next section. The energy terms may be classified as follows:
	\begin{enumerate}[label=(\roman*)]\itemsep 0.5em
		\item \textbf{Decoupled energy contributions.} First of all, we have the energetic terms referring to a single species in the mixture, namely
		\[ 
		\mathcal{E}_1(\b{c}) := \ii{\Omega}{}{\sum_{i \in \mathcal{I}} \left[  \dfrac{\gamma_i}{2}|\nabla c_i|^2 + \theta_i\Psi(c_i) \right]\!}{x},
		\]
		where the potential energy density function $\Psi : [0,1] \to \mathbb{R}$ is given and the constants $\gamma_i$ and $\theta_i$ are strictly positive for all $i \in \mathcal{I}$. On account on the convexity assumptions holding on the function $\Psi$, which will be later precised, let us already point out that $\mathcal{E}_1$ is the convex part of the usual Ginzburg--Landau free energy functional appearing in the standard theory of diffuse interface models.
		\item \textbf{Local interaction terms.} Local interaction terms are modeled using a function of three variables 
		\[
		I: [0,1]^3 \to \mathbb{R},
		\]
		hence giving rise to the energy contribution
		\[
		\mathcal{E}_2(\b{c}) = \ii{\Omega}{}{I(c_A,\, c_B,\, c_S)}{x}.
		\]
		\item \textbf{Long-range interaction terms.} In the modeling of diblock copolymers, long-range interactions are often taken into account. In diffuse interface models, this is usually achieved through the non-local energy terms defining the so-called Ohta--Kawasaki energy functional (see, for instance, \cite{Choksi03, Imai99, Nishiura84} and references therein), whose gradient flow is the conserved Cahn--Hilliard--Oono system. These terms read
		\[
		\mathcal{E}_3(c_A,\, c_B) := \sum_{(i,j) \in \{A,B\}^2} \dfrac{\alpha_{ij}}{2}\int_\Omega \int_\Omega G(\b{x},\b{y})(c_i(\b{x}) - \overline{c_i})(c_j(\b{y}) - \overline{c_j}) \:\d y \: \d x,
		\]
		where an overline stands for the integral average over $\Omega$ and the constants $\alpha_{ij} \in\mathbb{R}$ for all couples $(i,j) \in \{A,B\}^2$, and satisfy $\alpha_{AB} = \alpha_{BA}$. Moreover, $G$ denotes the Green function associated to the negative Laplace operator with homogeneous Neumann boundary conditions. In this way, the structure of the resulting Cahn--Hilliard--Oono system ensures conservation of mass (this is not always the case when dealing with general forms of the Cahn--Hilliard--Oono equation, see \cite{GGM}).
		\item \textbf{Electrical energy.} Let $\b{E}$ denote the electric field, and let $\b{D}$ be the electric displacement field. The two quantities are related by the formula
		\[
		\b{D} = \epsilon(c_A,\,c_B)\b{E},
		\]
		where $\epsilon$ is the electrical permittivity of the diblock copolymer, which depends on the copolymer order parameters $c_A$ and $c_B$. The corresponding energy contribution is given by
		\[
		\mathcal{E}_4(c_A,\,c_B, \b E) := \dfrac{1}{2}\ii{\Omega}{}{\b{D} \cdot \b{E}}{\boldsymbol{x}}= \dfrac{1}{2}\ii{\Omega}{}{\epsilon(c_A,\,c_B)|\b{E}|^2}{x}.
		\]
		In the following, instead of solving the problem for the electric field $\b{E}$, we will focus on the related scalar electric potential $\Phi$, and thus set $\b{E} = \b{E}_0 - \nabla \Phi$, where $\b{E}_0$ is some external applied, and thus given, electric field.
	\end{enumerate}
	Collecting the above terms, we define the energy functional
	\[
	\mathcal{E}(\b{c},\, \Phi) := \mathcal{E}_1(\b{c}) + \mathcal{E}_2(\b{c}) + \mathcal{E}_3(c_A,\, c_B) + \mathcal{E}_4(c_A,\,c_B,\, \Phi).
	\] 
	As shown in \cite[Section 2]{Shen21}, computing the variational derivatives of this energy functional and owing to the Onsager principle, we arrive to the following set of partial differential equations holding in $\Omega \times (0,T)$ for some given final time $T > 0$:
	\begin{equation} \label{eq:strongform} \small
		\begin{cases}
			\pd{c_A}{t} = \div\left[ M_{AA}(\b{c})\nabla \mu_A + M_{AB}(\b{c})\nabla \mu_B + M_{AS}(\b{c})\nabla \mu_S  \right] , \\[0.2cm]
			\pd{c_B}{t} = \div\left[ M_{AB}(\b{c})\nabla \mu_A + M_{BB}(\b{c})\nabla \mu_B + M_{BS}(\b{c})\nabla \mu_S  \right] , \\[0.2cm]
			\pd{c_S}{t} = \div\left[ M_{AS}(\b{c})\nabla \mu_A + M_{BS}(\b{c})\nabla \mu_B + M_{SS}(\b{c})\nabla \mu_S  \right] , \\[0.2cm]
			\mu_A = -\gamma_A\Delta c_A + \theta_A\Psi'(c_A) + \pd{I}{c_A} + \alpha_{AA}\mathcal{N}(c_A - \overline{c_A}) + \alpha_{AB}\mathcal{N}(c_B - \overline{c_B}) + \dfrac{1}{2}\pd{\epsilon}{c_A}|\b{E}_0 - \nabla \Phi|^2 , \\[0.3cm] 
			\mu_B = -\gamma_B\Delta c_B + \theta_B\Psi'(c_B) + \pd{I}{c_B} + \alpha_{BA}\mathcal{N}(c_A - \overline{c_A}) + \alpha_{BB}\mathcal{N}(c_B - \overline{c_B}) + \dfrac{1}{2}\pd{\epsilon}{c_B}|\b{E}_0 - \nabla \Phi|^2 , \\[0.3cm]
			\mu_S = -\gamma_S\Delta c_S + \theta_S\Psi'(c_S) + \pd{I}{c_S}, \\[0.2cm]
			\div\left[ \varepsilon(c_A,\, c_B)(\b{E}_0 - \nabla \Phi)\right] = 0, 
		\end{cases}
	\end{equation}
	where the mobility matrix $M = M(\b{c}) = (M_{ij}(\b{c}))_{(i,j) \in \mathcal{I} \times \mathcal{I}}$ is given. Above, the operator $\mathcal{N}$ is the inverse of the negative Laplace operator with homogeneous Neumann boundary conditions acting on functions with zero spatial average. Let now $\b{\mu}$ denote the vector of chemical potentials. System \eqref{eq:strongform} is subject to the following boundary and initial conditions:
	\begin{equation} \label{eq:conditions}
		\begin{cases}
			(M(\b{c})\nabla\b{\mu}) \cdot \b{n} = \b{0} & \quad \text{on }\partial\Omega \times (0,T), \\
			\nabla c_i \cdot \b{n} = 0 & \quad \text{on }\partial\Omega \times (0,T),\quad  \forall \: i \in \mathcal{I},\\
			\Phi = 0 & \quad \text{on }\partial\Omega \times (0,T), \\
			c_i(\cdot, 0) = c_{0i}& \quad \text{in }\Omega,\quad  \forall \: i \in \mathcal{I}.
		\end{cases}
	\end{equation}
	\textbf{Plan of the paper.}
	The contents of this work are structured as follows. Section \ref{sec:notation} contains the necessary preliminary information, namely notation, functional setting, assumptions and statements of main results. The existence result is then proved throughout Section \ref{sec:existence}, while uniqueness properties are discussed in Section \ref{sec:uniqueness}. Section \ref{sec:regularization} concerns instantaneous regularization results, while the longtime behavior is analyzed in Section \ref{sec:longtime}.
	
	\section{Preliminaries and main results} \label{sec:notation}
	\subsection{Notation and functional setting} The following paragraphs are devoted to clarify the notation employed throughout the present work. Nonetheless, we first point out once and for all that we may sometimes identify the index sets $\mathcal{I}$ and $\mathcal{J}$ through the map $\ell$ and use their elements interchangeably, i.e., with slight abuse of notation, we may identify any indexed variable $a_{j}$ with $a_{\ell(j)}$ (and viceversa) for any $j \in \mathcal{J}$. Furthermore, the (possibly indexed) symbol $C > 0$ denotes a positive constant, possibly depending on some structural parameters of the problem that are explicited out, if necessary. Its value may change within the same argument without relabeling.
	\paragraph{\textit{Functional spaces.}}
	For any Banach space $X$, the topological dual of $X$ is denoted by $X^*$, and the duality product between $X^*$ and $X$ by $\langle \cdot, \cdot \rangle_{X^*,X}$. If $X$ is a Hilbert space, then its scalar product is denoted by $(\cdot, \cdot)_X$.  Let $\Omega \subset \mathbb{R}^d$ be a bounded domain and let $Z$ denote either $\mathbb{R}^m$ or $\mathbb{R}^{m \times n}$ for given positive integers $m$ and $n$. For any $p \in [1, +\infty]$ and $k \in \mathbb{N}$, the ordinary Sobolev space of order $(k,p)$ is denoted by $W^{k,p}(\Omega;Z)$, with the convention that $W^{0,p}(\Omega;Z) := L^p(\Omega;Z)$, i.e., the Lebesgue space consituted by $p$-integrable classes of almost-everywhere equal functions with values in the space $Z$. All the former spaces are endowed with their canonical structure of Banach space, and we denote their norms with 
	\[
	\| \cdot \|_{W^{k,p}(\Omega;Z)} \text{ or } \| \cdot \|_{L^p(\Omega;Z)},
	\]
	accordingly. In the Hilbert space case $p = 2$, we set $H^k(\Omega; Z):=W^{k,2}(\Omega; Z)$ and we equip it with its classical scalar product. Moreover, we set 
	\[
	(\cdot,\,\cdot)_{H^k(\Omega;Z)} :=(\cdot,\,\cdot)_{W^{k,2}(\Omega;Z)}, \qquad \| \cdot \|_{H^k(\Omega;Z)} := \| \cdot \|_{W^{k,2}(\Omega; Z)}.
	\]
	For the sake of convenience, we introduce further notation to label the following variational structures:
	\begin{align*}
		V := H^1(\Omega; Z) & \hookrightarrow H := L^2(\Omega; Z) \hookrightarrow V^* := (H^1(\Omega; Z))^*,\\
		V_0 := H^1_0(\Omega; Z) := \{v \in V,\, v = 0 \text{ on } \partial \Omega\} & \hookrightarrow H := L^2(\Omega; Z) \hookrightarrow V^*_0 := (H^1_0(\Omega; Z))^*, \\
		V_{(0)} := H^1_{(0)}(\Omega; Z) := \{ v \in V, \: \overline{v} = 0\} & \hookrightarrow H_{(0)} := \{ v \in H, \: \overline{v} = 0\} \hookrightarrow V^*_{(0)} := (H^1_{(0)}(\Omega; Z))^*.
	\end{align*}
	The three triplets are endowed with their classical Hilbert structures (in particular, $V_0$ and $V_{(0)}$ are endowed with the usual $H^1$-seminorm). Let us point out that all the embeddings above are dense and compact in both two and three spatial dimensions. We shall leave to the context to clarify the nature of the space $Z$, with the convention that, as usual, if $Z$ is omitted, then $Z = \mathbb{R}$. We shall denote norms in the three spaces $V,\, H$ and $V^*$ with the symbols  $\| \cdot \|_V, \| \cdot \|_H$ and $\| \cdot \|_{V^*}$, respectively (scalar products are also denoted following the same logic), and the same is also applied to the triplets involving zero-trace and zero-mean functions. Finally, for any (possibly unbounded) interval $[a,b] \subset \mathbb R$, any Banach space $X$ and $p \geq 1$, we also denote the usual Bochner spaces by $L^p(a,b; X)$ (for the local variant, we shall use $L^p_{\text{loc}}(a,b; X)$). Analogously, the space of $X$-valued continuous functions (resp. $\alpha$-H\"{o}lder continuous functions, for $\alpha \in (0,1]$) is denoted by $C^0([a,b]; X)$ (resp. $C^{0, \alpha}([a,b]; X)$).
	\paragraph{\textit{The Gibbs simplex.}}
	Consider the affine hyperplane of $\mathbb{R}^3$
	\[
	\Sigma := \left\{ \b{y} \in \mathbb{R}^3 : \sum_{i=1}^3 y_i = 1 \right\},
	\]
	and the so-called Gibbs simplex
	\[
	\b{G} := \left\{ \b{y} \in \mathbb{R}^3 : \sum_{i=1}^3 y_i = 1, \: y_i \geq 0 \text{ for } i = 1,\,2,\,3 \right\} \subset \Sigma,
	\]
	the latter representing the set of physically relevant values for the order parameter vector $\b c$.
	The tangent space to $\Sigma$ is given by the vector subspace of $\mathbb{R}^3$
	\[
	T\Sigma := \left\{ \b{y} \in \mathbb{R}^3 : \sum_{i=1}^3 y_i = 0 \right\},
	\]
	and we denote by $\b{P}$ be the orthogonal projector onto $T\Sigma$. The action of $\b{P}$ is described by
	\begin{equation} \label{eq:projector}
		(\b{P}\b{v})_i = \dfrac{1}{3}\sum_{j=1}^3 (v_i - v_j), \qquad i=1,\,2,\,3,
	\end{equation}
	for all $\b{v} \in \mathbb{R}^3$. Let us introduce the set
	\[
	X := \left\{ \b{\eta} \in V : \b{\eta}(\b{x}) \in \Sigma \text{ for almost every } \b{x} \in \Omega \right\},
	\]
	and, for any fixed $\b{m} = [m_1,\,m_2,\,m_3]^T \in \mathbb{R}^3$ such that $m_1 + m_2 + m_3 = 1$, we set 
	\[
	X^{\b{m}} := \left\{ \b{\eta} \in X: \overline{\eta_{i}} = m_i \text{ for all }i \in \{1,2,3\} \right\},
	\]
	whose tangent space is
	\[
	TX^{\b{m}} = \left\{ \b{\eta} \in V_{(0)} :  \b\eta(\b{x}) \in T\Sigma \text{ for almost every }\b{x} \in \Omega\right\}.
	\]
	Furthermore, we define the space of elements of $(H^1(\Omega;\mathbb{R}^3))^*$ being the annihilator of elements of $H^1(\Omega;\mathbb{R}^3)$ that are either constant or have equal components, namely
	\[
	Y := \left\{ \b{\eta} \in  V^* : \left\langle \b{\eta}, \b{v}\right\rangle_{V^*, V} = 0, \, \b{v} = h(\b{x})[1,1,1]^T \text{ for some } h \in H^1(\Omega) \text{ or } \b{v} = \b{e}_i \text{ for some } i \in \{1,2,3\} \right\},
	\]
	where $\{\b{e}_i\}_{i=1}^3$ denotes the canonical basis of $\mathbb{R}^3$.
	\paragraph{\textit{Weak Laplace-like operators.}} Consider the subset of 3-by-3 matrices
	\[
	\mathcal{M} := \left\{ \mathfrak{M} \in \text{Sym}_\mathbb{R}(3), \, \ker \mathfrak{M} = \operatorname{span}([1,1,1]^T),\, \mathfrak{M} \text{ is positive definite over }T\Sigma \right\}
	\]
	and let $F: \mathbb{R}^3 \to \mathcal{M}$ be measurable and such that 
	\begin{align*}
		F(\b{s})\b{\zeta} \cdot \b{\zeta} \geq \lambda_F \b{\zeta}\cdot \b{\zeta},& \qquad \forall\: \b{\zeta} \in T\Sigma, \\
		|F_{ij}(\b{s})| \leq F^*, & \qquad \forall \: (i,j) \in \mathcal{J} \times \mathcal{J}, 
	\end{align*}
	for some $\lambda_F > 0$, $F^* > 0$ and all $\b{s} \in \mathbb{R}^3$. Let now $\Omega \subset \mathbb{R}^d$, with $d \in \{2,3\}$ be given,  be a bounded domain with Lipschitz boundary. For any measurable function $\b{\varphi} : \Omega \to \mathbb{R}^3$, define the operator
	\begin{equation} \label{eq:generalA}
		A_{\b{\varphi}} : TX^{\b{m}} \to Y,\qquad \b{u} \mapsto \left( \b{v} \mapsto \ii{\Omega}{}{F(\b{\varphi})\nabla \b{u} : \nabla \b{v}}{x} \right),
	\end{equation}
	which is also invertible owing to the Lax--Milgram theorem. Its inverse will be denoted by $\mathcal{N}_{\b{\varphi}}$, and we therefore have
	\begin{equation} \label{eq:weakformA}
		(F(\b{\varphi})\nabla\mathcal{N}_{\b{\varphi}}\b{\eta},\, \nabla \b{v})_H = \langle \b{\eta}, \b{v} \rangle_{V^*,V}
	\end{equation}
	for all $\b{\eta} \in Y$ and $\b{v} \in V$. The following is a basic result concerning the operators $A_{\b{\varphi}}$.
	\begin{lem} \label{lem:controlA}
		If every component of $F(\b{s})$ is differentiable and there exists $C_F > 0$ such that
		\[
		\left| \pd{F_{ij}}{s_k}\right| \leq C_F, \qquad \forall \: (i,j) \in \mathcal{J} \times \mathcal{J}, \quad \forall \: k \in \mathcal{J},
		\]
		then it holds
		\[
		\|\nabla\mathcal{N}_{\b{\varphi}_1}\b{\eta} - \nabla \mathcal{N}_{\b{\varphi}_2}\b{\eta}\|_H \leq \dfrac{C_F}{\lambda_F}\|\b{\varphi}_1 - \b{\varphi}_2\|_{L^4(\Omega;\mathbb{R}^3)}\|\nabla \mathcal{N}_{\b{\varphi}_2} \b{\eta}\|_{L^4(\Omega;\mathbb{R}^{d \times 3})},
		\]
		for all $\b{\varphi}_1,\,\b{\varphi}_2 \in L^4(\Omega;\mathbb{R}^3)$ and $\b{\eta} \in Y$.
	\end{lem} 
	\begin{proof}
		From \eqref{eq:weakformA}, we have that for any fixed $\b{\eta} \in Y$,
		\[
		(F(\b{\varphi}_1)\nabla\mathcal{N}_{\b{\varphi}_1}\b{\eta} - F(\b{\varphi}_2)\nabla\mathcal{N}_{\b{\varphi}_2}\b{\eta},\, \nabla \b{v})_H = 0,
		\]
		or equivalently
		\[
		(F(\b{\varphi}_1)[\nabla\mathcal{N}_{\b{\varphi}_1}\b{\eta} - \nabla \mathcal{N}_{\b{\varphi}_2}\b{\eta}], \nabla \b{v})_H = - ([F(\b{\varphi}_1) - F(\b{\varphi}_2)]\nabla\mathcal{N}_{\b{\varphi}_2}\b{\eta},\, \nabla \b{v})_H 
		\] 
		for all $\b{v} \in V$. Choosing $\b{v} = \mathcal{N}_{\b{\varphi}_1}\b{\eta} - \mathcal{N}_{\b{\varphi}_2}\b{\eta}$, we have,
		\[
		\lambda_F\|\nabla\mathcal{N}_{\b{\varphi}_1}\b{\eta} - \nabla \mathcal{N}_{\b{\varphi}_2}\b{\eta}\|_H^2 \leq C_F\|\nabla \mathcal{N}_{\b{\varphi}_1}\b{\eta} - \nabla \mathcal{N}_{\b{\varphi}_2}\b{\eta}\|_H\|\b{\varphi}_1 - \b{\varphi}_2\|_{L^4(\Omega;\mathbb{R}^3)}\|\nabla \mathcal{N}_{\b{\varphi}_2} \b{\eta}\|_{L^4(\Omega;\mathbb{R}^{d \times 3})}
		\]
		as a straightforward application of the H\"{o}lder inequality. The result then follows immediately.
	\end{proof} \noindent
	For $\b L \in Y$, we also introduce the norm
	\[ \| \b{L} \|_{*, \b{\varphi}} := \sqrt{(F(\b{\varphi})\nabla\mathcal{N}_{\b{\varphi}}\b{L}, \nabla\mathcal{N}_{\b{\varphi}}\b{L} )}_H,  \]
	which is equivalent to the usual $V_{(0)}^*$ norm in $Y$. Let us now give special notation for the relevant case of $F$ being a constant function, and hence identifiable with some element of $\mathcal{M}$. Let $\mathfrak{M} \in \mathcal{M}$ and assume $F \equiv \mathfrak{M}$. Then, as before, we shall define the linear operator
	\begin{equation*}
		A_\mathfrak{M}:  TX^{\b{m}} \to Y,\qquad \b{u} \mapsto \left( \b{v} \mapsto \ii{\Omega}{}{\mathfrak{M}\nabla \b{u} : \nabla \b{v}}{x} \right).
	\end{equation*}
	It is readily shown, owing to the Lax--Milgram theorem, that $A_\mathfrak{M}$ is an isomorphism. Therefore, the inverse operator $\mathcal{N}_\mathfrak{M} := A^{-1}_\mathfrak{M}$ is well defined. Once again, we consider the family of equivalent norms
	\[ \| \b{L} \|_{*, \mathfrak{M}} := \sqrt{( \mathfrak{M}\nabla\mathcal{N}_\mathfrak{M}\b{L}, \nabla\mathcal{N}_\mathfrak{M}\b{L} )_H} \]
	for any $\b{L} \in Y$. Finally, if $\mathfrak{M}$ is instead a positive definite matrix over the whole space $\mathbb{R}^3$, then the analogous version of $A_\mathfrak{M}$ actually defines an isomorphism between $V_{(0)}$ and its dual $V_{(0)}^*$. In particular, we shall be interested only in the case $\mathfrak{M} = I$, where $I$ is the identity matrix of dimension 3. The associated operator $A := A_I$ is the Riesz map of the Hilbert space $V_{(0)}$. Its inverse is denoted by $\mathcal{N}$ and, furthermore, we define the equivalent norm 
	\[ \| \b{L} \|_{*} := \|\nabla\mathcal{N}\b{L} \|_H = \sqrt{(\nabla\mathcal{N}\b{L}, \nabla\mathcal{N}\b{L})_H}  \]
	for any $\b{L} \in V_{(0)}^*$.
	\subsection{Basic assumptions.} The following assumptions are valid throughout the present work:
	\begin{enumerate}[label=\textbf{(A\arabic*)}] \itemsep 0.5em
		\item \label{hyp:potential} The function $\Psi: [0,1] \to \mathbb{R}$ belongs to the space $C^0([0,1]) \cap C^4((0,1])$. Furthermore, we assume the following properties: \vspace{0.1cm} 
		\begin{enumerate}[label = (\roman*)]
			\item \label{hyp:potpositive} the function $\Psi$ satisfies
			\[
			\Psi(s) \geq 0, \qquad \forall \: s \in [0,1];
			\]
			\item \label{hyp:potder1} the first derivative of $\Psi$ satisfies
			\[
			\lim_{s \to 0^+} \Psi'(s) = -\infty,
			\]
			and we denote by $\varepsilon_0$ some positive value such that
			\[
			\Psi'(s) \leq 0, \qquad \forall \: s \in (0,\varepsilon_0);
			\]	
			\item \label{hyp:potder1-2} there exists a constant $C_\Psi > 0$ such that
			\[
			\Psi'(s) \leq C_\Psi(1 + \Psi(s)), \qquad \forall \: s \in (0,1);
			\]	 
			\item \label{hyp:potder2} there exists a constant $\Theta > 0$ such that
			\[
			\Psi''(s) \geq \Theta > 0, \qquad \forall \: s \in (0,1);
			\]
			\item \label{hyp:potder4-2} there exists $c_\Psi > 0$ and $\varepsilon_1 \in (0,1)$ such that
			\[
			\Psi^{(3)}(s) \leq 0, \qquad \Psi^{(4)}(s) \geq c_\Psi, \qquad \forall \: s \in (0,\varepsilon_1);
			\]
			\item \label{hyp:potder4} there exists a constant $\varepsilon_2 \in (0,1)$ such that $\Psi^{(4)}$ is non-increasing in $(0, \varepsilon_2)$.
		\end{enumerate}
		Moreover, we extend the function $\reviewF{\Psi}$ with value $+\infty$ for every $s < 0$. For $s > 1$, we extend $\Psi$ in such a way that $\Psi \in C^1([1,+\infty))$ and both $\reviewF{\Psi}$ and $\reviewF{\Psi'}$ are bounded for $s \geq 1$.
		\item \label{hyp:interaction} The function $I: [0,1]^3 \to \mathbb{R}$ belongs to the space $C^2([0,1]^3)$.
		\item \label{hyp:mobility} The matrix-valued function $M : \mathbb{R}^3 \to \mathcal{M}$ and we set 
		\[
		M(\b{s}) = \begin{bmatrix}
			M_{AA}(\b{s}) & M_{AB}(\b{s}) & M_{AS}(\b{s}) \\
			M_{BA}(\b{s}) & M_{BB}(\b{s}) & M_{BS}(\b{s}) \\
			M_{SA}(\b{s}) & M_{SB}(\b{s}) & M_{SS}(\b{s})
		\end{bmatrix} \qquad \forall \: \b s \in \mathbb{R}^3,
		\]
		with the understanding that $M_{ij} = M_{ji}$ for all $(i,j) \in \mathcal{I} \times \mathcal{I}$. Moreover, $M_{ij} : \mathbb{R}^3 \to \mathbb{R}$ is of class $C^0(\mathbb{R}^3)$ for all $(i,j) \in \mathcal{I} \times \mathcal{I}$ and there exist a strictly positive constant $C_M$ such that 
		\begin{align*}
			|M_{ij}(\b{s})| \leq C_M & \qquad \forall \: \b{s} \in \mathbb{R}^3, \quad \forall \: (i,j) \in \mathcal{I} \times \mathcal{I}.
		\end{align*}
		Let $\b{\xi}$ denote the vector $[1,1,1]^T$. Then, we recall that $\b{\xi} \in \ker M(\b{s})$ for all $\b{s} \in \mathbb{R}^3$ and that $M(\b{s})$ is uniformly positive definite over $T\Sigma$, namely there exists $\lambda_0 > 0$ such that
		\[
		M(\b{s})\b{\zeta} \cdot \b{\zeta} \geq \lambda_0\b{\zeta} \cdot \b{\zeta}, \qquad \forall \: \b{\zeta} \in T\Sigma, \quad \forall \: \b{s} \in \mathbb{R}^3.
		\]
		\item \label{hyp:epsilon} The function $\epsilon : \mathbb{R}^2 \to \mathbb{R}$ is of class $C^1(\mathbb{R}^2)$. Furthermore, there exist two constants $\epsilon_*$ and $\epsilon^*$ such that
		\[
		0 < \epsilon_* \leq \epsilon(\b{s}) \leq \epsilon^*, \qquad \forall \: \b{s} \in \mathbb{R}^2,
		\]
		and the derivatives of $\epsilon$ vanish outside a compact set containing the set
		\[
		\{ (x,y) \in \mathbb{R}^2: x \geq 0,\,y \geq 0,\, x+y\leq 1\},
		\]for instance we set
		\[
		\operatorname{supp}\pd{\epsilon}{s_i} \subset \{\b x \in \mathbb{R}^2,\, |\b x| \leq 2\}, \quad \forall \: i \in \{1,2\}.
		\]
		\item \label{hyp:data} The structural data and parameters satisfy the following properties: \vspace{0.1cm}
		\begin{enumerate}[label = (\roman*)]
			\item \label{hyp:domain} the domain $\Omega \subset \mathbb{R}^d$ is bounded, with a boundary $\partial \Omega$ of class $C^{1,1}$;
			\item \label{hyp:param} the scalar parameters satisfy
			\[
			\gamma_i > 0, \qquad \theta_i > 0, \qquad \alpha_{jk} \in \mathbb{R}, \qquad \alpha_{AB} = \alpha_{BA},
			\]
			for all $i \in \mathcal{I}$, and $(j,k) \in \{A,B\}^2$;
			\item \label{hyp:initial} the initial state $\b{c}_0 = [c_{0A}, \: c_{0B}, \: c_{0S}]^T \in X$. Moreover, $\Psi(c_{0i}) \in L^1(\Omega)$ for all $i \in \mathcal{I}$;
			\item \label{hyp:E0} the external electric field $\b{E}_0 \in L^r(\Omega; \mathbb{R}^d)$ for some $r \geq 2$. Furthermore, we denote by $\Phi_0 \in V_0$ the unique weak solution to the problem
			\[
			\begin{cases}
				\div(\epsilon(c_{0A}, c_{0B})(\b{E}_0 - \nabla \Phi_0)) = 0 & \quad \text{ in }\Omega, \\
				\Phi = 0 & \quad \text{ on } \partial\Omega.
			\end{cases}
			\]
		\end{enumerate}
		Finally, we set 
		\begin{align*}
			\theta_* := \min_{i \in \mathcal{I}} \theta_i, \qquad \theta^* := \max_{i \in \mathcal{I}} \theta_i, \\ 
			\gamma_* := \min_{i \in \mathcal{I}} \gamma_i, \qquad \gamma^* := \max_{i \in \mathcal{I}} \gamma_i,  
		\end{align*}
		all of which are strictly positive, and define $\Gamma := \begin{bmatrix}
			\gamma_A & 0 & 0 \\
			0 & \gamma_B & 0 \\
			0 & 0 & \gamma_S
		\end{bmatrix}$.
	\end{enumerate}
	\begin{remark}
		Let us briefly comment on the above assumptions. First, observe that the thermodynamically relevant choice for the function $\Psi$, namely the logarithmic part of the Flory--Huggins potential $\Psi(s) = s \log s + e^{-1}$ (cfr. \cite{Flory42, Huggins41}), satisfies Assumption \ref{hyp:potential}. Indeed, it is sufficient to observe that its first four derivatives are
		\[
		\Psi'(s) = 1 + \log s, \qquad \Psi''(s) = \frac{1}{s}, \qquad \Psi^{(3)}(s) = -\dfrac{1}{s^2}, \qquad \Psi^{(4)}(s) = \dfrac{2}{s^3}.
		\]
		Assumption \ref{hyp:interaction} enables the modeling of local interactions with a wide range of functions. Usually, a typical choice is given by some high-order multivariate polynomial (refer for example to \cite{Avalos16, Shen21, Engblom13}), of which the simplest example is given by
		\[
		I(\b{s}) = s_1s_2 + s_2s_3 + s_1s_3, \qquad \forall \: \b{s} = [s_1,s_2,s_3]^T \in \mathbb{R}^3.
		\]
		Nonetheless, it is meaningful not to assume a too restrictive growth condition on the local interaction function $I$, as cubic or general higher-order polynomials are physically relevant choices for $I$. The matrix-valued function $M$ satisfies standard assumptions allowing for a wide class of mobility functions (see, for instance, \cite{Barrett99}). In particular, the property $\b{\xi} \in \ker M(\b{s})$ for all $\b{s} \in \mathbb{R}^3$ ensures conservation of mass. The electric permittivity complies with assumptions that are inspired from \cite{Garcke05}. Since the order parameter $\b{c}$ turns out to take values in the Gibbs simplex $\b{G}$, the assumption is not restrictive. Finally, the assumptions on the known parameters listed in \ref{hyp:data} are standard in the context of Cahn--Hilliard type equations.
	\end{remark} \noindent
	\subsection{Main results.} Before stating the main results of this work, we exploit the projection $\b{P}$ over the tangent space $T\Sigma$ to restate the problem in an equivalent alternative form. Defining the projected chemical potential vector
	\[
	\b{w} := \b{P}\b{\mu},
	\]
	and observing that
	\begin{equation} \label{eq:mudecomp}
		\b\mu = \b w + \dfrac{1}{3}(\b{\mu} \cdot \b\xi)\b\xi,
	\end{equation}
	we recast \eqref{eq:strongform}-\eqref{eq:conditions} in terms of $\b{w}$. Let us first show the following basic result.
	\begin{lem}
		It holds the equality $\div(M(\b{c}) \nabla \b{\mu}) = \div(M(\b{c}) \nabla \b P\b{\mu})$.
	\end{lem}
	\begin{proof}
		By \eqref{eq:mudecomp}, it amounts to show that $\div(M(\b{c}) \nabla [(\b{\mu} \cdot \b\xi)\b\xi]) = 0$. By components, for any $i \in \{1,2,3\}$, we find
		\[
		\begin{split}
			[\div(M(\b{c}) \nabla [(\b{\mu} \cdot \b\xi)\b\xi])]_i & = \sum_{j=1}^3 \pd{}{x_j} \sum_{k = 1}^3 M_{\ell(i)\ell(k)}(\b{c}) \left[ \sum_{m =1}^3 \pd{\mu_{\ell(m)}}{x_j} \right] = 0,
		\end{split}
		\]
		by the properties of $M(\b c)$, and the proof is complete.
	\end{proof} \noindent
	Therefore, we can recast the problem as follows:
	\begin{equation} \label{eq:strongform2} \small
		\begin{cases}
			\pd{\b{c}}{t} = \div(M(\b{c}) \nabla \b w) , \\[0.2cm]
			\b w= \b P \begin{bmatrix}
				-\gamma_A\Delta c_A + \theta_A\Psi'(c_A) + \pd{I}{c_A} + \alpha_{AA}\mathcal{N}(c_A - \overline{c_A}) + \alpha_{AB}\mathcal{N}(c_B - \overline{c_B}) + \dfrac{1}{2}\pd{\epsilon}{c_A}|\b{E}_0 - \nabla \Phi|^2 \\[0.3cm]
				-\gamma_B\Delta c_B + \theta_B\Psi'(c_B) + \pd{I}{c_B} + \alpha_{BA}\mathcal{N}(c_A - \overline{c_A}) + \alpha_{BB}\mathcal{N}(c_B - \overline{c_B}) + \dfrac{1}{2}\pd{\epsilon}{c_B}|\b{E}_0 - \nabla \Phi|^2 \\[0.3cm]
				-\gamma_S\Delta c_S + \theta_S\Psi'(c_S) + \pd{I}{c_S}
			\end{bmatrix} , \\[0.3cm] 
			\div\left[ \varepsilon(c_A,\, c_B)(\b{E}_0 - \nabla \Phi)\right] = 0, 
		\end{cases}
	\end{equation}
	with boundary and initial conditions given by
	\begin{equation} \label{eq:conditions2}
		\begin{cases}
			(M(\b{c})\nabla\b{w}) \cdot \b{n} = \b{0} & \quad \text{on }\partial\Omega \times (0,T), \\
			\nabla c_i \cdot \b{n} = 0 & \quad \text{on }\partial\Omega \times (0,T),\quad  \forall \: i \in \mathcal{I},\\
			\Phi = 0 & \quad \text{on }\partial\Omega \times (0,T), \\
			\b{c}(\cdot, 0) = \b{c}_0 & \quad \text{in }\Omega.
		\end{cases}
	\end{equation}
	The first results deals with the existence of a weak solution to \eqref{eq:strongform2}-\eqref{eq:conditions2}, which we define hereafter.
	\begin{defin} \label{def:weaksol}
		A triplet $(\b{c},\, \b{w},\, \Phi)$ satisfying 
		\begin{enumerate}[(i)]
			\item $\b{c} \in L^2(0,T;V) \cap H^1(0,T;V^*)$;
			\item $\b{c} \in \b{G}$ almost everywhere in $\Omega \times (0,T)$;
			\item $\b{w} \in L^2(0,T;V)$;
			\item $\Phi \in L^2(0,T;V_0)$;
			\item the following equalities hold
			\begin{equation*}
				\begin{cases}
					\left\langle \pd{\b{c}}{t}, \: \b{v} \right\rangle_{V^*, V} + \left( M(\b{c})\nabla \b{w}, \nabla \b{v} \right)_H = 0, & \quad \forall \: \b{v} \in V, \\[0.5cm]
					\left(\b{w}, \b{z} \right)_H = \left( \Gamma \nabla \b{c} , \nabla \b{z} \right)_H + \left( \b P\widetilde{\b{\mu}} , \b{z} \right)_H, & \quad \forall \: \b{z} \in V \cap L^\infty(\Omega;\mathbb{R}^3), \\[0.5cm]
					(\epsilon({c}_A,\,{c}_B)(\b{E}_0 - \nabla \Phi), \nabla \Lambda)_H = 0 & \quad \forall \: \Lambda \in V_0,
				\end{cases}
			\end{equation*}
			where
			\[
			\widetilde{\b{\mu}} := \begin{bmatrix}
				\theta_A\Psi'(c_A) + \pd{I}{c_A} + \alpha_{AA}\mathcal{N}(c_A - \overline{c_A}) + \alpha_{AB}\mathcal{N}(c_B - \overline{c_B}) + \dfrac{1}{2}\pd{\epsilon}{c_A}|\b{E}_0 - \nabla \Phi|^2 \\[0.5cm]
				\theta_B\Psi'(c_B) + \pd{I}{c_B} + \alpha_{BA}\mathcal{N}(c_A - \overline{c_A}) + \alpha_{BB}\mathcal{N}(c_B - \overline{c_B}) + \dfrac{1}{2}\pd{\epsilon}{c_B}|\b{E}_0 - \nabla \Phi|^2 \\[0.3cm]
				\theta_S\Psi'(c_S) + \pd{I}{c_S}
			\end{bmatrix};
			\]
			\item the initial condition $\b{c}(0) = \b{c}_0$ is attained in a strong sense;
		\end{enumerate}
		is called a weak solution to \eqref{eq:strongform}-\eqref{eq:conditions}.
	\end{defin}
	\begin{remark}
		As a straightforward consequence of Definition \ref{def:weaksol}, we have $\b{c} \in C^0([0,T];H)$, hence the initial condition is indeed attained in the strong sense.
	\end{remark}
	\begin{thm} \label{thm:wellposed} Let $d \in\{2,3\}$ and assume Assumptions \ref{hyp:potential}-\ref{hyp:data}. Then, there exists a weak solution $(\b{c},\, \b{w},\, \Phi)$ to \eqref{eq:strongform2}-\eqref{eq:conditions2}. Moreover, there exists $p = p(\Omega,\,\epsilon_*,\,\epsilon^*)> 2$ such that the solution satisfies
		\begin{enumerate}[(i)]
			\item $\b{c} \in C^{0,\frac{1}{4}}([0,T]; H) \cap L^\infty(0,T;V) \cap L^\frac q2(0,T;W^{2,\frac q2}(\Omega;\mathbb{R}^3)) \cap H^1(0,T;V^*)$;
			\item $\b{w} \in L^2(0,T;V)$;
			\item $\Phi \in L^\infty(0,T;W^{1,p}(\Omega))$,
		\end{enumerate}
	where $q := \min\{p,4\}$.
	\end{thm}
	
	\begin{remark}
		In two dimensions, Theorem \ref{thm:wellposed} gives H\"{o}lder regularity for $\Phi$ (with respect to the spatial variable). In particular, we have
		\[
		\Phi \in L^\infty(0,T;C^{0,1-\frac{2}{p}}(\Omega))
		\]
		by Sobolev embeddings. If $p > 3$, then the same result holds in three dimensions as well, and in that case we have
		\[
		\Phi \in L^\infty(0,T;C^{0,1-\frac{3}{p}}(\Omega)).
		\]
	\end{remark} \noindent
	Let us further state an immediate corollary which can be proved with minor modifications.
	\begin{cor}
		\label{cor:wellposed} Let $d \in\{2,3\}$ and assume Assumptions \ref{hyp:potential}-\ref{hyp:data}. Then, there exists a weak solution $(\b{c},\, \b{w},\, \Phi)$ to \eqref{eq:strongform2}-\eqref{eq:conditions2} on the time \reviewF{interval} $[0,+\infty)$. Moreover, there exists $p = p(\Omega,\,\epsilon_*,\,\epsilon^*)> 2$ such that the solution satisfies
		\begin{enumerate}[(i)]
			\item $\b{c} \in C^{0,\frac{1}{4}}([0,+\infty); H) \cap L^\infty(0,+\infty;V) \cap L^\frac q2_{\text{loc}}(\reviewF{[0,+\infty)};W^{2,\frac q2}(\Omega;\mathbb{R}^3)) \cap H^1_{\text{loc}}(\reviewF{[0,+\infty)};V^*)$;
			\item $\b{w} \in L^2_{\text{loc}}(\reviewF{[0,+\infty)};V)$;
			\item $\Phi \in L^\infty(0,+\infty;W^{1,p}(\Omega))$,
		\end{enumerate}
		where $q := \min\{p,4\}$.
	\end{cor} \noindent
	We are able to show additional regularity and an energy identity under slightly stronger assumptions.
	\begin{enumerate}[resume, label = \textbf{(A\arabic*)}]
		\item \label{hyp:EIadd} Assume that 
		\begin{enumerate}[(i)]
			\item \label{hyp:EIq} the external electric field satisfies $\b E_0 \in L^r(\Omega;\mathbb{R}^d)$ for some $r \geq 4$ and $\Phi \in L^\infty(0,T;W^{1,4}(\Omega))$;
		\end{enumerate}
		and one of the following:
		\begin{enumerate}[label=(ii)$_\text{\alph*}$]
			\item \label{hyp:EIE0_a} the external electric field $\b{E}_0 \in W^{1,\alpha}(\Omega;\mathbb{R}^d)$, with $\alpha > \frac 65$;
			\item \label{hyp:EIE0_b} the external electric field $\b{E}_0 \in W^{1,\alpha}(\Omega;\mathbb{R}^d)$, with $\alpha \geq 4$.
		\end{enumerate}
	\end{enumerate}
	
	\begin{thm} \label{thm:energyid}
		Let $d \in\{2,3\}$ and assume Assumptions \ref{hyp:potential}-\ref{hyp:data}. If Assumption \ref{hyp:EIadd}-\ref{hyp:EIq} and \ref{hyp:EIadd}-\ref{hyp:EIE0_a}  hold, then $\b{c} \in L^4(0,T;H^2(\Omega;\mathbb{R}^3))$ and $\Phi \in L^4(0,T;W^{2,\beta}(\Omega))$ for some $\beta = \beta(\alpha,p)> 1$. Furthermore, if also Assumption \ref{hyp:EIadd}-\ref{hyp:EIE0_b} holds, then $\b{c} \in C^0([0,T];V)$ and the energy identity
		\[
		\mathcal{E}(\b{c}(t), \Phi(t)) + \int_{0}^{t} \int_\Omega M(\b{c}(s))\nabla \b{w}(s) \cdot \nabla \b{w}(s) \: \d x \: \d s =	 \mathcal{E}(\b{c}_0, \Phi_0),
		\]
		holds for all $t$ in $[0,T]$. 
	\end{thm} \noindent
	\begin{remark}
		The additional assumption \ref{hyp:EIadd}-\ref{hyp:EIq} is satisfied, for instance, if the exponent $p$ appearing in Theorem \ref{thm:wellposed} is at least four. In general, one would expect $p$ to match the optimal regularity given by the integrability order of $\b E_0$ (i.e., that $p = r$, with $r$ the one appearing in Assumption \ref{hyp:data}-\ref{hyp:E0}). Unfortunately, one is not able to retrieve, in general, optimal integrability properties for $\Phi$, even if $r = +\infty$ and $\partial \Omega$ is smooth. Indeed, the exponent $p$ given by Theorem \ref{thm:wellposed} can turn out to be arbitrarily close to $2$. Actually, it holds $p \to 2^+$ as $\epsilon_* \to 0^+$ or $\epsilon^* \to +\infty$ and $p \to +\infty$ as the ratio between $\epsilon_*$ and $\epsilon^*$ approaches 1 (see \cite[Section 5]{Meyers63}). On account of this observation, if the oscillation of $\epsilon$ is sufficiently close to 0, we may assume $p$ to be arbitrarily large.
	\end{remark} \noindent
	Next, we show an a posteriori uniqueness result for problem \eqref{eq:strongform2}-\eqref{eq:conditions2}. In order to prove it, we need a set of different additional assumptions, mainly requiring some higher differentiablity order for the coefficients:
	\begin{enumerate}[resume, label = \textbf{(A\arabic*)}]
		\item \label{hyp:additional} Assume that
		\begin{enumerate}[(i)]
			\item the external electric field satisfies $\b E_0 \in L^r(\Omega;\mathbb{R}^d)$ for some $r \geq 4$;
			\item the mobility coefficients satisfy $M_{ij} \in C^1(\mathbb{R}^3)$ for all $(i,j) \in \mathcal{I} \times \mathcal{I}$;
			\item the electric permittivity fulfills $\epsilon \in C^2(\mathbb{R}^2)$.
		\end{enumerate}
	\end{enumerate}
	The uniqueness result follows as a corollary to the next theorem.
	\begin{thm}\label{thm:uniqueness}
		Let $d \in \{2,3\}$ and assume Assumptions \ref{hyp:potential}-\ref{hyp:data} and \ref{hyp:additional}. Let $\b{c}_{01}$ and $\b{c}_{02}$ be functions complying with Assumption \ref{hyp:data}-\ref{hyp:initial}.  Moreover, for $k \in \{1,2\}$, denote by $(\b{c}_k, \b{w}_k, \Phi_k)$ a weak solution given by Theorem \ref{thm:wellposed} originating from the initial state $\b{c}_{0k}$. If the additional regularity assumptions
		\[
		\pd{\b{c}_k}{t} \in L^\frac{8}{8-d}(0,T;H), \quad \b{c}_k \in L^{2^d}(0,T;H^2(\Omega;\mathbb{R}^3)), \quad \b{w}_2 \in L^\frac{8}{6-d}(0,T;V), \quad \Phi_k \in L^\infty(0,T;W^{1,4}(\Omega))
		\]
		hold for $k \in \{1,2\}$, then the following continuous dependence estimate
		\[
		\|\b{c}_1 - \b{c}_2\|_{L^\infty(0,T;V^*)} + \|\b{c}_1 - \b{c}_2\|_{L^2(0,T;V)}+ \|\Phi_1 - \Phi_2 \|_{L^2(0,T;V_0)} \leq C\left[ \|\b{c}_{01} - \b{c}_{02}\|_{V^*} + \|\overline{\b{c}_{01}} - \overline{\b{c}_{02}}\|^\frac 12_{\mathbb{R}^3}\right],
		\]
		holds, where the constant $C$ depends on the structural parameters of the problem, including $T$. Furthermore, if the mobility coefficients are constant, and the resulting mobility matrix $M \in \mathcal{M}$, then the continuous dependence estimate holds without any additional regularity assumption on the solutions.
	\end{thm}
	\begin{remark} 
		The extra requirements stated in Theorem \ref{thm:uniqueness} on the order parameter vector $\b c$ and the chemical potential $\b w$ are the same stated in \cite{Barrett99}, and reduce to the requirement on the time derivatives and on the electric potentials when $d = 2$.
	\end{remark} \noindent
	Next, we tackle the regularization of weak solutions.
	The following result shows that, under the previously stated regularity conditions on the coefficients, any weak solution regularizes instantaneously.
	\begin{thm} \label{thm:regularization} 	Let $d \in \{2,3\}$ and assume Assumptions \ref{hyp:potential}-\ref{hyp:data} and \ref{hyp:additional}. If the additional regularity assumptions
		\[
		\pd{\b{c}}{t} \in L^\frac{8}{8-d}(0,T;H), \quad \b{c} \in L^{2^d}(0,T;H^2(\Omega;\mathbb{R}^3)), \quad \b{w} \in L^\frac{8}{6-d}(0,T;V), \quad \Phi \in L^\infty(0,T;W^{1,4}(\Omega))
		\]hold, then for any $\sigma > 0$, there exists a constant $C > 0$, possibly depending on $\sigma$ and the structural parameters of the problem, such that
		\begin{align*}
			\|\partial_t\b{c}\|_{L^\infty(\sigma,t;V^*)} + \|\partial_t\b{c}\|_{L^2(t,t+1;V)} + \|\b{w}\|_{L^\infty(\sigma,t;V)} & \leq C
		\end{align*}
		for any $t \geq \sigma$. Additionally, under Assumption \ref{hyp:EIadd}-(ii)$_{a}$ we further have
		\begin{align*}
				\|\Phi\|_{L^\infty(\sigma,t;W^{2,\beta}(\Omega))}  + \|\b{c}\|_{L^\infty(\sigma,t;H^{2}(\Omega;\mathbb{R}^3))} + \sum_{i \in \mathcal{I}} \|\Psi'(c_i)\|_{L^\infty(\sigma,t;H)} & \leq C
		\end{align*}		
		for any $t \geq \sigma$, where $\beta > 1$ is the one appearing in Theorem \ref{thm:energyid}.	Furthermore, if the mobility coefficients are constant, and the resulting mobility matrix $M \in \mathcal{M}$, then the results hold without any additional regularity assumption on the solutions.
	\end{thm} \noindent
	Finally, we state a result concerning the asymptotic behavior of the weak solution of the system. \reviewF{Let us recall that choosing the initial condition $\b c_0 \in X^{\b m}$ uniquely defines the initial state $\Phi_0$ of the electric potential (cf. Assumption \ref{hyp:data}-\ref{hyp:E0}).} Let then $\b{m} \in \b G$ be arbitrary but fixed and let us define the phase space
	\[
	\mathcal{X}:= \left\{\b{c}_0 \in X^{\b{m}} : \reviewF{\mathcal{E}_1(\b{c}_0) < +\infty} \right\},
	\] 
	endowed with the $H^1$-distance, that we denote by $\operatorname{dist}_{\mathcal{X}}$. 
	Observe that it is immediate to prove that $(\mathcal{X},\, \operatorname{dist}_{\mathcal{X}})$ is a complete metric space.
	\reviewF{	
	\begin{remark}
		Owing to all the holding assumptions, it is straightforward to observe that $\mathcal{X}$ precisely contains all the initial states carrying finite energy.
	\end{remark}}
	\begin{remark}
		Let us point out that, given the trajectory of the order parameter $\b{c}$, the trajectory (in a suitable space) of the electric potential $\Phi$ is uniquely determined, as $\Phi(t) = \mathscr{S}(c_A(t), c_B(t))$ for all $t \geq 0$, where $\mathscr{S}$ denotes the solution map of the elliptic problem
		\[
		\begin{cases}
			\div(\epsilon(c_A(t),c_B(t))(\b E_0 - \nabla \Phi(t))) = 0 & \qquad \text{in } \Omega,\\
			\Phi = 0  & \qquad \text{on } \partial\Omega,
		\end{cases}
		\]
		analyzed in Appendix \ref{app:A}.
	\end{remark} \noindent
	From here onwards, we assume the solution $\b c$ to be uniquely defined. Consider the one-parameter family of solution maps $S$ as follows
	\[
	S(t): \mathcal{X} \to \mathcal{X}, \qquad S(t)(\b{c}_0) := \b{c}(t),
	\]
	where $\b{c}$ is the unique solution to \eqref{eq:strongform2}-\eqref{eq:conditions2} whose initial state equals $\b{c}_0$.
	\begin{remark}
		Since any weak solution exhibits instantaneous regularization (cf. Theorem \ref{thm:regularization}) and it is possible to investigate the longtime behavior of a solution by shifting it forwards in time (see also, for instance, \cite[Section 7]{Abels}), there is no loss of generality in considering a global higher degree of regularity for the solutions. In particular, such a solution satisfies, owing to all the proven theorems,
		\begin{align*}
			\b{c} &\in C^{0,\frac{1}{4}}([0,+\infty); H) \cap L^ \infty(0,+\infty;H^2(\Omega;\mathbb{R}^3)) \cap W^{1,\infty}(0,+\infty;V^*),\\
			\b{w} &\in L^\infty(0,+\infty;V),\\
			\Phi &\in L^\infty(0,+\infty;W^{2,\beta}(\Omega)),
		\end{align*}
		with $\beta$ defined in Theorem \ref{thm:energyid}.
	\end{remark} \noindent 
	Let us introduce the stationary version of the system \eqref{eq:strongform2}-\eqref{eq:conditions2}, namely
	\begin{equation} \label{eq:stationary} \footnotesize
		\begin{cases}
			\b{w}^\infty = \b{P}\begin{bmatrix}
				-\gamma_A\Delta c_A^\infty + \theta_A\Psi'(c_A^\infty) + \pd{I}{c_A^\infty} + \alpha_{AA}\mathcal{N}(c_A^\infty - \overline{c_A^\infty}) + \alpha_{AB}\mathcal{N}(c_B^\infty - \overline{c_B^\infty}) + \dfrac{1}{2}\pd{\epsilon}{c_A^\infty}|\b{E}_0 - \nabla \Phi^\infty|^2 \\[0.3cm] 
				-\gamma_B\Delta c_B^\infty + \theta_B\Psi'(c_B^\infty) + \pd{I}{c_B^\infty} + \alpha_{BA}\mathcal{N}(c_A^\infty - \overline{c_A^\infty}) + \alpha_{BB}\mathcal{N}(c_B^\infty - \overline{c_B^\infty}) + \dfrac{1}{2}\pd{\epsilon}{c_B^\infty}|\b{E}_0 - \nabla \Phi^\infty|^2  \\[0.3cm]
				-\gamma_S\Delta c_S^\infty + \theta_S\Psi'(c_S^\infty) + \pd{I}{c_S^\infty}
			\end{bmatrix} \\
			\displaystyle \sum_{i \in \mathcal{I}} c_i^\infty = 1, \\
			\div\left[ \varepsilon(c_A^\infty,\, c_B^\infty)(\b{E}_0 - \nabla \Phi^\infty)\right] = 0, 
		\end{cases}
	\end{equation}
	with $\b{w}^\infty := [w_A^\infty,\,w_B^\infty,\,w_S^\infty]^T \in \mathbb{R}^3$ being a constant vector. Problem \eqref{eq:stationary} is subject to the boundary conditions
	\begin{equation} \label{eq:stationaryBCs}
		\begin{cases}
			\nabla c_i^\infty \cdot \b{n} = 0 & \quad \text{on }\partial \Omega, \quad \forall \: i \in \mathcal{I}, \\
			\Phi^\infty = 0 & \quad \text{on } \partial\Omega.
		\end{cases} 
	\end{equation}
	For the next result we need a stronger version of the previous Assumptions, namely
	\begin{enumerate}[resume, label = \textbf{(A\arabic*)}]
		\item \label{hyp:additionalrefined} Assume that
		\begin{enumerate}[(i)]
			\item \review{the mobility matrix $M( \b s) \equiv M \in \mathcal{M}$};
			\item \review{letting $\alpha > 0$ be such that $\b E_0 \in W^{1,\alpha}(\Omega;\mathbb{R}^d)$ and $p > 1$ be such that the electric potential fulfills $\Phi \in L^\infty(\tau,+\infty;W^{1,p}(\Omega))$ for some $\tau \geq 0$, it holds $\min\{\alpha, p\} \geq \frac{6d}{6-d}$};
			\item the potential energy density function $\Psi$ is analytic in $(0,1)$,
			\item the local interaction function $I$ is analytic in $(0,1)^3$;
			\item the electric permittivity $\epsilon$ is \review{affine-linear in an open neighbourhood $\mathcal{U}_1$ of the set
			\[
			\mathcal{U}_0 := \{ (x,y) \in \mathbb{R}^2: x \geq 0,\,y \geq 0,\, x+y\leq 1\}.
			\]}.
		\end{enumerate} 
	\end{enumerate}
	\review{\begin{remark}
			The additional assumption on the electric permittivity $\epsilon$ enables to avoid cumbersome technicalities linked to analyticity issues. Nonetheless, it should be possible to only ask for local analyticity around $\mathcal{U}_0$. 
	\end{remark}}
	\begin{thm} \label{thm:longtime}
		Let $d \in \{2,3\}$ and assume Assumptions \ref{hyp:potential}-\ref{hyp:data} and \ref{hyp:additionalrefined}. Let $\b c_0 \in V$ comply with the hypotheses of Theorem \ref{thm:wellposed}. Consider the trajectory $(\b c(t), \Phi(t))$ originating from $\b c_0$ and the corresponding initial electric potential $\Phi_0$. Then, there exist a pair $(\b c^\infty,\, \Phi^\infty)$ being a strong solution of \eqref{eq:stationary}-\eqref{eq:stationaryBCs} such that
		\[
		\lim_{t \to +\infty} \b c(t) = \b c^\infty \qquad \text{ in } H^{2\vartheta}(\Omega;\mathbb{R}^3),
		\]
		for all $\vartheta \in (0,1)$ and
		\[
		\lim_{t \to +\infty} \Phi(t) = \Phi^\infty \qquad \text{ in } V_0.
		\]
	\end{thm}
	
\section{Proofs of Theorem \ref{thm:wellposed} and Theorem \ref{thm:energyid}} \label{sec:existence}
The present section is devoted to proving the existence of weak solutions stated by Theorem \ref{thm:wellposed} and the additional properties listed in Theorem \ref{thm:energyid}. For the sake of clarity, we split the proof in several steps.
\subsection{Regularization of the singular potential.} First of all, we exhibit a family of regular functions approximating the singular potential $\Psi$. To this end, we employ a fourth-order Taylor approximation of $\Psi$. Let $\delta \in (0,1)$ be given, and define 
\begin{equation} \label{eq:regularization}
	\Psi_\delta : \mathbb{R} \to \mathbb{R}, \qquad \Psi_\delta(s) := 
	\begin{cases}	
		\Psi(s) & \quad \text{if } s \geq \delta, \\
		\displaystyle\sum_{i=0}^4 \dfrac{\Psi^{(i)}(\delta)}{i!}(s-\delta)^i & \quad \text{if }  s < \delta. \\
	\end{cases}
\end{equation}
\begin{remark} \label{rem:approx}
	Observe that, given Assumption \ref{hyp:potential}, the function $\Psi_\delta$ belongs to the space $C^4(\mathbb{R})$ and is non-negative and convex for $\delta$ sufficiently small.
\end{remark} \noindent
Moreover, in the following, we denote by $\widetilde{I}$ a $C^2(\mathbb{R}^3)$-extension of $I$ so that there exists a positive constant $C_{\widetilde{I}}'$ such that
\begin{equation} \label{eq:boundIder}
	\left|\pd{\widetilde{I}}{s_i}(\b{s})\right| \leq C_{\widetilde{I}}' (1 + |\b{s}|) \qquad \forall \: \b{s} \in \mathbb{R}^3, \quad \forall \: i \in \{1,2,3\},
\end{equation}
which in turn implies the existence of some constant $C_{\widetilde{I}} > 0$ such that
\begin{equation} \label{eq:boundI}
	|\widetilde{I}(\b{s})| \leq C_{\widetilde{I}} (1 + |\b{s}|^2) \qquad \forall \: \b{s} \in \mathbb{R}^3, \quad \forall \: i \in \{1,2,3\}.
\end{equation}
Introducing the regularized free energy density

\[
F_\delta: \mathbb{R}^3 \to \mathbb{R}, \qquad F_\delta(\b{s}) := \sum_{i=1}^3 \theta_{\ell(i)}\Psi_\delta(s_i) + \widetilde{I}(\b{s}) +\delta |\b s|^4, \qquad \forall \: \b{s} = [s_1, s_2, s_3]^T \in \mathbb{R}^3,
\]we are in a position to state a preliminary result stating a uniform bound from below for $F_\delta$.
\begin{lem} \label{lem:bound}
	There exists some $\delta_0 > 0$ and a constant $C_1 = C_1(\delta)> 0$ such that
	\[
	F_\delta(\b{s}) \geq \frac \delta 2|\b{s}|^4 - C_1 \qquad \forall \: \b{s} \in \mathbb{R}^3,
	\]
	for all $\delta \in (0,\delta_0)$.
\end{lem}
\begin{proof}
	First, we recall that $\theta_*$ must be strictly positive, and that by \eqref{eq:boundI}
	\[
	\widetilde{I}(\b{s}) \geq - C_{\widetilde{I}}(1 + |\b{s}|^2)
	\]
	for all $\b{s} \in \mathbb{R}^3$. Moreover, if $\delta \leq \min\{\varepsilon_0,\,\varepsilon_1\}$, as a consequence of Assumption \ref{hyp:potential} and the definition of $\Psi_\delta$ we have that
	\[
	\Psi_\delta(r) \geq \dfrac{c_\Psi}{48}r^4 - C_0 \qquad \forall \: r < \delta,
	\]
	for some positive constant $C_0$ independent of $r$, and, still by Assumption \ref{hyp:potential},
	\[
	\Psi_\delta(r) \geq 0 \qquad \forall \: r \geq \delta.
	\]
	Then, the two properties above imply that
	\[
	\begin{split}
		F_\delta(\b{s}) & \geq \delta |\b s|^4-\sum_{i=1}^3 \theta_{\ell(i)}C_0 - C_{\widetilde{I}}(1 + |\b{s}|^2) \\
		& \geq \delta |\b s|^4 - 3\theta^*C_0 - C_{\widetilde{I}}(1 + |\b{s}|^2) \\
		& \geq \dfrac{\delta}{2}|\b s|^4- 3\theta^*C_0- C_{\widetilde{I}}-\dfrac{C_{\widetilde{I}}^2}{2\delta}
	\end{split}
	\]
	for all $\b s \in \mathbb{R}^3$. The claim follows by setting the constant $C_1 := 3\theta^*C_0+ C_{\widetilde{I}}+2\delta^{-1}C_{\widetilde{I}}^2$, while $\delta_0 := \min\{\varepsilon_0, \: \varepsilon_1\}$.
\end{proof}
\subsection{The semi-implicit time discretization.} Let $T > 0$ be a given final time. For an arbitrary $N \in \mathbb{N}$, define the time step
\[
\tau := \dfrac{T}{N}, 
\]
and, for every $n \in \{1, \dots, N\}$ and $\delta \in (0,\delta_0)$, consider the following regularized system:
\begin{equation} \label{eq:discform} \small
	\begin{cases}
		\dfrac{\b{c}^{\delta, n} - \b{c}^{\delta, n-1}}{\tau} = \div(M(\b{c}^{\delta, n- 1 }) \nabla \b{w}^{\delta, n}),\\[0.2cm]
		\b w^{\delta, n} = \b P \begin{bmatrix}
			-\gamma_A\Delta c_A^{\delta, n} + \pd{F_\delta}{s_1}(c_A^{\delta, n}) + \alpha_{AA}\mathcal{N}(c_A^{\delta, n} - \overline{c_A^{\delta, n}}) + \alpha_{AB}\mathcal{N}(c_B^{\delta, n} - \overline{c_B^{\delta, n}}) + \dfrac{1}{2}\pd{\epsilon}{s_1}(c_A^{\delta, n})|\b{E}_0 - \nabla \Phi^{\delta, n}|^2 \\[0.4cm]
			-\gamma_B\Delta c_B^{\delta, n} + \pd{F_\delta}{s_2}(c_B^{\delta, n}) + \alpha_{BA}\mathcal{N}(c_A^{\delta, n} - \overline{c_A^{\delta, n}}) + \alpha_{BB}\mathcal{N}(c_B^{\delta, n} - \overline{c_B^{\delta, n}}) + \dfrac{1}{2}\pd{\epsilon}{s_2}(c_B^{\delta, n})|\b{E}_0 - \nabla \Phi^{\delta, n}|^2 \\[0.4cm]
			-\gamma_S\Delta c_S^{\delta, n} + \pd{F_\delta}{s_3}(c_S^{\delta, n})
		\end{bmatrix} , \\[0.3cm] 
		\div\left[ \varepsilon(c_A^{\delta, n},\, c_B^{\delta, n})(\b{E}_0 - \nabla \Phi^{\delta, n})\right] = 0,
	\end{cases}
\end{equation}
subject to the boundary conditions
\begin{equation} \label{eq:disccond}
	\begin{cases}
		M(\b{c}^{\delta, n-1})\nabla\b w^{\delta, n} \cdot \b{n} = \b{0} & \quad \text{on }\partial\Omega \times (0,T), \\
		\nabla c_i^{\delta, n} \cdot \b{n} = 0 & \quad \text{on }\partial\Omega \times (0,T),\quad  \forall \: i \in \mathcal{I},\\
		\Phi^{\delta, n} = 0 & \quad \text{on }\partial\Omega \times (0,T).
	\end{cases}
\end{equation}
Here, the two indexes $\delta$ and $n$ refer to the value of $\delta$ in the potential regularization \eqref{eq:regularization} and to the evaluation of the approximated solution at the time instant $n\tau$, respectively. In order to find a solution to \eqref{eq:discform}-\eqref{eq:disccond}, we observe that the system can be obtained as a discrete version of the gradient flow of the regularized energy functional
\[
\mathcal{E}^\delta(\b{c},\, \Phi) := \mathcal{E}_1^\delta(\b{c}) + \mathcal{E}_2(\b{c}) + \mathcal{E}_3(c_A,\, c_B) + \mathcal{E}_4(c_A,\,c_B,\, \Phi),
\] 
where 
\[
\mathcal{E}_1^\delta(\b{c}) := \ii{\Omega}{}{\sum_{i \in \mathcal{I}} \left[  \dfrac{\gamma_i}{2}|\nabla c_i|^2 + \theta_i\Psi_\delta(c_i) \right]}{x}.
\]\review{
Indeed, defining the functional
\begin{equation} \label{eq:regenergy}
	\mathcal{E}^{\delta,n}(\b{v}, \Lambda) = \mathcal{E}^\delta(\b{v}, \Lambda) + \dfrac{1}{2\tau}\|\b{v} - \b{c}^{\delta,n-1}\|_{*, M(\b{c}^{\delta,n- 1})}^2,
\end{equation}
we obtain that a minimizer of $\mathcal{E}^{\delta,n}$ fulfills
\begin{equation} \label{eq:eullag1}
	\left\langle \dfrac{\delta \mathcal{E}^\delta}{\delta \b{c}},\, \b{v} \right\rangle_{V^*, V} = - \left( \b{v},\, \dfrac{\b{c}^{\delta, n} - \b{c}^{\delta, n-1}}{\tau} \right)_{*, M(\b{c}^{\delta,n- 1})},
\end{equation}
for all $\b{v} \in TX^{\b{m}}$. Moreover, we also have
\begin{equation} \label{eq:eullag2}
	\left\langle \dfrac{\delta \mathcal{E}^\delta}{\delta \Phi},\, \Lambda \right\rangle_{V^*_0, V_0} = 0,
\end{equation}
for all $\Lambda \in V_0$, which yields \eqref{eq:discform}. Observe that \eqref{eq:eullag1} and \eqref{eq:eullag2} are the Euler--Lagrange equations of the functional $\mathcal{E}^{\delta,n}$. Assuming $\b{c}^{\delta, n-1}$ to be known (starting from the vector of initial states for $n = 1$), we solve the system inductively to find $\b{c}^{\delta,n}$ and $\Phi^{\delta, n}$. Since the system conserves mass, if $\b{c}_0 \in X^{\b{m}}$ for some given $\b{m}$, then $\b{c}^{\delta, n}$ should also belong to the same space, regardless of the values of $\delta$ and $n$. The following is a first existence result.}
\begin{lem}
	Let $\delta \in (0,\delta_0)$, $N \in \mathbb{N}$ and $n \in \{1, \dots, N\}$. For any $\tau > 0$, there exists a minimizer of $\mathcal{E}^{\delta, n}$ in $X^{\b{m}} \times V_0$, where $\b{m} = \overline{\b{c}_0}$.
\end{lem}
\begin{proof}
	Here, we use the direct method of calculus of variations. Let $(\b{v}, \Lambda) \in X^{\b{m}}$ be arbitrary but fixed. First of all, on account of Assumptions \ref{hyp:potential}-\ref{hyp:initial} and Lemma \ref{lem:bound}, we have the following inequality
	\begin{multline*}
		\mathcal{E}^{\delta, n}(\b{v}, \Lambda) \geq  \gamma_*\|\nabla \b{v}\|^2_H + \frac \delta 2\|\b{v}\|^4_{L^4(\Omega; \mathbb{R}^3)} - \sum_{(i,j) \in \{1,2\}^2} |\alpha_{\ell(i)\ell(j)}|\|\nabla \mathcal{N}(v_i - \overline{v_i})\|_H \|\nabla \mathcal{N}(v_j - \overline{v_j})\|_H  \\
		+ \epsilon_*\|\b{E}_0 - \nabla \Lambda\|^2_H + \dfrac{1}{2\tau}\|\b{v} - \b{c}^{\delta,n-1}\|_{*, M(\b{c}^{\delta, n-1})}^2 - C_1,
	\end{multline*}
	where we also used the Cauchy--Schwarz inequality. Using the Young inequality, owing also to the embedding $H_0 \hookrightarrow V^*_0$, we have, for any fixed $(i,j) \in \{1,2\}^2$
	\[
	\begin{split}
		|\alpha_{\ell(i)\ell(j)}|\|\nabla \mathcal{N}(v_i - \overline{v_i})\|_H \|\nabla \mathcal{N}(v_j - \overline{v_j})\|_H & = |\alpha_{\ell(i)\ell(j)}|\|v_i - \overline{v_i}\|_{V_0^*}\|v_j - \overline{v_j}\|_{V_0^*} \\
		& \leq C\|v_i\|_{L^4(\Omega)}\|v_j\|_{L^4(\Omega)} \\
		& \leq C\|\b v\|_{L^4(\Omega; \mathbb{R}^3)}^2 \\
		& \leq \dfrac{\delta}{4}\|\b v\|_{L^4(\Omega; \mathbb{R}^3)}^4 + C,
	\end{split}
	\]
	with the constant $C > 0$ depending on $\delta$ and the parameters $\alpha_{ij}$. Therefore, we arrive at
	\begin{equation*}
		\mathcal{E}^{\delta, n}(\b{v}, \Lambda) \geq \gamma_*\|\nabla \b{v}\|^2_H + \dfrac{\delta}{4}\|\b{v}\|^4_{L^4(\Omega)}  
		+ \epsilon_*\|\b{E}_0 - \nabla \Lambda\|^2_H + \dfrac{1}{2\tau}\|\b{v} - \b{c}^{\delta,n-1}\|_{*, M(\b{c}^{\delta, n-1})}^2 - C_2,
	\end{equation*}
	where the constant $C_2$ depends on $\delta$, on the parameters $\alpha_{ij}$ and on $\overline{\b{c}_0}$. This proves that the regularized energy functional is bounded from below on $X^{\b{m}} \times V_0$ for any $\tau > 0$. Therefore, we can consider a minimizing sequence $\{(\b{v}_k, \Lambda_k)\}_{k \in \mathbb{N}} \subset X^{\b{m}} \times V_0$ such that
	\[
	\lim_{k \to +\infty} \mathcal{E}^{\delta, n}(\b{v}_k, \Lambda_k) = \inf_{(\b{w}, \Xi) \in X^{\b{m}} \times V_0}\mathcal{E}^{\delta, n}(\b{z}, \Xi).	
	\]
	Moreover, since 
	\[
	\epsilon_*\|\b{E}_0 - \nabla \Lambda\|^2_H \geq \dfrac{\epsilon_*}{2}\|\nabla \Lambda\|^2_H - C_3, 
	\]
	where $C_3$ depends on $\epsilon_*$ and the $H$-norm of $\b{E}_0$, we also conclude that the regularized energy functional is actually coercive on $X^{\b{m}} \times V_0$, and therefore the sequence $\{(\b{v}_k, \Lambda_k)\}_{k \in \mathbb{N}}$ is uniformly bounded in $X^{\b{m}} \times V_0$. By reflexivity, we can replace it, without relabeling, with some weakly convergent subsequence to some limit $(\b{v}, \Lambda)$. Furthermore, we can also choose the subsequence in such a way that 
	\[
	\b{v}_k \to \b{v} \text{ in $H$ and almost everywhere in } \Omega.\] We are thus only left to prove that $(\b{v}, \Lambda)$ is indeed a minimizer for $\mathcal{E}^{\delta, n}$. To this end, we pass to the limit in the expression for $\mathcal{E}^{\delta, n}(\b{v}_k, \Lambda_k)$. It is straightforward to pass to the limit in the convex terms by lower semicontinuity, therefore we only focus on the remaining ones, i.e., the interaction terms. Here, we achieve our aim through the Lebesgue dominated convergence theorem. For the local interaction term, it is enough to recall the quadratic bound given by \eqref{eq:boundI}. Instead, for the nonlocal interaction terms, observe that $\mathcal{N}$ is a continuous linear operator in $H$, hence 
	\[
	\mathcal{N}(\b{v}_k - \b{m}) \to \mathcal{N}(\b{v} - \b{m}) \text{ in } H \text{ and almost everywhere in } \Omega,\]
	up to a subsequence. We conclude that, up to a subsequence which we do not relabel,
	\[
	\mathcal{E}^{\delta,n}(\b{v}, \Lambda) \leq \liminf_{k \to +\infty} \mathcal{E}^{\delta,n}(\b{v}_k, \Lambda_k) =  \inf_{(\b{z}, \Xi) \in X^{\b{m}} \times V_0}\mathcal{E}^{\delta, n}(\b{z}, \Xi),
	\]
	and the proof is complete.
\end{proof} \noindent
At each step, i.e., for every $n \in \{1, \dots, N\}$, the minimizer of $\mathcal{E}^{\delta,n}$ is the sought solution $(\b{c}^{\delta, n}, \Phi^{\delta, n})$. The next result states the Euler--Lagrange equations satisfied by the minimizer.
\begin{lem} \label{lem:eul-lag}
	The couple $(\b{c}^{\delta, n}, \Phi^{\delta, n})$ satisfies
	\begin{equation*}
		\begin{cases}
			\left(\dfrac{\b{c}^{\delta, n} -\b{c}^{\delta, n-1}}{\tau}, \: \b{v} \right)_H + \left( M(\b{c}^{\delta, n-1})\nabla \b{w}^{\delta, n}, \nabla \b{v} \right)_H = 0, & \quad \forall \: \b{v} \in V, \\[0.5cm]
			\left(\b{w}^{\delta, n}, \b{z} \right)_H =  \left( \Gamma \nabla \b{c}^{\delta, N} , \nabla \b{z} \right)_H + \left( \b P\widetilde{\b{\mu}}^{\delta,n} , \b{z} \right)_H, & \quad \forall \: \b{z} \in V \cap L^\infty(\Omega;\mathbb{R}^3), \\[0.5cm]
			(\epsilon({c}_A^{\delta,n},\,{c}_B^{\delta,n})(\b{E}_0 - \nabla \Phi^{\delta,n}), \nabla \Lambda)_H = 0 & \quad \forall \: \Lambda \in V_0,
		\end{cases}
	\end{equation*}
	where
	\[
	\widetilde{\b{\mu}}^{\delta,n} := \begin{bmatrix}
		\pd{F_\delta}{s_1}(c_A^{\delta, n}) + \alpha_{AA}\mathcal{N}(c_A^{\delta, n} - \overline{c_A^{\delta, n}}) + \alpha_{AB}\mathcal{N}(c_B^{\delta, n} - \overline{c_B^{\delta, n}}) + \dfrac{1}{2}\pd{\epsilon}{s_1}(c_A^{\delta, n})|\b{E}_0 - \nabla \Phi^{\delta, n}|^2 \\[0.5cm]
		\pd{F_\delta}{s_2}(c_B^{\delta, n}) + \alpha_{BA}\mathcal{N}(c_A^{\delta, n} - \overline{c_A^{\delta, n}}) + \alpha_{BB}\mathcal{N}(c_B^{\delta, n} - \overline{c_B^{\delta, n}}) + \dfrac{1}{2}\pd{\epsilon}{s_2}(c_B^{\delta, n})|\b{E}_0 - \nabla \Phi^{\delta, n}|^2 \\[0.3cm]
		\pd{F_\delta}{s_3}(c_S^{\delta, n})
	\end{bmatrix}.
	\]
\end{lem}
\begin{proof}
	First of all, we compute the first variation of the regularized energy functional using only increments in $(TX^{\b{m}} \cap L^\infty(\Omega; \mathbb{R}^3)) \times V_0$, namely for any fixed $(\b{v}, \Lambda) \in (TX^{\b{m}} \cap L^\infty(\Omega; \mathbb{R}^3)) \times V_0$, we compute
	\[
	\lim_{h \to 0} \dfrac{\mathcal{E}^{\delta,n}(\b{c}^{\delta,n} + h\b{v},\,\Phi^{\delta,n} + h\Lambda) - \mathcal{E}^{\delta,n}(\b{c}^{\delta,n},\,\Phi^{\delta,n})}{h}.
	\]
	Firstly, we deal with the terms involving the regularization $\Psi_\delta$ exploiting convexity (recall Remark \ref{rem:approx}). Indeed, we have
	\[
	\sum_{i =1}^3 \theta_i \Psi_\delta(c_i^{\delta,n}) \geq  \sum_{i=1}^3 \theta_i \Psi_\delta(c_i^{\delta,n} + hv_{i	}) - h\theta_iv_{i}\Psi_\delta'(c_i^{\delta,n} + hv_i),
	\]
	and by Assumption \ref{hyp:potential}-\ref{hyp:potder1-2} we easily get 
	\[
	\sum_{i =1}^3 \theta_i \Psi_\delta(c_i^{\delta,n} + hv_{i}) \leq \sum_{i =1}^3 \theta_i \Psi_\delta(c_i^{\delta,n}) + |h|\theta_i\|v_i\|_{L^\infty(\Omega)}C_\Psi(1 + \Psi_\delta(c_i^{\delta,n} + hv_{i}))
	\]
	and if
	\[
	|h| < \dfrac{\theta_*}{\theta^*\|\b{v}\|_{L^\infty(\Omega;\mathbb{R}^3)}C_\Psi},
	\]
	we arrive at
	\[
	\sum_{i =1}^3 \theta_i \Psi_\delta(c_i^{\delta,n} + hv_{i}) \leq \sum_{i \in \mathcal{I}} \theta_i \Psi_\delta(c_i^{\delta,n}) + C|h|,
	\]
	for some fixed given constant $C > 0$, hence also giving
	\[
	\left|\dfrac{1}{h}\sum_{i =1}^3 \theta_i\left(  \Psi_\delta(c_i^{\delta,n} + hv_{i}) - \Psi_\delta(c_i^{\delta,n}) \right) \right| \leq C\left( 1 + \sum_{i =1}^3 \theta_i \Psi_\delta(c_i^{\delta,n}) \right).
	\]
	Therefore, the dominated convergence theorem entails the equality
	\[
	\lim_{h \to 0}	\dfrac{1}{h}\ii{\Omega}{}{\sum_{i =1}^3 \theta_i\left(  \Psi_\delta(c_i^{\delta,n} + hv_{i}) - \Psi_\delta(c_i^{\delta,n}) \right)}{x} = \ii{\Omega}{}{\sum_{i =1}^3 \theta_i  \Psi_\delta'(c_i^{\delta,n})v_{i}}{x}.
	\]
	The quadratic bound on $\widetilde{I}$ given by \eqref{eq:boundI} enables a similar argument for the local interaction term. For the nonlocal terms, a direct computation yields
	\[
	\begin{split}
		& \lim_{h \to 0}\dfrac{1}{h}\bigg[ \sum_{(i,j) \in \{1,2\}^2} \dfrac{\alpha_{ij}}{2}\bigg( \left(\mathcal{N}(c_i^{\delta,n} + hv_{i} - \overline{c_i^{\delta,n} + hv_{i}}),\, c_j^{\delta,n} + hv_{j} - \overline{c_j^{\delta,n} + hv_{j}}\right)_H \\ & \hspace{9cm} - \left(\mathcal{N}(c_i^{\delta,n}- \overline{c_i^{\delta,n}}), c_j^{\delta,n} - \overline{c_j^{\delta,n}}\right)_H\bigg) \bigg] \\
		& \quad = \lim_{h \to 0}\dfrac{1}{h}\bigg[ \sum_{(i,j) \in \{1,2\}^2} \dfrac{\alpha_{ij}}{2}\bigg( h\left(\mathcal{N}(v_i-\overline{v_i}),\, c_j^{\delta,n} - \overline{c_j^{\delta,n}}\right)_H + h\left(\mathcal{N}(c_i^{\delta,n} - \overline{c_i^{\delta,n}}),\,v_{j}  -\overline{v_{j}}\right)_H \\ & \hspace{9cm} +h^2 \left(\mathcal{N}(c_i^{\delta,n}- \overline{c_i^{\delta,n}}), c_j^{\delta,n} - \overline{c_j^{\delta,n}}\right)_H\bigg) \bigg] \\
		& = \sum_{(i,j) \in \{A,B\}^2} \alpha_{ij} \left(\mathcal{N}(c_i^{\delta,n}- \overline{c_i^{\delta,n}}), \, v_{j} - \overline{v_{j}}\right)_H + \lim_{h \to 0} \sum_{(i,j) \in \{1,2\}^2} h \left(\mathcal{N}(c_i^{\delta,n}- \overline{c_i^{\delta,n}}), c_j^{\delta,n} - \overline{c_j^{\delta,n}}\right)_H \\
		& = \sum_{(i,j) \in \{A,B\}^2} \alpha_{ij} \left(\mathcal{N}(c_i^{\delta,n}- \overline{c_i^{\delta,n}}), \, v_{j} - \overline{v_{j}}\right)_H,
	\end{split}
	\]
	where we also exploited the fact that $\mathcal{N}$ is self-adjoint. Next, we address the electrical term. Computing the derivative of the electrical energy contribution leads to 
	\[
	\begin{split}
		& \lim_{h \to 0} \dfrac{1}{h}\left[\mathcal{E}_4(c_A+hv_1,\,c_B+hv_2,\,\Phi + h\Lambda) -\mathcal{E}_4(c_A,\,c_B,\,\Phi) \right] \\
		& \quad = - \left(\epsilon(c_A^{\delta,n}, c_B^{\delta,n})(\b{E}_0 - \nabla \Psi), \nabla \Lambda \right)_H + \sum_{i =1}^2 \dfrac 12 \left(\pd{\epsilon}{s_i}(c_i^{\delta, n})v_{i}, |\b{E}_0 - \nabla \Psi|^2\right)_H.
	\end{split}
	\]
	\review{
	For the perturbation term arising from the regularization scheme, we have
	\[
	\begin{split}
		& \lim_{h \to 0} \dfrac{\delta(c_{\ell(j)}+hv_j)^4-\delta c_{\ell(j)}^4}{h} = 4\delta c_{\ell(j)}^3 v_j
	\end{split}
	\]
	for all $j \in \mathcal{J}$.}
	The additional term arising from the discretization scheme gives the contribution 
	\[
	\begin{split}
		& \lim_{h \to 0} \dfrac{1}{2\tau h}\left[\|\b{c}^{\delta, n} + h \b{v} -\b{c}^{\delta, n-1}\|^2_{*, M(\b{c}^{\delta, n-1})} - \|\b{c}^{\delta, n} -\b{c}^{\delta, n-1}\|^2_{*, M(\b{c}^{\delta, n-1})}  \right] \\
		& \hspace{1cm} = \lim_{h \to 0} \dfrac{1}{\tau h}(\b{c}^{\delta, n} -\b{c}^{\delta, n-1}, h\b{v})_{*, M(\b{c}^{\delta, n-1})} \\
		&  \hspace{1cm} = \dfrac{1}{\tau}(\b{c}^{\delta, n} -\b{c}^{\delta, n-1}, \b{v})_{*, M(\b{c}^{\delta, n-1})} \\
		& \hspace{1cm} = (\b{w}^{\delta, n} - \overline{\b{w}^{\delta, n}}, \b{v})_H.
	\end{split}
	\]
	The conclusion follows from the minimizing property of $(\b{c}^{\delta,n}, \Phi^{\delta,n})$, as it makes the Gateaux derivative of the energy functional vanish. The same result for generic $\b{v}$ readily follows.
\end{proof}
\subsection{Uniform estimates.} \label{ssec:unifest}Here we give uniform estimates with respect to the parameters $N$ and $\delta$. To this end, we first need to define the approximated solutions on the whole interval $[0,T]$ by suitably interpolating the $N$ known values obtained by inductively solving \eqref{eq:discform}-\eqref{eq:disccond}. First, we define the constant interpolant. Let $\b{c}^{\delta,N}$, $\b{w}^{\delta, N}$ and $\Phi^{\delta,N}$ be such that for any $n \in \{1, \dots, N\}$ and $t \in ((n-1)\tau, n\tau]$ we have
\[
\b{c}^{\delta,N}(t) := \b{c}^{\delta,n}(t), \qquad \b{w}^{\delta, N}(t) := \b{w}^{\delta, n}(t), \qquad
\Phi^{\delta,N}(t) := \Phi^{\delta,n}(t).
\]
Of course, such a solution is not continuous in time. In order to recover such a property, we also define the piecewise linear interpolant as follows. For any $\beta \in [0,1]$, $n \in \{1, \dots, N\}$ and $t = \beta n \tau + (1-\beta)(n-1)\tau$, we set
\begin{align*}
	\widetilde{\b{c}}^{\delta,N}(t) & := \beta \b{c}^{\delta,n}(t) + (1-\beta)\b{c}^{\delta,n-1}(t), \\
	\widetilde{\b{w}}^{\delta,N}(t) & := \beta \b{w}^{\delta,n}(t) + (1-\beta)\b{w}^{\delta,n-1}(t), \\
	\widetilde{\Phi}^{\delta,N}(t) & := \beta \Phi^{\delta,n}(t) + (1-\beta)\Phi^{\delta,n-1}(t).
\end{align*}
Moreover, regardless of the kind of interpolation, we extend the regularized order parameters and electric potentials backwards in time in such a way that
\[
\b{c}^{\delta,N}(t) = \widetilde{\b{c}}^{\delta,N}(t) = \b{c}_0, \qquad \Phi^{\delta,N}(t) = \widetilde{\Phi}^{\delta,N}(t) = \Phi_0, \qquad  \forall \: t \in [-\tau,0].
\]
The first two results concern a discrete energy inequality and its basic consequences.
\begin{lem} \label{lem:energy}
	For any $\delta \in (0,\delta_0)$ and $N \in \mathbb{N}$, the constant interpolant satisfies the energy inequality
	\[
	\mathcal{E}^\delta(\b{c}^{\delta, N}(t), \Phi^{\delta,N}(t)) + \int_{0}^{t} \int_\Omega M(\b{c}^{\delta, N}(s-\tau))\nabla \b{w}^{\delta, N}(s) \cdot \nabla \b{w}^{\delta,N}(s) \: \d x \: \d s \leq \mathcal{E}^\delta(\b{c}_0, \Phi_0),
	\]
	for all $t \in [0,T]$.
\end{lem}
\begin{proof}
	Fix $\delta \in (0,\delta_0)$ and $N \in \mathbb{N}$. Since $(\b{c}^{\delta, n}, \Phi^{\delta,n})$ minimizes \eqref{eq:regenergy}, we easily have
	\[
	\mathcal{E}^{\delta,n}(\b{c}^{\delta, n}, \Phi^{\delta,n}) \leq \mathcal{E}^{\delta,n}(\b{c}^{\delta, n-1}, \Phi^{\delta,n-1}) 
	\]
	for all $n \in \{1,\dots,N\}$. This is equivalent to
	\[
	\mathcal{E}^\delta(\b{c}^{\delta, n}, \Phi^{\delta,n}) + \dfrac{1}{2\tau}\|\b{c}^{\delta, n} - \b{c}^{\delta,n-1}\|_{*, M(\b{c}^{\delta,n- 1})}^2 \leq \mathcal{E}^\delta(\b{c}^{\delta, n-1}, \Phi^{\delta,n-1}),
	\]
	and hence, recalling that the solutions are obtained by constant interpolation, we have
	\[
	\mathcal{E}^\delta(\b{c}^{\delta, N}(t), \Phi^{\delta,N}(t)) + \int_{(n-1)\tau}^{n\tau} \int_\Omega M(\b{c}^{\delta, N}(s-\tau))\nabla \b{w}^{\delta, N}(s) \cdot \nabla \b{w}^{\delta,N}(s) \: \d x \: \d s \leq \mathcal{E}^\delta(\b{c}^{\delta, N}(t-\tau), \Phi^{\delta,N}(t-\tau))
	\]
	for all $t \in (0,T]$. Here, $n$ is chosen so that $t \in ((n-1)\tau, n\tau]$. Iterating the same argument yields (in a finite number of steps)
	\[
	\mathcal{E}^\delta(\b{c}^{\delta, N}(t), \Phi^{\delta,N}(t)) + \int_{0}^{n\tau} \int_\Omega M(\b{c}^{\delta, N}(s-\tau))\nabla \b{w}^{\delta, N}(s) \cdot \nabla \b{w}^{\delta,N}(s) \: \d x \: \d s \leq \mathcal{E}^\delta(\b{c}_0, \Phi_0),
	\]
	and the proof is complete observing that $n\tau \geq t$.
\end{proof}
\begin{cor} \label{cor:uniform1}
	There exists some $\delta_1 > 0$ and a constant $K_1 > 0$ independent of $\delta$, $N$ and $t$ such that 
	\[
	\sup_{t \in [0,T]} \|\b{c}^{\delta, N}(t)\|^2_V + \sup_{t \in [0,T]} \int_\Omega \sum_{i \in \mathcal{I}} \Psi^\delta(c_i^{\delta,N}(t)) \: \d x + \int_0^T \|\nabla\b{w}^{\delta, N}(t)\|^2_H \:\d t + \sup_{t \in [0,T]} \|\Phi^{\delta, N}(t)\|^2_V \leq K_1
	\]
	for all $\delta \in (0,\delta_1)$ and $N \in \mathbb{N}$.
\end{cor}
\begin{proof}
	In light of Lemma \ref{lem:energy}, we are only left to prove that the initial regularized energy is uniformly bounded. Indeed, then the result is a straightforward consequence of the definition of the energy functional, on account of Assumptions \ref{hyp:mobility}, \ref{hyp:epsilon} and the Poincaré inequality. The hypotheses on the initial state $\b{c}_0$ given in Assumption \ref{hyp:data}-\ref{hyp:initial} imply that $\b{c}_0$ takes values in the Gibbs simplex almost everywhere in $\Omega$. Furthermore, Assumption \ref{hyp:potential}-\ref{hyp:potder4} implies that (see also \cite{FG12})
	\[
	\Psi_\delta(s) \leq \Psi(s), \qquad \forall \: s \in (0,1), \quad \forall \: \delta \in (0,\varepsilon_2],
	\] 
	and therefore we have
	\review{
	\[
	\mathcal{E}^\delta(\b{c}_0, \Phi_0) \leq \mathcal{E}(\b{c}_0, \Phi_0) + \delta_1\|\b c_0\|_{\b L^4(\Omega;\mathbb{R}^3)}^4,
	\]
	}for all $\delta \in (0,\delta_1)$ with $\delta_1 := \min \{ \varepsilon_0, \varepsilon_1, \varepsilon_2\}$, and the proof is complete.
\end{proof} \noindent
Considering the linear interpolant, we can give an H\"{o}lder estimate in time on the approximated order parameter vector.
\begin{lem} \label{lem:holder}
	There exists a constant $K_2 > 0$ independent of $\delta$ and $N$ such that the linear interpolant satisfies
	\[
	\|\widetilde{\b{c}}^{\delta,N}(t_1) - \widetilde{\b{c}}^{\delta,N}(t_2)\|_H \leq K_2|t_1-t_2|^\frac 14
	\]
	for all $\delta \in (0,\delta_1)$, $N \in \mathbb{N}$ and $t_1,\,t_2 \in [0,T]$.
\end{lem}
\begin{proof}
	We follow, for instance, the approach in \cite{Gar03}. Without loss of generality, assume that $0 \leq t_1 < t_2 \leq T$. Choosing $\b{v} = \widetilde{\b{c}}^{\delta,N}(t_1) - \widetilde{\b{c}}^{\delta,N}(t_2)$ in the first equation given in Lemma \ref{lem:eul-lag} and observing that the finite difference equals the time derivative of the linear interpolant almost everywhere in $[0,T]$, one gets
	\[
	\left(\pd{\widetilde{\b{c}}^{\delta,N}}{t}(t), \: \widetilde{\b{c}}^{\delta,N}(t_1) - \widetilde{\b{c}}^{\delta,N}(t_2) \right)_H + \left( M(\widetilde{\b{c}}^{\delta, N}(t-\tau))\nabla \b{w}^{\delta, n}(t), \nabla (\widetilde{\b{c}}^{\delta,N}(t_1) - \widetilde{\b{c}}^{\delta,N}(t_2)	) \right)_H = 0,
	\]
	and integrating over $[t_1, t_2]$ yields
	\[
	\|\widetilde{\b{c}}^{\delta,N}(t_1) - \widetilde{\b{c}}^{\delta,N}(t_2)\|^2_H + 2\ii{t_1}{t_2}{\left( M(\widetilde{\b{c}}^{\delta, N}(t-\tau))\nabla \widetilde{\b{w}}^{\delta, n}(t), \nabla (\widetilde{\b{c}}^{\delta,N}(t_1) - \widetilde{\b{c}}^{\delta,N}(t_2)	) \right)_H}{t} = 0
	\]
	Recalling Corollary \ref{cor:uniform1}, we have
	\[
	\begin{split}
		\left| \ii{t_1}{t_2}{\left( M(\widetilde{\b{c}}^{\delta, N}(t-\tau))\nabla \widetilde{\b{w}}^{\delta, n}(t), \nabla (\widetilde{\b{c}}^{\delta,N}(t_1) - \widetilde{\b{c}}^{\delta,N}(t_2)	) \right)_H}{t}\right| \leq C\ii{t_1}{t_2}{\|\nabla \widetilde{\b{w}}^{\delta, n}(t)\|_H}{t} 
	\end{split}
	\]
	and the H\"{o}lder inequality in time yields the claim. The constant $K_2$ depends on $K_1$ and the parameter $m^*$ appearing in Assumption \ref{hyp:mobility}.
\end{proof} \noindent
Next, we give a uniform bound on the full $V$-norm of the chemical potentials.
\begin{lem} \label{lem:uniform2}
	There exists a constant $K_3 > 0$ independent of $\delta$ and $N$, but possibly depending on $T$, such that
	\[
	\int_0^T \|\b{w}^{\delta,N}(t)\|^2_V \: \d t \leq K_3
	\]
	for all $\delta \in (0,\delta_1)$ and $N \in \mathbb{N}$.
\end{lem}
\begin{proof}
	In light of Corollary \ref{cor:uniform1}, it is enough to give a uniform bound on the integral average of $\b{w}^{\delta,N}$. To this end, we exploit the same idea given in \cite{Garcke05}. Let us test the second equation of Lemma \ref{lem:eul-lag} by $\b{z} = \b{k} - \b{c}^{\delta,N}$ for some $\b{k} \in \Sigma$ to be defined later, possibly depending only on time. Moreover, we decompose 
	\[
	\b{w}^{\delta,N} = \left[\b{w}^{\delta,N} - \overline{\b{w}^{\delta,N}}\right] +\overline{\b{w}^{\delta,N}},
	\]
	arriving at
	\begin{equation} \label{eq:chemical0}
		\left(\left[\b{w}^{\delta,N} - \overline{\b{w}^{\delta,N}}\right] +\overline{\b{w}^{\delta,N}},\: \b{k} - \b{c}^{\delta,N} \right)_H = -\left( \Gamma \nabla \b{c}^{\delta, N} , \nabla \b{c}^{\delta, N} \right)_H  + \left( \b P\widetilde{\b{\mu}}^{\delta,N} , \b{k} - \b{c}^{\delta,N} \right)_H,
	\end{equation}
	with
	\[
	\widetilde{\b{\mu}}^{\delta,N} := \begin{bmatrix}
		\pd{F_\delta}{s_1}(c_A^{\delta, n}) + \alpha_{AA}\mathcal{N}(c_A^{\delta, N} - \overline{c_A^{\delta, N}}) + \alpha_{AB}\mathcal{N}(c_B^{\delta, N} - \overline{c_B^{\delta, N}}) + \dfrac{1}{2}\pd{\epsilon}{s_1}(c_A^{\delta, n})|\b{E}_0 - \nabla \Phi^{\delta, N}|^2 \\[0.5cm]
		\pd{F_\delta}{s_2}(c_B^{\delta, n}) + \alpha_{BA}\mathcal{N}(c_A^{\delta, N} - \overline{c_A^{\delta, N}}) + \alpha_{BB}\mathcal{N}(c_B^{\delta, N} - \overline{c_B^{\delta, N}}) + \dfrac{1}{2}\pd{\epsilon}{s_2}(c_B^{\delta, n})|\b{E}_0 - \nabla \Phi^{\delta, N}|^2 \\[0.3cm]
		\pd{F_\delta}{s_3}(c_S^{\delta, n})
	\end{bmatrix}.
	\]
	Now we give controls for the terms appearing to the right hand side of \eqref{eq:chemical0}. Owing to the convexity of $\Psi_\delta$, we have
	\begin{equation} \label{eq:chemical1}
		\begin{split}
			\sum_{i=1}^3 \ii{\Omega}{}{\theta_i \Psi_\delta(k_i)}{x} & \geq  \sum_{i=1}^3 \ii{\Omega}{}{\theta_i \Psi_\delta(c_i^{\delta,N})}{x} + \sum_{i=1}^3 \ii{\Omega}{}{\theta_i \Psi_\delta'(c_i^{\delta,N})(k_i - c_i^{\delta,N})}{x} \\
			& = \sum_{i=1}^3 \ii{\Omega}{}{\theta_i \Psi_\delta(c_i^{\delta,N})}{x} + \ii{\Omega}{}{\theta_i \b{P}\left[ \Psi_\delta'(\b{c^{\delta,N}}) \right] \cdot (\b{k} - {\b{c}^{\delta,N}})}{x},
		\end{split}
	\end{equation}
	where the last equality follows from the fact that $\b{k} - \b{c}^{\delta,N} \in T\Sigma$ almost everywhere. The bound given by \eqref{eq:boundI} and the fact that $\b{k} \in \Sigma \Rightarrow \|\b{k}\|_H \leq |\Omega|^\frac 12$ gives the straightforward control
	\begin{equation} \label{eq:chemical2}
		\left|\ii{\Omega}{}{\b{P}\td{\widetilde{I}}{\b s}(\b{c}^{\delta,N}) \cdot (\b{k} - {\b{c}^{\delta,N}})}{x}\right| \leq C(1 + \|\b{c}^{\delta,N}\|^2_H).
	\end{equation}
	For the nonlocal terms, observe that for any $i$ and $j$ in $\{1,2\}$, we have by continuity of $\mathcal{N}$ and the Poincarè and Cauchy-Schwarz inequalities
	\begin{equation} \label{eq:chemical3}
		([\b{P}\mathcal{N}(\b c^{\delta, N} - \overline{\b{c}^{\delta, N}})]_i,\, k_{j} - c_j^{\delta, N})_H \leq C\|\nabla c_i^{\delta, N}\|_H(1 +\|c_j^{\delta, N}\|_H),
	\end{equation}
	whereas by Assumption \ref{hyp:epsilon}
	\begin{equation} \label{eq:chemical4}
		\left( \dfrac{1}{2}\b{P}\left[\td{\epsilon}{\b s}(\b{c}^{\delta, N})|\b{E}_0 - \nabla \Phi^{\delta, N}|^2\right]_i,\, k_{i} - c_i^{\delta, N} \right)_H \leq C\left( 1+ \|k_{i}\|_{L^\infty(\Omega)} \right).
	\end{equation}
	In \eqref{eq:chemical4}, we understand that $\frac{\partial \epsilon}{\partial s_3} = 0$. Let us now give a precise form to the vector $\b{k}$. The assumptions on the initial condition given in Assumption \ref{hyp:data}-\ref{hyp:initial} also yield the existence of a value $\rho \in (0,1)$ such that
	\[
	\rho \leq \dfrac{1}{|\Omega|}\int_\Omega c_i^{\delta,N} \: \d x \leq 1-\rho,
	\]
	for all $i \in \mathcal{I}$. Fixed any $p \in \{1,2,3\}$ and any $l \in \{1,2,3\}$, let
	\[
	\b{k} := \overline{\b c^{\delta,N}} + \rho \operatorname{sign}\left(\overline{w_p^{\delta,N}} -\overline{w_l^{\delta,N}} \right)(\b{e}_p - \b{e}_l) \in \b{G}\subset \Sigma.	
	\]
	Since the Gibbs simplex is contained in the unitary ball of $\mathbb{R}^3$, we have $\|\b{k}\|_{L^\infty(\Omega; \mathbb{R}^3)} \leq 1$.
	Collecting the above and \eqref{eq:chemical1}-\eqref{eq:chemical4} into \eqref{eq:chemical0}, we have, as a consequence of the Poincaré inequality and Corollary \ref{cor:uniform1}, 
	\begin{equation}
		\left( \overline{\b{w}^{\delta,N}},\: \b{k} - \b{c}^{\delta,N} \right)_H \leq C\left( 1 + \|\nabla \b{w}^{\delta,N}\|_H \right)
	\end{equation}	
	and the definition of $\b{k}$ gives (cfr. \cite[Proof of Lemma 4.3]{Garcke05})
	\[
	\ii{0}{T}{|\overline{\b{w}^{\delta,N}(s)}|^2}{s}\leq C,
	\]
	and the proof is complete.
\end{proof}
\subsection{Higher integrability properties.} In this subsection, we address higher-order estimates for the regularized solutions.
\begin{lem} \label{lem:higher2}
	There exist a constant $K_4 > 0$ independent of $\delta$, $N$ and $T$ as well as a constant $p = p(\Omega,\epsilon_*,\epsilon^*) > 2$ such that
	\[
	\sup_{t \in [0,T]}\|\Phi^{\delta,N}(t)\|_{W^{1,p}(\Omega)} \leq K_4
	\]
	for all $\delta \in (0,\delta_1)$ and $N \in \mathbb{N}$.
\end{lem}
\begin{proof}
	The result follows from applying, for instance, \cite[Theorem 1]{Meyers63} (see also \cite[Theorem 1]{Groger89}), implying the estimate
	\[
	\|\nabla \Phi\|_{L^p(\Omega;\mathbb{R}^3)} \leq C\|\b{E}_0\|_{L^p(\Omega;\mathbb{R}^3)}
	\]
	for some $p > 2$ only depending on $\Omega$ and the ellipticity constants of $\epsilon$. The thesis follows from the Poincaré inequality and taking supremums over time.
\end{proof} \noindent
Now we give a higher-order integrability result on the derivatives of the regularized singular potential $\Psi_\delta$, following \cite[Lemma 5.1]{Garcke05}.
\begin{lem} \label{lem:higher3}
	There exists a constant $K_5 > 0$ independent of $\delta$ and $N$, but possibly depending on $T$, such that
	\[
	\ii{0}{T}{\|\Psi'_\delta(c_i^{\delta,N}(t))\|_{L^\frac q2(\Omega)}^\frac{q}{2}}{t} \leq K_5
	\]
	where $q = \min\{p,4\}$, for all $i \in \mathcal{I}$, $N \in \mathbb{N}$ and $\delta \in (0,\delta_2)$, where $\delta_2 := \min\{\frac 13, \delta_1\}$.
\end{lem}
\begin{proof}
	For any arbitrary but fixed $b > 0$, let us define the function
	\[
	\Xi_\delta : \mathbb{R} \to \mathbb{R}, \qquad \Xi_\delta(s) :=  \Psi'_\delta(s) |\Psi'_\delta(s)|^{b-1},
	\]
	with the understanding that, if $b \in (0,1)$, then $\Xi_\delta$ vanishes at the zeros of $\Psi_\delta'$. If $b \in (0,1)$, then the map $x \mapsto x|x|^{b-1}$ is not differentiable at zero, and therefore the function $\Xi_\delta$ turns out not to be differentiable at the zeros of $\Psi'_\delta$. Therefore, we smooth it considering a family of functions $\{\Xi_{\delta, \omega}\}_{\omega \in [0,1]}$ with the following properties:
	\begin{enumerate}[label=(\arabic*)]
		\item the function $\Xi_{\delta, \omega} \in C^1(\mathbb{R})$ and is monotone increasing;
		\item the equality $\Xi_{\delta, \omega} = \Xi_{\delta}$ holds outside the interval $[0,1]$;
		\item the convergence $\Xi_{\delta, \omega} \to \Xi_{\delta} \text{ in }C^0(\mathbb{R})$ holds as $ \omega \to 0^+$.
	\end{enumerate}
	In the following, we shall work within the smoothing scheme even if $b \geq 1$.
	Looking back to Lemma \ref{lem:eul-lag} and taking into account Lemma \ref{lem:higher2}, we can enlarge the set of admissible test functions $\b{z}$ to the space $V \cap L^\frac{p}{p-2}(\Omega;\mathbb{R}^d)$ (after extending the solutions in time).	Observe that we have $\frac{p}{p-2} > 1$. Let us take as test function $\b{z}$ the vector such that
	\[
	z_i = \Xi_{\delta, \omega}(c_i^{\delta,N}(t)),
	\] 
	which is admissible provided that $\frac{bp}{p-2} \leq 2$ owing to the fact that $\Psi_\delta'$ is bounded above by a cubic and the Sobolev embedding $H^1(\Omega) \hookrightarrow L^6(\Omega)$. Recalling \eqref{eq:projector}, this results in \review{
	\begin{multline} \label{eq:singpot1}
		\sum_{i=1}^{3} \left( {w}_i^{\delta,N},  z_i \right)_H = \sum_{i=1}^{3} \left[ \gamma_i\left( \nabla c_i^{\delta,N}, \nabla z_i \right)_H + \theta_i\left(\Psi'_\delta(c_i^{\delta,N}) - \dfrac{1}{3}\sum_{k=1}^3\Psi'_\delta(c_k^{\delta,N}),\, z_i  \right)_H \right] \\ + \left(\b{P}\td{\widetilde{I}}{\b s}(\b{c}^{\delta,N}), \b{z}\right)_H + \sum_{j=1}^2\sum_{k=1}^2 \alpha_{jk}(\mathcal{N}(c_j^{\delta,N} - \overline{c_j^{\delta,N}}) - \dfrac{1}{3}\sum_{l=1}^3\mathcal{N}(c_l^{\delta,N} - \overline{c_l^{\delta,N}}), z_k)_H \\
		+ (\b{P}\b{\epsilon}'|\b{E}_0 - \nabla \Phi^{\delta,N}|^2, \b{z})_H + \sum_{m=1}^{3} \left(4\delta (c_m^{\delta, N})^3 - \dfrac{4\delta}{3}\sum_{n=1}^{3}(c_n^{\delta, N})^3,\,z_m \right)_H,
	\end{multline}
	where we set 
	\[
	\b{\epsilon}' := \left[ \pd{\epsilon}{s_1}(c_A^{\delta,N}),\, \pd{\epsilon}{s_2}(c_B^{\delta,N}),\, 0 \right]^T.
	\]
	Observe that since $\Xi_{\delta, \omega}$ is monotone increasing, the first term to the right hand side is positive. Next, we have, following computations in \cite[Proof of Lemma 5.1]{Garcke05}, 
	\[
	\begin{split}
		&\sum_{i=1}^3\theta_i\left[\Psi'_\delta(c_i^{\delta,N}) - \dfrac{1}{3}\sum_{k=1}^3\Psi'_\delta(c_k^{\delta,N}) \right] \Xi_{\delta, \omega}(c_i^{\delta,N}) \\ & \hspace{2cm}= \dfrac{1}{3}\sum_{\substack{i,k = 1 \\ i < k}}^3\theta_i\left[\Psi'_\delta(c_i^{\delta,N}) - \Psi'_\delta(c_k^{\delta,N}) \right]\left[ \Xi_{\delta, \omega}(c_i^{\delta,N}) - \Xi_{\delta, \omega}(c_k^{\delta,N})\right] \geq 0,
	\end{split}
	\]
	once again exploiting monotonicity of both $\Xi_{\delta, \omega}$ and $\Psi'_\delta$ (the latter following from convexity). The same goes for the perturbation term arising from the regularization scheme (i.e., the term on the last line of \eqref{eq:singpot1}), as the map $x \mapsto x^3$ is monotone increasing in $\mathbb{R}$.} Therefore, \eqref{eq:singpot1} implies
	\[
	\begin{split}
		& \dfrac{1}{3}\sum_{\substack{i,k = 1 \\ i < k}}^3\theta_i\left( \Psi'_\delta(c_i^{\delta,N}) - \Psi'_\delta(c_k^{\delta,N}),\, \Xi_{\delta, \omega}(c_i^{\delta,N}) - \Xi_{\delta, \omega}(c_k^{\delta,N})\right)_H\\
		& \hspace{0.5cm} \leq \sum_{i=1}^{3} \left( {w}_i^{\delta,N},  \Xi_{\delta, \omega}(c_i^{\delta,N}) \right)_H - \left(\b{P}\td{\widetilde{I}}{\b s}(\b{c}^{\delta,N}), \b{z}\right)_H + (\b{P}\b{\epsilon}'|\b{E}_0 - \nabla \Phi^{\delta,N}|^2, \b{z})_H \\
		& \hspace{6cm}+ \sum_{j=1}^2\sum_{k=1}^2 \alpha_{jk}(\mathcal{N}(c_j^{\delta,N} - \overline{c_j^{\delta,N}}) - \dfrac{1}{3}\sum_{l=1}^3\mathcal{N}(c_l^{\delta,N} - \overline{c_l^{\delta,N}}), \Xi_{\delta, \omega}(c_k^{\delta,N}))_H, \\
		& \hspace{0.5cm} \leq C \left[ \max_{i \in \mathcal{I}} \|\Xi_{\delta, \omega}(c_i^{\delta,N})\|_H\left(1 + \|\b{w}^{\delta,N}\|_H + \|\b{c}^{\delta,N}\|_V\right) + \|\b{E}_0 - \nabla \Phi^{\delta,N}\|_{L^p(\Omega; \mathbb{R}^d)} \max_{i \in \mathcal{I}}\|\Xi_{\delta, \omega}(c_i^{\delta,N})\|_{L^\frac{p}{p-2}}\right].
	\end{split}
	\]
	Integrating the inequality above in time, and exploiting Corollary \ref{cor:uniform1} jointly with Lemma \ref{lem:higher2}, we can now follow line by line \cite[Proof of Lemma 5.1]{Garcke05} to yield 
	\[
	\int_0^T \int_\Omega \max_{i \in \mathcal{I}}|\Xi_{\delta, \omega}(c_i^{\delta,N})|^{b+1} \: \d x \: \d t \leq \dfrac{3\sigma}{2\theta_*}\int_0^T \int_\Omega \max_{i \in \mathcal{I}}|\Xi_{\delta, \omega}(c_i^{\delta,N})|^{\frac{bp}{p-2}} \: \d x \: \d t + C_\sigma\]
	for any $\sigma > 0$. The constant $C_\sigma$ depends on $\sigma$. Choosing
	\[
	b = \dfrac{p}{2} - 1 \Rightarrow \dfrac{bp}{p-2} = b + 1 = \dfrac{p}{2}, \qquad \sigma < \dfrac{2\theta_*}{3}
	\]
	yields the claim, keeping in mind that here we must have
	\[
	\dfrac{bp}{p-2} \leq 2 \Leftrightarrow p \leq 4,
	\]
	and the proof is complete.
\end{proof}
\begin{cor} \label{cor:higher4}
	There exists a constant $K_6 > 0$ independent of $\delta$ and $N$, but possibly depending on $T$, such that
	\[
	\ii{0}{T}{\|\b{c}^{\delta,N}(t)\|^\frac q2_{W^{2,\frac q2}(\Omega;\mathbb{R}^d)}}{t} \leq K_6
	\]
	for all $N \in \mathbb{N}$ and $\delta \in (0,\delta_2)$, 	where $q = \min \{p,4\}$.
\end{cor}
\begin{proof}
	Follows by comparison in \eqref{eq:strongform2} and elliptic regularity.
\end{proof}
\subsection{Passage to the limit.} The previously proven uniform estimates enable us to retrieve a solution to the original problem. First of all, we sum up the results obtained throughout the previous subsection. We showed that for any $\delta \in (0,\delta_2)$ and $N \in \mathbb{N}$
\begin{align*}
	\{\b{c}^{\delta,N}\}_{\delta,N} & \text{ is uniformly bounded in } L^\infty(0,T;V) \cap L^\frac q2(0,T;W^{2,\frac q2}(\Omega;\mathbb{R}^d)), \\
	\{\widetilde{\b{c}}^{\delta,N}\}_{\delta,N} & \text{ is uniformly $\alpha$-H\"{o}lder-equicontinuous for } \alpha \in \left(0,\frac 14\right], \\
	\{\b{w}^{\delta,N}\}_{\delta,N} & \text{ is uniformly bounded in } L^2(0,T;V), \\
	\{\Phi^{\delta,N}\}_{\delta,N} & \text{ is uniformly bounded in } L^\infty(0,T;W^{1,p}(\Omega) \cap V_0).
\end{align*}
Let us now pass to the limit as $N \to +\infty$, i.e., as the time step $\tau \to 0^+$.
\begin{prop} \label{prop:limitN1}
	Let $\delta \in (0,\delta_2)$ be fixed. There exists a triple $(\b{c}^\delta,\, \b{w}^\delta,\,\Phi^\delta)$ such that
	\begin{align*}
		\b{c}^\delta & \in C^{0,\reviewF{\frac 14}}([0,T]; H) \cap L^\infty(0,T;V) \cap L^\frac q2(0,T;W^{2,\frac q2}(\Omega;\mathbb{R}^d)) \cap H^1(0,T;V^*), \\
		\b{w}^\delta & \in L^2(0,T;V), \\
		\Phi^\delta & \in L^\infty(0,T;V_0),
	\end{align*} 
	and the following hold:
	\begin{enumerate}[label = (\roman*)] \itemsep 0.2em
		\item $\widetilde{\b{c}}^{\delta, N}\to \b{c}^\delta$ in $C^{0,\alpha}([0,T]; H)$ for all $\alpha \in (0,\frac 14)$;
		\item $\b{c}^{\delta, N}\to \b{c}^\delta$ in $L^\infty(0,T; H)$ and almost everywhere in $\Omega \times (0,T)$;
		\item $\b{c}^{\delta, N} \rightharpoonup \b{c}^\delta$ weakly$^*$ in $L^\infty(0,T; V)$ and weakly in $L^\frac q2(0,T;W^{2,\frac q2}(\Omega;\mathbb{R}^d))$;
		\item $\Psi_\delta(c_i^{\delta, N}) \to \Psi_\delta(c_i^{\delta})$ in $L^1(\Omega \times (0,T))$ for all $i \in \mathcal{I}$;
		\item $\pd{\widetilde{I}}{s_i}(c_i^{\delta,N}) \to \pd{\widetilde{I}}{s_i}(c_i^{\delta})$ in $L^1(\Omega \times (0,T))$ for all $i \in \mathcal{I}$;
		\item $\b{w}^{\delta,N} \rightharpoonup \b{w}^{\delta}$ weakly in $L^2(0,T;V)$;
		\item $\Phi^{\delta,N} \to \Phi^\delta$ in $L^2(0,T;V_0)$.
	\end{enumerate}
	as $N \to +\infty$, up to a subsequence.
\end{prop}
\begin{proof}
	Claim $(i)$ follows from the Ascoli--Arzelà theorem. Let us show that it implies $(ii)$ and the that the limits coincide. Let $t = \beta n\tau + (1-\beta)(n-1)\tau$ for some $n \in \{1,\dots,N\}$ and $\beta \in [0,1]$. Then,
	\[
	\|\widetilde{\b{c}}^{\delta,N}(t) - \b{c}^{\delta,N}(t)\|_H = (1-\beta)\|\b{c}^{\delta,N}(t) - \b{c}^{\delta,N}(t-\tau)\|_H \leq C\tau^\frac{1}{4}
	\]
	by Lemma \ref{lem:holder}. Claim $(iii)$ follows from the Banach--Alaoglu theorem and by reflexivity of the Banach space $L^\frac q2(0,T;W^{2,\frac q2}(\Omega;\mathbb{R}^d))$ for $q > 2$. The same goes for $(vi)$. As for $(iv)$ and $(v)$, they are consequences of \ref{hyp:potential}-\ref{hyp:potder1-2}, the Vitali convergence theorem, as well as the bound \eqref{eq:boundIder}. We are left to show $(vii)$. The uniform bound in $L^\infty(0,T;W^{1,p}(\Omega))$ entails a uniform bound in $L^2(0,T;V)$, hence weak convergence is immediately achieved by reflexivity. Consider back the last equation in the weak formulation given in Lemma \ref{lem:eul-lag}, and choose $\Lambda = \Phi^{\delta,N} - \Phi^\delta$. By suitably summing and subtracting $\nabla \Phi^\delta$, we obtain
	\[
	\left(\epsilon(c_A^{\delta,N},\,c_B^{\delta,N})(\b E_0 - \nabla \Phi^\delta + \nabla \Phi^\delta - \nabla \Phi^{\delta,N}),\,\nabla\Phi^{\delta,N} - \nabla\Phi^\delta\right)_H = 0,
	\]
	or, equivalently,
	\[
	\left(\epsilon(c_A^{\delta,N},\,c_B^{\delta,N})(\nabla\Phi^{\delta,N} - \nabla\Phi^\delta),\,\nabla\Phi^{\delta,N} - \nabla\Phi^\delta\right)_H = \left(\epsilon(c_A^{\delta,N},\,c_B^{\delta,N})(\b E_0 - \nabla \Phi^\delta),\,\nabla\Phi^{\delta,N} - \nabla\Phi^\delta\right)_H.
	\]
	This yields, owing to Assumption \ref{hyp:epsilon},
	\begin{equation} \label{eq:wtsc}
		\begin{split}
			\epsilon_*\|\nabla \Phi^{\delta,N} - \nabla \Phi^\delta\|^2_H \leq   \left(\epsilon(c_A^{\delta,N},\,c_B^{\delta,N})(\b E_0 - \nabla \Phi^\delta),\,\nabla\Phi^{\delta,N} - \nabla\Phi^\delta\right)_H.
		\end{split}
	\end{equation}
	As $\nabla \Phi^{\delta,N} \rightharpoonup \nabla \Phi^\delta$ weakly in $L^2(0,T;H)$ as $N \to +\infty$, we infer $\nabla \Phi^{\delta,N}-\nabla \Phi^\delta\rightharpoonup 0$ in $L^2(0,T;H)$ as $N \to +\infty$. By the mean value and the Lebesgue dominated convergence theorem, we easily show that, since $c_i^{\delta, N} \to c_i^\delta$ in $L^2(0,T;H)$ and almost everywhere in $\Omega \times (0,T)$, then
\[
\epsilon(c_A^{\delta,N},\,c_B^{\delta,N})(\b E_0 - \nabla \Phi^\delta) \to \epsilon(c_A^{\delta},\,c_B^{\delta})(\b E_0 - \nabla \Phi^\delta) \quad \text{in } L^2(0,T;H), \quad \text{as } N \to +\infty.
\]
	This implies, by the general theory, that the right hand side is \eqref{eq:wtsc} tends to zero as $N \to +\infty$. The proof is complete.
\end{proof} \noindent
Now we are in a position to identify the limit triplet $(\b{c}^\delta,\, \b{w}^\delta,\,\Phi^\delta)$ as a solution to a continuous regularized version of \eqref{eq:strongform2}-\eqref{eq:conditions2}, obtained as a limit version of the Euler--Lagrange equations given by Lemma \ref{lem:eul-lag}.
\begin{prop} \label{prop:limitN2}
	The triplet $(\b{c}^\delta,\, \b{w}^\delta,\,\Phi^\delta)$ satisfies
	\begin{equation*}
		\begin{cases}
			\left\langle \pd{\b{c}^\delta}{t}, \: \b{v} \right\rangle_{V^*, V} + \left( M(\b{c}^{\delta})\nabla \b{w}^{\delta}, \nabla \b{v} \right)_H = 0, & \quad \forall \: \b{v} \in V, \\[0.5cm]
			\left(\b{w}^{\delta}, \b{z} \right)_H = \left( \Gamma \nabla \b{c}^{\delta} , \nabla \b{z} \right)_H + \left( \b P\widetilde{\b{\mu}}^{\delta} , \b{z} \right)_H, & \quad \forall \: \b{z} \in V \cap L^\infty(\Omega;\mathbb{R}^3), \\[0.5cm]
			(\epsilon({c}_A^{\delta},\,{c}_B^{\delta})(\b{E}_0 - \nabla \Phi^{\delta}), \nabla \Lambda)_H = 0 & \quad \forall \: \Lambda \in V_0,
		\end{cases}
	\end{equation*}
	almost everywhere in $[0,T]$, where
	\[
	\widetilde{\b{\mu}}^{\delta} := \begin{bmatrix}
		\pd{F_\delta}{s_1}(c_A^{\delta}) + \alpha_{AA}\mathcal{N}(c_A^{\delta} - \overline{c_A^{\delta}}) + \alpha_{AB}\mathcal{N}(c_B^{\delta} - \overline{c_B^{\delta}}) + \dfrac{1}{2}\pd{\epsilon}{s_1}(c_A^{\delta})|\b{E}_0 - \nabla \Phi^{\delta}|^2 \\[0.5cm]
		\pd{F_\delta}{s_2}(c_B^{\delta}) + \alpha_{BA}\mathcal{N}(c_A^{\delta} - \overline{c_A^{\delta}}) + \alpha_{BB}\mathcal{N}(c_B^{\delta} - \overline{c_B^{\delta}}) + \dfrac{1}{2}\pd{\epsilon}{s_2}(c_B^{\delta})|\b{E}_0 - \nabla \Phi^{\delta}|^2 \\[0.3cm]
		\pd{F_\delta}{s_3}(c_S^{\delta})
	\end{bmatrix}.
	\]
\end{prop}
\begin{proof}
	Let us address first the evolution equation. Recall that for the linear interpolant, as shown in the proof of Lemma \ref{lem:holder}, it holds
	\[
	\left(\pd{\widetilde{\b{c}}^{\delta,N}}{t}(t), \: \b{v} \right)_H + \left( M(\widetilde{\b{c}}^{\delta, N}(t-\tau))\nabla \b{w}^{\delta, N}(t), \nabla \b{v} \right)_H = 0,
	\]
	almost everywhere in $[0,T]$ and for all $\b{v} \in V$. The weak convergence properties of $\b{c}^{\delta,N}$ are enough to pass to the limit in the first term, while for the second term we argue as follows. Assumption \ref{hyp:mobility} and the weak and almost everywhere convergences given by Proposition \ref{prop:limitN1} allow us to use the Lebesgue dominated convergence theorem to find 
	\[
	\lim_{N\to +\infty} \int_0^T \left( M(\widetilde{\b{c}}^{\delta, N}(t-\tau))\nabla \b{w}^{\delta, N}(t), \nabla \b{v} \right)_H \: \d t = \int_0^T \left( M(\b{c}^{\delta}(t))\nabla \b{w}^{\delta}(t), \nabla \b{v} \right)_H \: \d t,
	\] 
	for all $\b v \in L^2(0,T;V)$, recalling that $N \to +\infty$ implies $\tau \to 0^+$. Hence,
	\[
	\int_0^T \left\langle \pd{ \b c^\delta(t)}{t}, \b v \right\rangle_{V^*,V} + \left( M(\b{c}^{\delta}(t))\nabla \b{w}^{\delta}(t), \nabla \b{v} \right)_H \: \d t = 0
	\] 
	for all $\b v \in L^2(0,T;V)$. Choosing $\b v(\b x,t) = \b z(\b x) \mathbbm{1}_{(t,t+h)}(t)$, where $\b z \in V$ and $t,\,h$ are such that $(t,t+h) \subset [0,T]$ yields, by dividing the result by $h$ and applying the Lebesgue differentiation theorem, the fact that 
\[
\left\langle\pd{ \b c^\delta(t)}{t}, \b v \right\rangle_{V^*,V} + \left( M(\b{c}^{\delta}(t))\nabla \b{w}^{\delta}(t), \nabla \b{v} \right)_H = 0 \quad \text{for a.a. }t \in [0,T].
\]
A similar argument also works for the third equation. Finally, we address the equation for the projected chemical potentials. After an integration on $[0,T]$, given the convergence properties of $\b{w}^{\delta, N}$ and $\b{c}^{\delta,N}$, it is immediate to pass to the limit to the left hand side, in the gradient terms and in the linear nonlocal terms. For the electrical term, we can exploit the strong convergence of $\Phi^{\delta,N} \to \Phi^\delta$ in $L^2(0,T;V_0)$ and the Lebesgue dominated convergence theorem. Finally, for the local interaction term we can exploit the growth bound given by \eqref{eq:boundIder}, while for the regularization, for any $i \in \mathcal{I}$, we can use the Vitali convergence theorem and Corollary \ref{cor:uniform1}.
\end{proof} \noindent
In order to prove the existence of weak solutions for the original problem, we shall prove the analogous versions of Propositions \ref{prop:limitN1} and \ref{prop:limitN2} in the limit $\delta \to 0^+$. Since the proofs of these results are in many regards similar to the ones given above, we highlight only the main differences. Let us point out that, since the limit system has the same structure of the original one, uniform estimates can be proven for the triplets $(\b{c}^\delta,\, \b{w}^\delta,\,\Phi^\delta)$ exactly as shown throughout Subsection \ref{ssec:unifest}, with minor modifications.
\begin{prop} \label{prop:limitd1}
	There exists a triple $(\b{c},\, \b{w},\,\Phi)$ such that
	\begin{align*}
		\b{c} & \in C^{0,\reviewF{\frac 14}}([0,T]; H) \cap L^\infty(0,T;V) \cap L^\frac q2(0,T;W^{2,\frac q2}(\Omega;\mathbb{R}^3)) \cap H^1(0,T;V^*), \\
		\b{w} & \in L^2(0,T;V), \\
		\Phi & \in L^\infty(0,T;V_0),
	\end{align*} 
	and the following hold:
	\begin{enumerate}[label = (\roman*)] \itemsep 0.2em
		\item $\b{c}^{\delta}\to \b{c}$ in $C^{0,\alpha}([0,T]; H)$ for all $\alpha \in (0,\frac 14)$;
		\item $\b{c}^{\delta}\to \b{c}$ in $L^\infty(0,T; H)$ and almost everywhere in $\Omega \times (0,T)$;
		\item $\b{c}^{\delta} \rightharpoonup \b{c}$ weakly$^*$ in $L^\infty(0,T; V)$ and weakly in $L^\frac q2(0,T;W^{2,\frac q2}(\Omega;\mathbb{R}^3))$;
		\item $\Psi_\delta(c_i^{\delta}) \to \Psi(c_i)$ in $L^1(\Omega \times (0,T))$ for all $i \in \mathcal{I}$;
		\item $\pd{\widetilde{I}}{s_i}(c_i^{\delta}) \to \pd{\widetilde{I}}{c_i}$ in $L^1(\Omega \times (0,T))$ in $L^1(\Omega \times (0,T))$ for all $i \in \mathcal{I}$;
		\item $\b{w}^{\delta} \rightharpoonup \b{w}$ weakly in $L^2(0,T;V)$;
		\item $\Phi^{\delta} \to \Phi$ in $L^2(0,T;V_0)$.
	\end{enumerate}
	as $\delta \to 0^+$, up to a subsequence. Moreover, $c_i(\b x, t) \in (0,1)$ for almost any $(\b x, t) \in \Omega \times (0,T)$ and for all $i \in \mathcal{I}$.
\end{prop}
\begin{proof}
	The results follow mainly as displayed in Proposition \ref{prop:limitN1}, except for $(iv)$ and $(v)$. By the Fatou lemma we have
	\begin{equation} \label{eq:liminf}
		\ii{\Omega}{}{\liminf_{\delta\to 0^+} |\Psi'_\delta(c_i^\delta)|^\frac q2}{x} \leq \liminf_{\delta\to 0^+} \ii{\Omega}{}{ |\Psi'_\delta(c_i^\delta)|^\frac q2}{x}.
	\end{equation}
	Let us compute the left hand side for fixed $i \in \mathcal{I}$. Let us define the partition of $\Omega \times (0,T)$
	\[
	\Omega \times (0,T) = Q_+ \cup Q_- \cup N
	\]
	with
	\begin{align*}
		Q_+ := \{(\b{x},t) \in \Omega \times (0,T): c_i^\delta(\b{x},t) \to c_i(\b{x},t) > 0\} \\
		Q_- := \{(\b{x},t) \in \Omega \times (0,T): c_i^\delta(\b{x},t) \to c_i(\b{x},t) \leq 0\}
	\end{align*}
	and $N$ being a null set with respect to the $(d+1)$--Lebesgue measure. On elements of $Q_+$, since $\Psi_\delta'(s)= \Psi'(s)$ for every $s \geq \delta$, the pointwise convergence $\Psi'_\delta(c_i^\delta) \to \Psi'(c_i)$ follows. On elements of $Q_-$, instead, by the structure of the first derivative of the regularized potential, for any fixed $\delta$, the following alternative holds
	\[
	\begin{cases}
		\Psi'_\delta(c_i^\delta) = \Psi'(c_i^\delta), & \quad \text{ if } c_i^\delta \geq \delta, \\
		\Psi'_\delta(c_i^\delta) = \Psi'(\delta)- \displaystyle\sum_{i=1}^{3} \dfrac{\Psi^{(i+1)}(\delta)}{i!}(\delta-c_i^\delta)^i, & \quad \text{ if } c_i^\delta < \delta.
	\end{cases} 
	\]
	If $\delta$ is fixed to be sufficiently small, then the signs of the derivatives of $\Psi$ evaluated at $\delta$ are known (see Assumption \ref{hyp:potential}). In particular, we have
	\[
	\begin{cases}
		|\Psi'_\delta(c_i^\delta)| = |\Psi'(c_i^\delta)|, & \quad \text{ if } c_i^\delta \geq \delta, \\
		|\Psi'_\delta(c_i^\delta)| \geq |\Psi'(\delta)|, & \quad \text{ if } c_i^\delta < \delta.
	\end{cases} 
	\]
	The two relations can be summarized to
	\[
	|\Psi'_\delta(c_i^\delta)| \geq |\Psi'(\max\{c_i^\delta,\,\delta\})|,
	\]
	and it is easy to show that, as the right hand side diverges at $+\infty$ when $\delta \to 0^+$, it must hold that $|\Psi'_\delta(c_i^\delta)|$ diverges to $+\infty$ on $Q_-$ when $\delta \to 0^+$. Therefore, the $(d+1)$--Lebesgue measure of $Q_-$ must be 0. Otherwise, Lemma \ref{lem:higher3} cannot hold independently of $\delta$. This proves that $\b{c} \in \b{G}$ almost everywhere in $\Omega \times (0,T)$ and thus $c_i \in (0,1)$ for all $i \in \mathcal{I}$ almost everywhere in $\Omega \times (0,T)$. In particular, we have equintegrability from \eqref{eq:liminf} and pointwise convergence
	\[
	\Psi_\delta'(c_i^\delta) \to \Psi'(c_i)
	\]
	almost everywhere. The Vitali convergence theorem leads the result. Claim $(v)$ is proven again exploiting \eqref{eq:boundI}, and $c_i \in (0,1)$ yields $\widetilde{I} = I$.
\end{proof} \noindent
Finally, we recover a solution to the original problem.
\begin{prop} \label{prop:limitd2}
	The triplet $(\b{c},\, \b{w},\,\Phi)$ satisfies
	\begin{equation*}
		\begin{cases}
			\left\langle \pd{\b{c}}{t}, \: \b{v} \right\rangle_{V^*, V} + \left( M(\b{c})\nabla \b{w}, \nabla \b{v} \right)_H = 0, & \quad \forall \: \b{v} \in V, \\[0.5cm]
			\left(\b{w}, \b{z} \right)_H = \left( \Gamma \nabla \b{c} , \nabla \b{z} \right)_H + \left( \b P\widetilde{\b{\mu}} , \b{z} \right)_H, & \quad \forall \: \b{z} \in V \cap L^\infty(\Omega;\mathbb{R}^3), \\[0.5cm]
			(\epsilon({c}_A,\,{c}_B)(\b{E}_0 - \nabla \Phi), \nabla \Lambda)_H = 0 & \quad \forall \: \Lambda \in V_0,
		\end{cases}
	\end{equation*}
	where
	\[
	\widetilde{\b{\mu}} := \begin{bmatrix}
		\theta_A\Psi'(c_A) + \pd{I}{c_A} + \alpha_{AA}\mathcal{N}(c_A - \overline{c_A}) + \alpha_{AB}\mathcal{N}(c_B - \overline{c_B}) + \dfrac{1}{2}\pd{\epsilon}{c_A}|\b{E}_0 - \nabla \Phi|^2 \\[0.5cm]
		\theta_B\Psi'(c_B) + \pd{I}{c_B} + \alpha_{BA}\mathcal{N}(c_A - \overline{c_A}) + \alpha_{BB}\mathcal{N}(c_B - \overline{c_B}) + \dfrac{1}{2}\pd{\epsilon}{c_B}|\b{E}_0 - \nabla \Phi|^2 \\[0.3cm]
		\theta_S\Psi'(c_S) + \pd{I}{c_S}
	\end{bmatrix}.
	\]
\end{prop}
\begin{proof}
	The proof follows line by line the one of Proposition \ref{prop:limitN2}.	
\end{proof} \noindent
Finally, we state two a-posteriori regularity results.

\begin{cor}
	There exist a constant $K_7 > 0$ independent of $T$ as well as a constant $p = p(\Omega,\epsilon_*,\epsilon^*) > 2$ such that
		\[
		\sup_{t \in [0,T]}\|\Phi(t)\|_{W^{1,p}(\Omega)} \leq K_7.
		\]
\end{cor}
\begin{proof}
	This is an iteration of the proof of Lemma \ref{lem:higher2}.
\end{proof}
\begin{prop} \label{prop:higher4}
	If $p \geq 4$, then $\b{c} \in L^4(0,T;H^2(\Omega;\mathbb{R}^3))$. If also $\b{E}_0 \in W^{1,\alpha}(\Omega;\mathbb{R}^d)$, with $\alpha > \frac 65$, then $\Phi \in L^4(0,T;W^{2,\beta}(\Omega))$ with
	\[
	\beta := \min\left\{\frac{6\alpha}{6+\alpha}, \frac{6p}{6+p}\right\} > 1.
	\]
\end{prop}
\begin{proof}
	Without loss of generality, we assume $p = 4$. Obviously, $q = \min \{p,4\} = 4$. Therefore, we have by Proposition \ref{prop:limitd1}
	\begin{align*}
		\Phi & \in L^\infty(0,T;W^{1,4}(\Omega)), \\
		\Psi'(c_i) & \in L^2(0,T;H) \quad \forall \: i \in \mathcal{I}, \\
		\b{c} & \in L^2(0,T;H^2(\Omega;\mathbb{R}^3)),
	\end{align*}
	implying, in turn, that the equation for $\b{w}$ can be interpreted strongly in $\Omega \times (0,T)$ as an equality in $H$. Multiplying it by $-\Delta \b{c}$ in $H$ yields 
	\begin{multline*}
		(\nabla \b{w}, \nabla \b{c})_H \geq \gamma_* \|\Delta \b{c}\|^2_H +  \left( \td{^2 I}{\b s^2}(\b{c}) \nabla \b{c}, \nabla\b{c}\right)_H  + \sum_{i \in \mathcal{I}}\theta_i(\Psi''(c_i)\nabla c_i, \nabla c_i)_H \\ + \sum_{(i,j) \in \{A,B\}^2} \alpha_{ij}(c_i-\overline{c_i},\, c_j)_H + \sum_{i \in \{A,B\}} \left( \pd{\epsilon}{c_i}|\b{E}_0-\nabla\Phi|^2,  -\Delta c_i\right)_H
	\end{multline*}
	thus, the convexity of $\Psi$ and the H\"{o}lder inequality yields
	\begin{equation} \label{eq:mu-ccomp}
		\gamma_*\|\Delta \b{c}\|^2_H \leq C\left( 1 + \|\nabla \b{w}\|_H \right),
	\end{equation}
	and the proof is complete recalling that $\nabla \b{w} \in L^2(0,T;H)$. In order to show the claim on the electric potential $\Phi$ we argue by comparison. Indeed, expanding computations in the Maxwell equation, we have
	\[
	-\epsilon(c_A,\,c_B)\Delta \Phi = \epsilon(c_A,\,c_B)\div \b E_0 + \td{\epsilon}{\b s}(c_A,c_B) \cdot \nabla \b E_0 + \td{\epsilon}{\b s}(c_A,c_B) \cdot \nabla \Phi.
	\]
	If $p \geq 4$ and $\alpha > \frac 65$, then the right hand side of the equality above defines, by interpolation and Assumption \ref{hyp:epsilon}, an element of $L^4(0,T;L^{\min\{\frac{6\alpha}{6+\alpha}, \frac{6p}{6+p}\}}(\Omega))$. In any case, we immediately obtain the claim by elliptic regularity.
\end{proof}
\subsection{Global solutions.} The existence of global solutions on the unbounded temporal \reviewF{interval} $[0,+\infty)$, given by Corollary \ref{cor:wellposed}, can be proved by similar means. However, the time approximation scheme has to be carried out slightly differently. For any $\tau > 0$ we analogously consider the time discretization of time step $\tau$, and still inductively construct the approximated solutions. Interpolants can be defined analogously, and uniform estimates depending on $T$ can only be carried out locally in time.
\subsection{The energy identity.} Finally, we are only left to prove the energy identity. We want to apply \cite[Lemma 4.1]{RS04}. Here, we adapt an argument given in \cite[Subsection 4.1]{GP23}. Fix $\b{m} \in \b{G}$ to be arbitrary. Let us consider a variation of the energy contribution $\mathcal{E}_1$ defined in the Introduction
\[
\mathcal{E}^*_1(\b{f}) := \ii{\Omega}{}{\sum_{i \in \mathcal{I}} \left[  \dfrac{\gamma_i}{2}|\nabla f_i|^2 + \theta_i\Psi(f_i+m_i) \right]}{x},
\]
defined on the domain
\[
\operatorname{dom}\mathcal{E}^*_1 := \left\{ \b{f} \in TX^{\b{m}} \text{ and } \Psi(f_i+m_i) \in L^1(\Omega) \:\: \forall \: i \in \{1,2,3\}  \right\}.
\]
As $\theta_i \Psi$ is still a convex function for all $i \in \mathcal{I}$, \cite[Lemma 4.1]{AW07} readily implies that $\mathcal{E}^*_1$ is the sum of three proper, convex and lower semi-continuous functionals. Therefore, we readily infer that $\mathcal{E}^*_1$ is a proper, convex and lower semi-continuous functional on its domain. Observe now that Theorem \ref{thm:wellposed} yields that the projection of the solution of the problem onto $H_{(0)}$, i.e., $\b{c}-\overline{\b{c}_0}$, clearly belongs to the domain $\operatorname{dom}\mathcal{E}^*_1$ as long as $\b{m} = \overline{\b{c}_0}$. Let us now prove the key property enabling the argument.
\begin{lem} \label{lem:subgradient}
	Assume that $p \geq 4$ and $\b{m} = \overline{\b{c}_0}$. The vector
	\[
	\b{\eta} := \b{P}\begin{bmatrix}
		\pd{I}{c_A} + \alpha_{AA}\mathcal{N}(c_A - \overline{c_A}) + \alpha_{AB}\mathcal{N}(c_B - \overline{c_B}) + \dfrac{1}{2}\pd{\epsilon}{c_A}|\b{E}_0 - \nabla \Phi|^2 \\[0.5cm]
		\pd{I}{c_B} + \alpha_{BA}\mathcal{N}(c_A - \overline{c_A}) + \alpha_{BB}\mathcal{N}(c_B - \overline{c_B}) + \dfrac{1}{2}\pd{\epsilon}{c_B}|\b{E}_0 - \nabla \Phi|^2 \\[0.3cm]
		\pd{I}{c_S}
	\end{bmatrix} 
	\]
	is such that $\b{w}(t) - \b{\eta}(t) - \overline{\b{w}(t) - \b{\eta}(t)} \in \partial\mathcal{E}^*_1(\b{c}(t)-\overline{\b{c}_0})$ for almost every $t \in [0,T]$.
\end{lem}
\begin{proof}
	We use the definition of subdifferential. The claim is equivalent to proving that
	\[
	(\b{z}-(\b{c}(t)-\overline{\b{c}_0}),\, \b{w}(t) - \b{\eta}(t) - \overline{\b{w}(t) - \b{\eta}(t)})_H \leq \mathcal{E}^*_1(\b{z}) - \mathcal{E}^*_1(\b{c}(t)-\overline{\b{c}_0})
	\]
	for almost every $t \geq 0$ and every $\b{z} \in \operatorname{dom} \mathcal{E}^*_1$. Let
	\[
	\b{\Psi}'(\b{c}) := \begin{bmatrix}
		\theta_A \Psi'(c_A) \\
		\theta_B \Psi'(c_B) \\
		\theta_S \Psi'(c_S) 
	\end{bmatrix}
	\]
	and observe that, since $\b{z}-(\b{c}(t)-\overline{\b{c}_0})$ has null mean value and takes values in $T\Sigma$ almost everywhere,	
	\[
	\begin{split} & (\b{z}-(\b{c}(t)-\overline{\b{c}_0}),\, \b{w}(t) - \b{\eta}(t) - \overline{\b{w}(t) - \b{\eta}(t)})_H \\ & \hspace{2cm} = (\b{z}-(\b{c}(t)-\overline{\b{c}_0}),\, -\Gamma\Delta \b{c} + \b{P}\b{\Psi}'(\b{c}))_H \\
		& \hspace{2cm} = (\b{z}-(\b{c}(t)-\overline{\b{c}_0}),\, -\Gamma\Delta \b{(\b{c}(t)-\overline{\b{c}_0})} + \b{P}\b{\Psi}'((\b{c}(t)-\overline{\b{c}_0}) +\overline{\b{c}_0}))_H\\
		& \hspace{2cm} =(\b{z}-(\b{c}(t)-\overline{\b{c}_0}),\, -\Gamma\Delta \b{(\b{c}(t)-\overline{\b{c}_0})} + \b{\Psi}'((\b{c}(t)-\overline{\b{c}_0}) +\overline{\b{c}_0}))_H  \\
		& \hspace{2cm} = \sum_{i \in \mathcal{I}} \gamma_i(\nabla z_i,\, \nabla c_i(t))_H - \gamma_i\|\nabla c_i(t)\|^2_H  + (z_i - (c_i(t)-\overline{c_{0i}}), \theta_i \Psi'((c_i(t)-\overline{c_{0i}}) + \overline{c_{0i}}))_H\\
		& \hspace{2cm} \leq \sum_{i \in \mathcal{I}} \dfrac{\gamma_i}{2}\|\nabla z_i \|^2_H - \dfrac{\gamma_i}{2}\|\nabla c_i (t)\|^2_H + \theta_i(\Psi(z_i+\overline{c_{0i}})-\Psi((c_i(t)-\overline{c_{0i}}) + \overline{c_{0i}}), 1)_H\\
		& \hspace{2cm} = \mathcal{E}_1^*(\b z) - \mathcal{E}_1^*(\b{c}(t)-\overline{\b{c}_0})
	\end{split}
	\]where in the last line we exploited the convexity of $\Psi$. Indeed,
	\[
	\Psi(s) - \Psi(t) \geq \Psi'(t)(s-t) 
	\]
	for all $s,\,t \in (-1,1)$. The proof is complete.
\end{proof} \noindent
In order to apply \cite[Lemma 4.1]{RS04}, we also have to check that
\[
\mathcal{E}^*_1(\b{f}) \geq C\|\b{f}\|^2_H - C,
\]
for some constant $C$ independent of $\b{f}$. This is an immediate consequence of Assumption \ref{hyp:potential}-\ref{hyp:potpositive} and the Poincaré inequality. Finally, we have
\begin{lem} 
	If $p \geq 4$ and $\b{E}_0 \in W^{1,\alpha}(\Omega;\mathbb{R}^d)$, with $\alpha \geq 4$, then the vector $\b{\eta}$ appearing in Lemma \ref{lem:subgradient} satisfies $\b{w} - \b{\eta} - \overline{\b{w} - \b{\eta}} \in L^2(0,T;V_{(0)})$.
\end{lem}
\begin{proof}
	Observe preliminarly that, by the action of $\b{P}$ given by \eqref{eq:projector}, we have 
	\[
	\|\nabla \b{P} \b{f}\|_H \leq C\|\nabla \b{f}\|_H.
	\] 
	Therefore, since $\b{P}$ is non-expansive in $H$, we can safely drop it for notational convenience. We have
	\[
	\|\b{w} - \b{\eta} - \overline{\b{w} - \b{\eta}}\|_V \leq C\left( \|\b{w}\|_V + \|\b{\eta}\|_V\right)
	\] 
	Proposition \ref{prop:limitd1} gives $\b{w} \in L^2(0,T;V)$, and, after straightforward controls, $\b{\eta} \in L^2(0,T;H)$, hence we focus on controlling the gradient of $\b{\eta}$. In particular, consider the gradient of the first component of $\b{\eta}$. The local interaction term satisfies
	\[
	\left\|\nabla \pd{I}{c_A} \right\|_H \leq C,
	\]
	directly from Assumption \ref{hyp:interaction}. The nonlocal interaction terms are easily bounded by virtue of the fact that the range of the continuous operator $\mathcal{N}$ acting on $H$ is contained in $H^2(\Omega)$ and we have
	\[
	\|\mathcal{N}(c_A - \overline{c_A})\|_{H^2(\Omega)} + \|\mathcal{N}(c_B - \overline{c_B})\|_{H^2(\Omega)} \leq C\|\b{c}\|_H \leq C.\qquad 
	\]Finally, we deal with the gradient of the electric term. Without loss of generality, we assume $p = \alpha = 4$ (hence $\beta = \frac{12}{5}$). Here we can exploit the additional regularity given by Proposition \ref{prop:higher4}
	\[ 
	\begin{split}
		& \left\| \nabla \left[  \pd{\epsilon}{c_A}|\b{E}_0 - \nabla \Phi|^2 \right] \right\|_H \\
		& \hspace{1cm} \leq \left\| \left(\pd{^2\epsilon}{c_A^2}\nabla c_A + \pd{^2\epsilon}{c_A\partial c_B}\nabla c_B \right)|\b{E}_0 - \nabla \Phi|^2 \right\|_H  + \left\| 2\pd{\epsilon}{c_A}(\nabla \b{E}_0 - D^2\Phi)(\b{E}_0 - \nabla \Phi) \right\|_H \\
		& \hspace{1cm} \leq C\bigg( \|\b{E}_0 - \nabla \Phi\|_{L^4(\Omega;\mathbb{R}^d)}\|\b{E}_0 - \nabla \Phi\|_{L^{12}(\Omega;\mathbb{R}^d)} \|\nabla \b{c}\|_{L^6(\Omega;\mathbb{R}^{3\times d})} \\
		& \hspace{6cm} + \|\nabla \b{E}_0 - D^2\Phi\|_{L^\frac{12}{5}(\Omega;\mathbb{R}^{d \times d})}\|\b{E}_0 - \nabla \Phi \|_{L^{12}(\Omega;\mathbb{R}^d)} \bigg)\\
		& \hspace{1cm} \leq C\|\b{E}_0 - \nabla \Phi\|_{L^{12}(\Omega;\mathbb{R}^d)}( \|\nabla \b{E}_0 - D^2\Phi\|_{L^\frac{12}{5}(\Omega;\mathbb{R}^{d \times d})} +  \|\nabla \b{c}\|_{L^6(\Omega;\mathbb{R}^{3\times d})}),
	\end{split}
	\]
	where we used the continuous embedding $W^{1,\frac{12}{5}}(\Omega) \hookrightarrow L^{12}(\Omega)$ in (two and) three dimensions. Collecting the results (after arguing similarly for the remaining components of $\b \eta$) and squaring the inequality we have
	\[
	\|\b{w} - \b{\eta} - \overline{\b{w} - \b{\eta}}\|_V^2 \leq C \left[ 1 + \|\nabla \b{E}_0 - D^2\Phi\|_{L^\frac{12}{5}(\Omega;\mathbb{R}^{d \times d})}^2 +\| \b{c}\|^2_{H^2(\Omega;\mathbb{R}^{3})}\left( 1 + \|\b{E}_0 - \nabla \Phi\|_{L^{12}(\Omega;\mathbb{R}^d)}^2 \right) \right].
	\] 
	Integrating in time yields the result, recalling Proposition \ref{prop:higher4}.
\end{proof} \noindent
Hence, we have that \cite[Lemma 4.1]{RS04} implies that $\mathcal{E}^*_1$ is absolutely continuous in $[0,T]$, and moreover,
\[
\int_s^t \left\langle \pd{\b{c}}{t}(\tau),\, \b{w}(\tau) - \b{\eta}(\tau) - \overline{\b{w}(\tau) - \b{\eta}(\tau)} \right\rangle_{V_{(0)}^*, V_{(0)}}\: \d s = \mathcal{E}^*_1(\b{c}(t)-\overline{\b c_0}) - \mathcal{E}^*_1(\b{c}(s)-\overline{\b c_0})
\]
for all $s,\,t \in [0,T]$. From here, it is straightforward (see \cite[Subsection 4.1]{GP23}) to deduce that also
\[
\int_s^t \left\langle \pd{\b{c}}{t}(\tau),\, \b{w}(\tau) - \b{\eta}(\tau) \right\rangle_{V^*, V}\: \d s = \mathcal{E}_1(\b{c}(t)) - \mathcal{E}_1(\b{c}(s)),
\]
where the functional $\mathcal{E}_1$ is one of the energy contributions defined in the Introduction. Of course, it still holds that $\mathcal{E}_1$ is absolutely continuous in $[0,T]$. Moreover, since $V \hookrightarrow H$ and $\b c \in C^0([0,T];H) \cap L^\infty(0,T;V)$, then it holds that $\b c$ is weakly continuous in time with values in $V$, namely $\b c \in C^0_{\text{w}}([0,T];V)$ (refer, for instance, to \cite[Chapter II, Lemma 3.3]{RTDDS}). By continuity of the energy functional, it is readily shown also that the map
\[
t \mapsto \|\nabla \b c(t)\|_{H}
\]
is continuous in $\mathbb{R}$. The two aforementioned properties are enough to infer that $\b{c} \in C([0,T];V)$. It is then possible to exploit the weak formulation of the problem to reconstruct the time derivative of the energy functional by testing the evolution equation for $\b c$ by $\b w(t)$ to show
\[
\td{}{t}\mathcal{E}(t) + \int_\Omega \b M(\b c(t)) \nabla \b w(t) \cdot \nabla \b w (t) \: \d x = 0,
\]
for almost all $t \in [0,T]$. Integrating the above in the interval $[0,t]$ yields the claim.
	\section{Proof of Theorem \ref{thm:uniqueness}}
	\label{sec:uniqueness}
	Let us now prove Theorem \ref{thm:uniqueness}, i.e., the uniqueness result. Before giving the details of the argument, let us introduce some convenient notation. For $k \in \{1,2\}$, let $\b{c}_{0k}$ denote a function complying with the hypotheses given in Assumption \ref{hyp:data}-\ref{hyp:initial}. Let $(\b{c}_k,\,\b{w}_k, \Phi_k)$ denote the corresponding solutions of \eqref{eq:strongform2}-\eqref{eq:conditions2} satisfying
	\[
	\b{c}_k(0) = \b{c}_{0k}, \qquad k \in \{1,2\}.	
	\]
	Consider the differences
	\[
	\b{c} := \b{c}_1 - \b{c}_2, \qquad \b{w} := \b{w}_1 - \b{w}_2, \qquad 
	\Phi := \Phi_1 - \Phi_2, \qquad \b{c}_0 := \b{c}_{01} - \b{c}_{02},
	\]
	satisfying the following equations:
	\begin{equation} \label{eq:uniqueness0}
		\begin{cases}
			\left\langle \pd{\b{c}}{t}, \: \b{v} \right\rangle_{V^*, V} + \left( M(\b{c}_1)\nabla \b{w}_1 -  M(\b{c}_2)\nabla \b{w}_2, \nabla \b{v} \right)_H = 0, & \quad \forall \: \b{v} \in V, \\[0.5cm]
			\left(\b{w}, \b{z} \right)_H = \left( \Gamma \nabla \b{c} , \nabla \b{z} \right)_H + \left( \b P( \widetilde{\b{\mu}}_1 - \widetilde{\b{\mu}}_2),\, \b{z} \right)_H, & \quad \forall \: \b{z} \in V \cap L^\infty(\Omega;\mathbb{R}^3), \\[0.5cm]
			(\epsilon({c}_{1A},\,{c}_{1B})(\b{E}_0 - \nabla \Phi_1) - \epsilon({c}_{2A},\,{c}_{2B})(\b{E}_0 - \nabla \Phi_2) , \nabla \Lambda)_H = 0 & \quad \forall \: \Lambda \in V_0,
		\end{cases}
	\end{equation}
	with $\b{c}(0) = \b{c}_0$ and
	\[
	\widetilde{\b{\mu}}_k := \begin{bmatrix}
		\theta_A\Psi'(c_{kA}) + \pd{I}{c_{kA}} + \alpha_{AA}\mathcal{N}(c_{kA} - \overline{c_{kA}}) + \alpha_{AB}\mathcal{N}(c_{kB} - \overline{c_{kB}}) + \dfrac{1}{2}\pd{\epsilon}{c_{kA}}|\b{E}_0 - \nabla \Phi|^2 \\[0.5cm]
		\theta_B\Psi'(c_{kB}) + \pd{I}{c_{kB}} + \alpha_{BA}\mathcal{N}(c_{kA} - \overline{c_{kA}}) + \alpha_{BB}\mathcal{N}(c_{kB} - \overline{c_{kB}}) + \dfrac{1}{2}\pd{\epsilon}{c_{kB}}|\b{E}_0 - \nabla \Phi|^2 \\[0.3cm]
		\theta_S\Psi'(c_{kS}) + \pd{I}{c_{kS}}
	\end{bmatrix}.
	\]
	Following an idea in \cite{Barrett99}, it is convenient to recast the first equation in an alternative form. It holds
	\[
	\left\langle \pd{\b{c}_k}{t}, \: \b{v} \right\rangle_{V^*, V} = -\left\langle A_{\b{c}_k} \b{w}_k, \: \b{v} \right\rangle_{V^*, V}
	\]
	choosing $F = M$ in \eqref{eq:generalA}. Therefore, we also have
	\[
	-\mathcal{N}_{\b{c}_k} \pd{\b{c}_k}{t} = \b{w}_k - \overline{\b{w}_k},
	\]
	and the above is an equality in $V$. Let us then write \eqref{eq:uniqueness0} as
	\begin{equation} \label{eq:uniqueness1}
		\begin{cases}
			\b{w} - \overline{\b{w}} = -\mathcal{N}_{\b{c}_1}\pd{\b{c}_1}{t} + \mathcal{N}_{\b{c}_2}\pd{\b{c}_2}{t}, \\
			\left(\b{w}, \b{z} \right)_H = \left( \Gamma \nabla \b{c} , \nabla \b{z} \right)_H + \left( \b P( \widetilde{\b{\mu}}_1 - \widetilde{\b{\mu}}_2),\, \b{z} \right)_H, & \quad \forall \: \b{z} \in V \cap L^\infty(\Omega;\mathbb{R}^3),\\
			(\epsilon({c}_{1A},\,{c}_{1B})(\b{E}_0 - \nabla \Phi_1) - \epsilon({c}_{2A},\,{c}_{2B})(\b{E}_0 - \nabla \Phi_2) , \nabla \Lambda)_H = 0 & \quad \forall \: \Lambda \in V_0.
		\end{cases}
	\end{equation}
	Firstly, we address the general case. Then, we investigate how the results can be improved in the constant mobility case. In both scenarios, without loss of generality, we set $r = 4$.
	\subsection{The general case.} \label{ssec:general} Let us multiply the first equation in \eqref{eq:uniqueness1} by $\b{c} - \overline{\b{c}} = \b{c} - \overline{\b{c}_0}$ in $H$. Taking into account the second equation in \eqref{eq:uniqueness1}, this results in
	\begin{multline} \label{eq:uniqueG1}
		(\Gamma \nabla \b{c}, \nabla \b{c})_H + \sum_{i \in \mathcal{I}} \theta_i (\Psi'(c_{1i}) - \Psi'(c_{2i}), c_i - \overline{c_{0i}})_H + \left(\td{I}{\b s} (\b{c}_1) - \td{I}{\b s} (\b{c}_2), \b{c} - \overline{\b{c}_0}\right) _H \\
		+ \sum_{(i,j) \in \{A,B\}^2} \alpha_{ij} (\mathcal{N}(c_i-\overline{c_i}), c_j-\overline{c_{0j}})_H + \dfrac{1}{2} \sum_{i \in \{A,B\}} \left(\pd{\epsilon}{c_{1i}}|\b{E}_0 - \nabla \Phi_1|^2 -\pd{\epsilon}{c_{2i}}|\b{E}_0 - \nabla \Phi_2|^2, c_i - \overline{c_{0i}} \right)_H \\
		= -\left( \mathcal{N}_{\b{c}_1} \pd{\b{c}}{t}, \b{c} - \overline{\b{c}_0} \right)_H -\left( [\mathcal{N}_{\b{c}_1} - \mathcal{N}_{\b{c}_2}] \pd{\b{c_2}}{t}, \b{c} - \overline{\b{c}_0} \right)_H.
	\end{multline}
	We now address each of the terms appearing in \eqref{eq:uniqueG1}. A basic estimate shows that
	\begin{equation} \label{eq:uniqueG2}
		(\Gamma \nabla \b{c}, \nabla \b{c})_H \geq \gamma_*\|\nabla \b{c}\|^2_H,
	\end{equation}
	whereas, owing to the convexity of $\Psi$ and the fact that $\b{c} \in \b{G}$,
	\begin{equation} \label{eq:uniqueG3}
		\sum_{i \in \mathcal{I}} \theta_i (\Psi'(c_{1i}) - \Psi'(c_{2i}), c_i)_H \geq 0.
	\end{equation}
	The remaining terms are controlled by
	\begin{equation} \label{eq:uniqueG3bis}
		\left| \sum_{i \in \mathcal{I}} \theta_i (\Psi'(c_{1i}) - \Psi'(c_{2i}),\, \overline{c_{0i}})_H \right| \leq \theta^*\sum_{i \in \mathcal{I}}|\overline{c_{0i}}|\left[\|\Psi'(c_{1i})\|_{L^1(\Omega)} + \|\Psi'(c_{2i})\|_{L^1(\Omega)} \right].
	\end{equation}
	Moreover, we have by the mean value theorem
	\begin{equation} \label{eq:uniqueG4}
		\left|\left(\td{I}{\b s} (\b{c}_1) - \td{I}{\b s} (\b{c}_2), \b{c} - \overline{\b{c}_0}\right) _H\right| \leq \|I\|_{C^2([0,1]^3)}\|\b{c}\|_H\|\b{c} - \overline{\b{c}_0}\|_H.
	\end{equation}
	The nonlocal terms to the left hand side are controlled by
	\begin{equation} \label{eq:uniqueG5}
		\left| \sum_{(i,j) \in \{A,B\}^2} \alpha_{ij} (\mathcal{N}(c_i-\overline{c_i}), c_j-\overline{c_{0j}})_H \right| \leq C\|\b{c} - \overline{\b{c}_0}\|_H^2
	\end{equation}
	by continuity of $\mathcal{N}$. Here, the constant $C$ depends on $\alpha_{ij}$. Let us now address the terms linked to the electric quantities. Observe that,
	\begin{equation} \label{eq:uniqueG6}
		\begin{split}
			& \dfrac{1}{2} \sum_{i \in \{A,B\}} \left(\pd{\epsilon}{c_{1i}}|\b{E}_0 - \nabla \Phi_1|^2 -\pd{\epsilon}{c_{2i}}|\b{E}_0 - \nabla \Phi_2|^2, c_i - \overline{c_{0i}} \right)_H \\
			& \hspace{2cm} = \dfrac{1}{2} \sum_{i \in \{A,B\}} \left(\pd{\epsilon}{c_{1i}}\left[|\b{E}_0 - \nabla \Phi_1|^2 -|\b{E}_0 - \nabla \Phi_2|^2 \right], c_i - \overline{c_{0i}} \right)_H  \\
			& \hspace{4cm}+ \dfrac{1}{2} \sum_{i \in \{A,B\}} \left(\left[\pd{\epsilon}{c_{1i}} -\pd{\epsilon}{c_{2i}} \right]|\b{E}_0 - \nabla \Phi_2|^2, c_i - \overline{c_{0i}} \right)_H.
		\end{split}
	\end{equation}
	Here we exploit the additional properties given by Assumption \ref{hyp:additional}. The first term can be handled observing that
	\[
	|\b{E}_0 - \nabla \Phi_1|^2 -|\b{E}_0 - \nabla \Phi_2|^2 = (2\b{E}_0 - \nabla \Phi_1 - \nabla \Phi_2) \cdot (- \nabla \Phi).
	\]
	Indeed,
	\begin{equation} \label{eq:uniqueG7}
		\begin{split}
			& \left| \dfrac{1}{2} \sum_{i \in \{A,B\}} \left(\pd{\epsilon}{c_{1i}}\left[|\b{E}_0 - \nabla \Phi_1|^2 -|\b{E}_0 - \nabla \Phi_2|^2 \right], c_i - \overline{c_{0i}} \right)_H \right| \\
			& \hspace{3cm} =  \left| \dfrac{1}{2} \sum_{i \in \{A,B\}} \left(\pd{\epsilon}{c_{1i}}(2\b{E}_0 - \nabla \Phi_1 - \nabla \Phi_2) \cdot (- \nabla \Phi), c_i - \overline{c_{0i}} \right)_H \right| \\
			& \hspace{3cm} \leq C\|2\b{E}_0 - \nabla \Phi_1 - \nabla \Phi_2\|_{L^4(\Omega;\mathbb{R}^d)}\|\nabla \Phi\|_{H}\|\b{c}-\overline{\b{c}_0}\|_{L^4(\Omega;\mathbb{R}^3)} \\
			& \hspace{3cm} \leq C\|\nabla \Phi\|_{H}\|\b{c}-\overline{\b{c}_0}\|_{L^4(\Omega;\mathbb{R}^3)},
		\end{split}
	\end{equation}
	whereas by the H\"{o}lder inequality
	\begin{equation} \label{eq:uniqueG8}
		\begin{split}
			& \left|\dfrac{1}{2} \sum_{i \in \{A,B\}} \left(\left[\pd{\epsilon}{c_{1i}} -\pd{\epsilon}{c_{2i}} \right]|\b{E}_0 - \nabla \Phi_2|^2, c_i - \overline{c_{0i}} \right)_H\right| \\
			& \hspace{3cm}\leq C\|\b{c}\|_{L^4(\Omega;\mathbb{R}^3)}\|\b{c}-\overline{\b{c}_0}\|_{L^4(\Omega;\mathbb{R}^3)}\|\b{E}_0 - \nabla \Phi_2\|^2_{L^4(\Omega;\mathbb{R}^d)} \\
			& \hspace{3cm} \leq C\|\b{c}\|_{L^4(\Omega;\mathbb{R}^3)}\|\b{c}-\overline{\b{c}_0}\|_{L^4(\Omega;\mathbb{R}^3)}.
		\end{split}
	\end{equation}
	Collecting estimates \eqref{eq:uniqueG2}-\eqref{eq:uniqueG8} in \eqref{eq:uniqueG1}, we have
	\begin{multline} \label{eq:uniqueG9}
		\gamma_*\|\nabla \b{c}\|^2_H + \left( \mathcal{N}_{\b{c}_1} \pd{\b{c}}{t}, \b{c} - \overline{\b{c}_0} \right)_H \leq \theta^*\sum_{i \in \mathcal{I}}|\overline{c_{0i}}|\left[\|\Psi'(c_{1i})\|_{L^1(\Omega)} + \|\Psi'(c_{2i})\|_{L^1(\Omega)} \right] \\
		-\left( [\mathcal{N}_{\b{c}_1} - \mathcal{N}_{\b{c}_2}] \pd{\b{c_2}}{t}, \b{c} - \overline{\b{c}_0}\right)_H \\ + C\left[ \|\b{c}\|_{L^4(\Omega;\mathbb{R}^3)}\|\b{c}-\overline{\b{c}_0}\|_{L^4(\Omega;\mathbb{R}^3)} + \|\b{c} - \overline{\b{c}_0}\|_H^2 + \|\nabla \Phi\|_{H}\|\b{c}-\overline{\b{c}_0}\|_{L^4(\Omega;\mathbb{R}^3)} \right].
	\end{multline}
	Now, we deal with the two remaining scalar products in \eqref{eq:uniqueG9}. Referring to \cite[(2.21)-(2.22)]{Barrett99}, thanks on the additional regularity assumptions on the time derivatives, we have
	\[
	\dfrac 12 \td{}{t}(\mathcal{N}_{\b{c}_1}(\b{c}-\overline{\b{c}_0}), \b{c}-\overline{\b{c}_0})_H = \left(\mathcal{N}_{\b{c}_1} \pd{\b{c}}{t},\, \b{c}-\overline{\b{c}_0}\right)_H - \dfrac 12 \left(\left[ \nabla M(\b{c}_1) \pd{\b{c}_1}{t} \right] \nabla\mathcal{N}_{\b{c}_1}(\b{c}-\overline{\b{c}_0}),\, \nabla\mathcal{N}_{\b{c}_1}(\b{c}-\overline{\b{c}_0})\right)_H.
	\]
	Moreover, by the Ladyzhenskaya and H\"{o}lder inequalities we also have the controls
	\begin{equation} \label{eq:uniqueG10}
		\begin{split}
			& \left| \dfrac 12 \left(\left[ \nabla M(\b{c}_1) \pd{\b{c}_1}{t} \right] \nabla\mathcal{N}_{\b{c}_1}(\b{c}-\overline{\b{c}_0}),\, \nabla\mathcal{N}_{\b{c}_1}(\b{c}-\overline{\b{c}_0})\right)_H \right| \\
			& \hspace{3cm} \leq C\left\|\pd{\b{c}_1}{t}\right\|_{H}\|\nabla\mathcal{N}_{\b{c}_1}(\b{c}-\overline{\b{c}_0})\|_{L^4(\Omega;\mathbb{R}^{3 \times d})}^2 \\
			& \hspace{3cm} \leq C\left\|\pd{\b{c}_1}{t}\right\|_{H}\|\nabla \mathcal{N}_{\b{c}_1}(\b{c}-\overline{\b{c}_0})\|_{H}^{2-\frac d2}\|\mathcal{N}_{\b{c}_1}(\b{c}-\overline{\b{c}_0})\|_{H^2(\Omega;\mathbb{R}^3)}^{\frac d2},
		\end{split} 
	\end{equation}
	and, thanks to Assumption \ref{hyp:additional}, by Lemma \ref{lem:controlA} we have,
	\begin{equation} \label{eq:uniqueG11}
		\begin{split}
			\left| \left( [\mathcal{N}_{\b{c}_1} - \mathcal{N}_{\b{c}_2}] \pd{\b{c_2}}{t}, \b{c} - \overline{\b{c}_0}\right)_H \right| & = \left| \left( \nabla[\mathcal{N}_{\b{c}_1} - \mathcal{N}_{\b{c}_2}] \pd{\b{c_2}}{t}, M(\b{c}_1)\nabla\mathcal{N}_{\b{c}_1}(\b{c} - \overline{\b{c}_0})\right)_H \right|  \\
			& \leq C\|\nabla \mathcal{N}_{\b{c}_1}(\b{c} - \overline{\b{c}_0})\|_H\|\b{c}\|_{L^4(\Omega;\mathbb{R}^3)}\left\|\nabla\mathcal{N}_{\b{c}_2}\pd{\b{c_2}}{t} \right\|_{L^4(\Omega;\mathbb{R}^3)} \\
			& \leq C\|\nabla \mathcal{N}_{\b{c}_1}(\b{c} - \overline{\b{c}_0})\|_H\|\b{c}\|^{1-\frac{d}{4}}_{H}\|\b{c}\|^\frac{d}{4}_{V}\|\nabla\b{w}_2\|_H^{1-\frac{d}{4}}\left\|\mathcal{N}_{\b{c}_2}\pd{\b{c_2}}{t}\right\|_{H^2(\Omega;\mathbb{R}^3)}^{\frac{d}{4}},
		\end{split}
	\end{equation}
	where we also used the Ladyzhenskaya inequality. Now, using \eqref{eq:uniqueG10}-\eqref{eq:uniqueG11} in \eqref{eq:uniqueG9} yields
	\begin{multline} \label{eq:uniqueG12}
		\gamma_*\|\nabla \b{c}\|^2_H + \dfrac 12 \td{}{t}(\mathcal{N}_{\b{c}_1}(\b{c}-\overline{\b{c}_0}), \b{c}-\overline{\b{c}_0})_H \leq \theta^*\max_{i \in \mathcal{I}}|\overline{c_{0i}}|\sum_{i \in \mathcal{I}}\left[\|\Psi'(c_{1i})\|_{L^1(\Omega)} + \|\Psi'(c_{2i})\|_{L^1(\Omega)} \right] \\ + C\bigg[ \|\b{c}\|_{L^4(\Omega;\mathbb{R}^3)}\|\b{c}-\overline{\b{c}_0}\|_{L^4(\Omega;\mathbb{R}^3)} + \|\b{c} - \overline{\b{c}_0}\|_H^2 + \|\nabla \Phi\|_{H}\|\b{c}-\overline{\b{c}_0}\|_{L^4(\Omega;\mathbb{R}^3)}  \\
		+ \left\|\pd{\b{c}_1}{t}\right\|_{H}\|\nabla \mathcal{N}_{\b{c}_1}(\b{c}-\overline{\b{c}_0})\|_{H}^{2-\frac d2}\|\mathcal{N}_{\b{c}_1}(\b{c}-\overline{\b{c}_0})\|_{H^2(\Omega;\mathbb{R}^3)}^{\frac d2} \\ + \|\nabla \mathcal{N}_{\b{c}_1}(\b{c} - \overline{\b{c}_0})\|_H\|\b{c}\|^{1-\frac{d}{4}}_{H}\|\b{c}\|^\frac{d}{4}_{V}\|\nabla\b{w}_2\|_H^{1-\frac{d}{4}}\left\|\mathcal{N}_{\b{c}_2}\pd{\b{c_2}}{t}\right\|_{H^2(\Omega;\mathbb{R}^3)}^{\frac{d}{4}}\bigg].
	\end{multline}
	Before tackling the norm products to the right hand side of \eqref{eq:uniqueG12}, we deal with the last equation of \eqref{eq:uniqueness1}. Choosing $\Lambda = -\Phi$, we have
	\[
	([\epsilon({c}_{1A},\,{c}_{1B})- \epsilon({c}_{2A},\,{c}_{2B})](\b{E}_0 - \nabla \Phi_1),\, -\nabla \Phi)_H + ( \epsilon({c}_{2A},\,{c}_{2B})\nabla \Phi ,\, \nabla \Phi)_H = 0.
	\]
	Therefore, we reach
	\[
	\epsilon_*\|\nabla \Phi\|^2_H \leq C\|\nabla \Phi\|_H\|\b{E}_0 - \nabla \Phi_1\|_{L^4(\Omega;\mathbb{R}^d)}\|\b{c}\|_{L^4(\Omega;\mathbb{R}^3)},
	\]
	implying
	\begin{equation} \label{eq:uniqueG13}
		\epsilon_*\|\nabla \Phi\|_H^2 \leq C\|\b{c}\|_{L^4(\Omega;\mathbb{R}^3)}^2.
	\end{equation}
	Summing \eqref{eq:uniqueG13} to \eqref{eq:uniqueG12} gives
	\begin{multline} \label{eq:uniqueG14}
		\gamma_*\|\nabla \b{c}\|^2_H + \epsilon_*\|\nabla \Phi\|^2_H+ \dfrac 12 \td{}{t}\|\b{c}-\overline{\b{c}_0}\|^2_{*,\b{c}_1} \leq \theta^*\max_{i \in \mathcal{I}}|\overline{c_{0i}}|\sum_{i \in \mathcal{I}}\left[\|\Psi'(c_{1i})\|_{L^1(\Omega)} + \|\Psi'(c_{2i})\|_{L^1(\Omega)} \right] \\ + C\bigg[ \|\b{c}\|_{L^4(\Omega;\mathbb{R}^3)}\|\b{c}-\overline{\b{c}_0}\|_{L^4(\Omega;\mathbb{R}^3)} + \|\b{c} - \overline{\b{c}_0}\|_H^2 + \|\b{c}\|^2_{L^4(\Omega;\mathbb{R}^3)} +  \|\nabla \Phi\|_{H}\|\b{c}-\overline{\b{c}_0}\|_{L^4(\Omega;\mathbb{R}^3)}  \\
		+ \left\|\pd{\b{c}_1}{t}\right\|_{H}\|\nabla \mathcal{N}_{\b{c}_1}(\b{c}-\overline{\b{c}_0})\|_{H}^{2-\frac d2}\|\mathcal{N}_{\b{c}_1}(\b{c}-\overline{\b{c}_0})\|_{H^2(\Omega;\mathbb{R}^3)}^{\frac d2} \\+ \|\nabla \mathcal{N}_{\b{c}_1}(\b{c} - \overline{\b{c}_0})\|_H\|\b{c}\|^{1-\frac{d}{4}}_{H}\|\b{c}\|^\frac{d}{4}_{V}\|\nabla\b{w}_2\|_H^{1-\frac{d}{4}}\left\|\mathcal{N}_{\b{c}_2}\pd{\b{c_2}}{t}\right\|_{H^2(\Omega;\mathbb{R}^3)}^{\frac{d}{4}} \bigg].
	\end{multline}
	Finally, we address the terms to the right hand side. To this end, we set
	\begin{align*}
		\mathcal{I}_1 & := \|\b{c}\|_{L^4(\Omega;\mathbb{R}^3)}\|\b{c}-\overline{\b{c}_0}\|_{L^4(\Omega;\mathbb{R}^3)} + \|\b{c} - \overline{\b{c}_0}\|_H^2 + \|\b{c}\|^2_{L^4(\Omega;\mathbb{R}^3)}, \\
		\mathcal{I}_2 & := \|\nabla \Phi\|_{H}\|\b{c}-\overline{\b{c}_0}\|_{L^4(\Omega;\mathbb{R}^3)}, \\
		\mathcal{I}_3 & := \left\|\pd{\b{c}_1}{t}\right\|_{H}\|\nabla \mathcal{N}_{\b{c}_1}(\b{c}-\overline{\b{c}_0})\|_{H}^{2-\frac d2}\|\mathcal{N}_{\b{c}_1}(\b{c}-\overline{\b{c}_0})\|_{H^2(\Omega;\mathbb{R}^3)}^{\frac d2}, \\
		\mathcal{I}_4 & := \|\nabla \mathcal{N}_{\b{c}_1}(\b{c} - \overline{\b{c}_0})\|_H\|\b{c}\|^{1-\frac{d}{4}}_{H}\|\b{c}\|^\frac{d}{4}_{V}\|\nabla\b{w}_2\|_H^{1-\frac{d}{4}}\left\|\mathcal{N}_{\b{c}_2}\pd{\b{c_2}}{t}\right\|_{H^2(\Omega;\mathbb{R}^3)}^{\frac{d}{4}}.
	\end{align*}
	The first two controls are rather straightforward. Indeed, we have, by repeated applications of the Ladyzhenskaya and Young inequalities
	\begin{equation} \label{eq:uniqueG15}
		\mathcal{I}_1 \leq \delta_1\|\nabla \b{c}\|^2_H + C\left(\|\b{c}-\overline{\b{c}_0}\|^2_{*,\b{c}_1} + |\overline{c_0}|^2\right),
	\end{equation}
	for any $\delta_1 > 0$. The constant $C$ depends on $\delta_1$. Moreover, a similar argument gives
	\begin{equation} \label{eq:uniqueG16}
		\mathcal{I}_2 \leq \delta_1\|\nabla \b{c}\|^2_H + \delta_2\|\nabla \Phi\|^2 +  C\|\b{c}-\overline{\b{c}_0}\|^2_{*,\b{c}_1}
	\end{equation}
	for any $\delta_2 > 0$. The constant $C$ depends on $\delta_1$ and $\delta_2$. In order to give controls for the last two terms, let us recall the inequality given in \cite[(2.26)]{Barrett99}
	\begin{equation} \label{eq:barrett}
		\left\|\mathcal{N}_{\b{c}_k}\b{f}\right\|_{H^2(\Omega;\mathbb{R}^3)}^{\alpha} \leq C\left(\|\nabla \b{c}_k\|_H^{\alpha \beta}\|\b{c}_k\|^{\alpha2^{d-2}}_{H^2(\Omega;\mathbb{R}^3)}\|\nabla \mathcal{N}_{\b{c}_k}\b{f}\|_H^\alpha + \|\b{f}\|_H^\alpha \right),
	\end{equation}
	for all $\b{f} \in H \cap Y$, $k \in \{1,2\}$ and $\alpha > 0$. Here, we have
	\[
	\beta := \begin{cases}
		\dfrac{b}{1-b} & \quad \text{if } d = 2, \text{ for any } b \in (0,1),\\
		0 & \quad \text{if } d = 3.
	\end{cases}
	\] The constant $C$ depends on $M$ and on $\alpha$. First, we have, using \eqref{eq:barrett} with $\alpha = \frac{d}{2}$ and $k = 1$,
	\begin{equation} \label{eq:uniqueG17}
		\begin{split}
			\mathcal{I}_3 & \leq C\left\|\pd{\b{c}_1}{t}\right\|_{H}\|\nabla \mathcal{N}_{\b{c}_1}(\b{c}-\overline{\b{c}_0})\|_{H}^{2-\frac d2}\left(\|\b{c}_1\|^{\frac{d2^{d}}{8}}_{H^2(\Omega;\mathbb{R}^3)}\|\nabla \mathcal{N}_{\b{c}_1}(\b{c}-\overline{\b{c}_0})\|_H^\frac d2 + \|\b{c}-\overline{\b{c}_0}\|_H^\frac d2 \right) \\ 
			& \leq C\|\nabla \mathcal{N}_{\b{c}_1}(\b{c}-\overline{\b{c}_0})\|_{H}^2\left( \left\|\pd{\b{c}_1}{t}\right\|_{H}\|\b{c}_1\|^{\frac{d2^{d}}{8}}_{H^2(\Omega;\mathbb{R}^3)} \right) + \|\nabla \mathcal{N}_{\b{c}_1}(\b{c}-\overline{\b{c}_0})\|_{H}^{2-\frac d2}\left(\left\|\pd{\b{c}_1}{t}\right\|_{H}\|\b{c}-\overline{\b{c}_0}\|_H^\frac d2 \right) \\
			& \leq C\|\nabla \mathcal{N}_{\b{c}_1}(\b{c}-\overline{\b{c}_0})\|_{H}^2\left( \left\|\pd{\b{c}_1}{t}\right\|_{H}^\frac{8}{8-d}+ \|\b{c}_1\|^{2^{d}}_{H^2(\Omega;\mathbb{R}^3)} \right) + \|\nabla \mathcal{N}_{\b{c}_1}(\b{c}-\overline{\b{c}_0})\|_{H}^{2-\frac d4}\left(\left\|\pd{\b{c}_1}{t}\right\|_{H}\|\nabla \b{c}\|_H^\frac d4 \right) \\
			& \leq \delta_1\|\nabla \b{c}\|^2_{H} + C\|\nabla \mathcal{N}_{\b{c}_1}(\b{c}-\overline{\b{c}_0})\|_{H}^2\left( \left\|\pd{\b{c}_1}{t}\right\|_{H}^\frac{8}{8-d}+ \|\b{c}_1\|^{2^{d}}_{H^2(\Omega;\mathbb{R}^3)} \right).
		\end{split}
	\end{equation}
	Next, we have, still using \eqref{eq:barrett} with $\alpha = \frac d4$ and $k = 2$,  
	\begin{equation} \label{eq:uniqueG18} \small
		\begin{split}
			\mathcal{I}_4 & \leq C\|\nabla \mathcal{N}_{\b{c}_1}(\b{c} - \overline{\b{c}_0})\|_H\|\b{c}\|^{1-\frac{d}{4}}_{H}\|\nabla \b{c}\|^\frac{d}{4}_{H}\|\nabla\b{w}_2\|_H^{1-\frac{d}{4}}\left(\|\b{c}_2\|^{\frac{d2^{d}}{16}}_{H^2(\Omega;\mathbb{R}^3)}\left\|\nabla \mathcal{N}_{\b{c}_2}\pd{\b{c}_2}{t}\right\|_H^\frac d4 + \left\|\pd{\b{c}_2}{t}\right\|_H^\frac d4 \right) \\ 
			& \leq  C\left[ \|\b{c}-\overline{\b{c}_0}\|^{\frac{12-d}{8}}_{*,\b{c}_1} + |\overline{\b{c}_0}|^{\frac{12-d}{8}}\ \right] \|\nabla \b{c}\|^{\frac{4+d}{8}}_{H}\left(\|\nabla\b{w}_2\|_H\|\b{c}_2\|^{\frac{d2^{d}}{16}}_{H^2(\Omega;\mathbb{R}^3)} + \left\|\pd{\b{c}_2}{t}\right\|_H^\frac d4 \right) \\
			& \leq \delta_1\|\nabla \b{c}\|^2_H + C\left[ \|\b{c}-\overline{\b{c}_0}\|^2_{*,\b{c}_1} + |\overline{\b{c}_0}|^2 \right]\left( \|\nabla\b{w}_2\|_H^{\frac{16}{12-d}}\|\b{c}_2\|^{\frac{d2^{d}}{12-d}}_{H^2(\Omega;\mathbb{R}^3)} + \left\|\pd{\b{c}_2}{t}\right\|_H^\frac{4d}{12-d} \right) \\
			& \leq \delta_1\|\nabla \b{c}\|^2_H + C\left[ \|\b{c}-\overline{\b{c}_0}\|^2_{*,\b{c}_1} + |\overline{\b{c}_0}|^2 \right]\left( 1+  \|\nabla\b{w}_2\|_H^{\frac{8}{6-d}} + \|\b{c}_2\|^{2^{d}}_{H^2(\Omega;\mathbb{R}^3)} + \left\|\pd{\b{c}_2}{t}\right\|_H^\frac{8}{8-d} \right), \\
		\end{split}
	\end{equation}
	where we observe that 
	\[
	\dfrac{4d}{12-d} < \dfrac{8}{8-d}, \qquad d \in \{2,3\}.
	\]
	Plugging then \eqref{eq:uniqueG15}-\eqref{eq:uniqueG16} and \eqref{eq:uniqueG17}-\eqref{eq:uniqueG18} into \eqref{eq:uniqueG14} gives
	\begin{multline} \label{eq:uniqueG19} \small
		\gamma_*\|\nabla \b{c}\|^2_H + \epsilon_*\|\nabla \Phi\|^2_H+ \dfrac 12 \td{}{t}\|\b{c}-\overline{\b{c}_0}\|^2_{*,\b{c}_1} \leq 4\delta_1\|\nabla \b{c}\|^2_H + \delta_2\|\nabla \Phi\|^2_H \\ \small + \theta^*\max_{i \in \mathcal{I}}|\overline{c_{0i}}|\sum_{i \in \mathcal{I}}\left[\|\Psi'(c_{1i})\|_{L^1(\Omega)} + \|\Psi'(c_{2i})\|_{L^1(\Omega)} \right] \\ \small
		+ C\left[ \|\b{c}-\overline{\b{c}_0}\|^2_{*,\b{c}_1} + |\overline{\b{c}_0}|^2 \right]\left( 1+  \|\nabla\b{w}_2\|_H^{\frac{8}{6-d}} + \|\b{c}_1\|^{2^{d}}_{H^2(\Omega;\mathbb{R}^3)}+ \|\b{c}_2\|^{2^{d}}_{H^2(\Omega;\mathbb{R}^3)} +\left\|\pd{\b{c}_1}{t}\right\|_{H}^\frac{8}{8-d}+ \left\|\pd{\b{c}_2}{t}\right\|_H^\frac{8}{8-d} \right).
	\end{multline}
	Choose now $\delta_1 = \frac{\gamma_*}{8}$ and $\delta_2 = \frac{\epsilon_*}{2}$ in \eqref{eq:uniqueG19}, and define the maps
	\begin{align*}
		\mathcal{G}(t) & := C\left( 1+  \|\nabla\b{w}_2\|_H^{\frac{8}{6-d}} + \|\b{c}_1\|^{2^{d}}_{H^2(\Omega;\mathbb{R}^3)}+ \|\b{c}_2\|^{2^{d}}_{H^2(\Omega;\mathbb{R}^3)} +\left\|\pd{\b{c}_1}{t}\right\|_{H}^\frac{8}{8-d}+ \left\|\pd{\b{c}_2}{t}\right\|_H^\frac{8}{8-d} \right), \\
		\mathcal{W}(t) & := \gamma_*\|\nabla \b{c}\|^2_H + \epsilon_*\|\nabla \Phi\|^2_H, \\
		\mathcal{Y}(t) & := \|\b{c}-\overline{\b{c}_0}\|^2_{*,\b{c}_1} + |\overline{\b{c}_0}|^2, \\
		\mathcal{Z}(t) & := \theta^*\sum_{i \in \mathcal{I}}\left[\|\Psi'(c_{1i})\|_{L^1(\Omega)} + \|\Psi'(c_{2i})\|_{L^1(\Omega)} \right],
	\end{align*}
	where $C$ is the constant appearing in \eqref{eq:uniqueG19}. The very same inequality now reads
	\[
	\td{\mathcal{Y}}{t} + \mathcal{W}(t) \leq  \mathcal{G}(t)\mathcal{Y}(t) + \max_{i \in \mathcal{I}}|\overline{c_{0i}}|\mathcal{Z}(t).
	\]
	Since $\mathcal{G}$ and $\mathcal{Z}$ belong to $L^1(0,T)$, the Gronwall lemma implies that
	\[
	\mathcal{Y}(t) + \int_0^t \mathcal{W}(s) \: \d s \leq C\left[ \mathcal{Y}(0) + \max_{i \in \mathcal{I}}|\overline{c_{0i}}|\int_0^t \mathcal{Z}(s) \: \d s \right].
	\]
	Hence, we have, by definition,
	\[
	\|\b{c}_1 - \b{c}_2\|_{L^\infty(0,T;V^*)}^2 + \|\b{c}_1 - \b{c}_2\|_{L^2(0,T;V)}^2+ \|\Phi_1 - \Phi_2 \|^2_{L^2(0,T;V_0)} \leq C\left[ \|\b{c}_{01} - \b{c}_{02}\|_{V^*}^2 + \|\overline{\b{c}_{01}} - \overline{\b{c}_{02}}\|_{\mathbb{R}^3}\right],
	\]
	and uniqueness follows standardly.
	\subsection{The constant mobility case.} \label{ssec:constmob} Here, we discuss how the proof given in the previous subsection can be modified in the case of constant mobility coefficients. Therefore, we assume that
	\[
	M(\b{s}) \equiv M, \qquad M \in \mathcal{M}.
	\]
	First, \eqref{eq:uniqueness0} simply reads
	\begin{equation} \label{eq:uniquenessC0}
		\begin{cases}
			\left\langle \pd{\b{c}}{t}, \: \b{v} \right\rangle_{V^*, V} + \left( M\nabla \b{w}, \nabla \b{v} \right)_H = 0, & \quad \forall \: \b{v} \in V, \\[0.5cm]
			\left(\b{w}, \b{z} \right)_H = \left( \Gamma \nabla \b{c} , \nabla \b{z} \right)_H + \left( \b P( \widetilde{\b{\mu}}_1 - \widetilde{\b{\mu}}_2),\, \b{z} \right)_H, & \quad \forall \: \b{z} \in V \cap L^\infty(\Omega;\mathbb{R}^3), \\[0.5cm]
			(\epsilon({c}_{1A},\,{c}_{1B})(\b{E}_0 - \nabla \Phi_1) - \epsilon({c}_{2A},\,{c}_{2B})(\b{E}_0 - \nabla \Phi_2) , \nabla \Lambda)_H = 0 & \quad \forall \: \Lambda \in V_0,
		\end{cases}
	\end{equation}
	or equivalently
	\begin{equation} \label{eq:uniquenessC1}
		\begin{cases}
			\b{w} - \overline{\b{w}} = -\mathcal{N}_{M}\pd{\b{c}}{t}, \\
			\left(\b{w}, \b{z} \right)_H = \left( \Gamma \nabla \b{c} , \nabla \b{z} \right)_H + \left( \b P( \widetilde{\b{\mu}}_1 - \widetilde{\b{\mu}}_2),\, \b{z} \right)_H, & \quad \forall \: \b{z} \in V \cap L^\infty(\Omega;\mathbb{R}^3),\\
			(\epsilon({c}_{1A},\,{c}_{1B})(\b{E}_0 - \nabla \Phi_1) - \epsilon({c}_{2A},\,{c}_{2B})(\b{E}_0 - \nabla \Phi_2) , \nabla \Lambda)_H = 0 & \quad \forall \: \Lambda \in V_0.
		\end{cases}
	\end{equation}
	Employing the same testing procedure as in the general case, the analogous of \eqref{eq:uniqueG9} is readily found:
	\begin{multline} \label{eq:uniqueC1}
		\gamma_*\|\nabla \b{c}\|^2_H + \left( \mathcal{N}_{M} \pd{\b{c}}{t}, \b{c} - \overline{\b{c}_0} \right)_H \leq \theta^*\sum_{i \in \mathcal{I}}|\overline{c_{0i}}|\left[\|\Psi'(c_{1i})\|_{L^1(\Omega)} + \|\Psi'(c_{2i})\|_{L^1(\Omega)} \right] \\
		+ C\left[ \|\b{c}\|_{L^4(\Omega;\mathbb{R}^3)}\|\b{c}-\overline{\b{c}_0}\|_{L^4(\Omega;\mathbb{R}^3)} + \|\b{c} - \overline{\b{c}_0}\|_H^2 + \|\nabla \Phi\|_{H}\|\b{c}-\overline{\b{c}_0}\|_{L^4(\Omega;\mathbb{R}^3)} \right].
	\end{multline}
	Therefore, following the rest of the proof line by line carries the result without additional regularity assumptions.
	\section{Proof of Theorem 
		\ref{thm:regularization} 
	} \label{sec:regularization}
	Let us now investigate regularization properties of the weak solution, i.e., we prove Theorem \ref{thm:regularization}. Let $h \geq 0$ and $X$ be a Banach space. For any $f: [0,T] \to X$, we denote by
	\[
	\dq{f} := \dfrac{f(t+h)-f(t)}{h}, \qquad 0\leq t \leq t+h \leq T,
	\]
	the (right) finite difference associated to $f(t)$. Moreover, for the sake of readability we shall denote time derivatives with the symbol $(\cdot)_t$ throughout this section only. The finite differences \review{of solutions }solve the problem
	\begin{equation*}
		\begin{cases}
			\left\langle(\dq{\b{c}})_t, \: \b{v} \right\rangle_{V^*, V} + \dfrac{1}{h}\left( M(\b{c}(t+h))\nabla \b{w}(t+h) -  M(\b{c}(t))\nabla \b{w}(t), \nabla \b{v} \right)_H = 0,\\
			\left(\dq{\b{w}}, \b{z} \right)_H = \left( \Gamma \nabla \dq{\b{c}} , \nabla \b{z} \right)_H + \dfrac{1}{h}\left( \b P( \widetilde{\b{\mu}}(t+h) - \widetilde{\b{\mu}}(t)),\, \b{z} \right)_H,\\
			\dfrac{1}{h}(\epsilon({c}_{A}(t+h),\,{c}_{B}(t+h))(\b{E}_0 - \nabla \Phi(t+h)) - \epsilon({c}_{A}(t),\,{c}_{B}(t))(\b{E}_0 - \nabla \Phi(t)) , \nabla \Lambda)_H = 0,
		\end{cases}
	\end{equation*}
	for every $\b{v} \in V$, $\b{z} \in V \cap L^\infty(\Omega;\mathbb{R}^3)$ and $\Lambda \in V_0$. Once again, we recast the problem in the form
	\begin{equation} \label{eq:reg0}
		\begin{cases}
			\dq{\b{w}} - \overline{\dq{\b{w}}} = -\dfrac{1}{h}\mathcal{N}_{\b{c}(t+h)}\b{c}_t(t+h) + \dfrac{1}{h}\mathcal{N}_{\b{c}(t)}\b{c}_t(t), \\
			\left(\dq{\b{w}}, \b{z} \right)_H = \left( \Gamma \nabla \dq{\b{c}} , \nabla \b{z} \right)_H + \dfrac{1}{h}\left( \b P( \widetilde{\b{\mu}}(t+h) - \widetilde{\b{\mu}}(t)),\, \b{z} \right)_H,\\
			\dfrac{1}{h}(\epsilon({c}_{A}(t+h),\,{c}_{B}(t+h))(\b{E}_0 - \nabla \Phi(t+h)) - \epsilon({c}_{A}(t),\,{c}_{B}(t))(\b{E}_0 - \nabla \Phi(t)) , \nabla \Lambda)_H = 0,
		\end{cases}
	\end{equation}
	Observe that $\dq{\b{c}} \in V_{(0)}$ by conservation of mass. Therefore, multiplying the first equation in \eqref{eq:reg0} by $\dq{\b{c}}$ yields
		\begin{multline} \label{eq:reg1}
			(\Gamma \nabla \dq{\b{c}}, \nabla \dq{\b{c}})_H + \dfrac{1}{h}\sum_{i \in \mathcal{I}} \theta_i (\Psi'(c_i(t+h)) - \Psi'(c_{i}(t)), \dq{c_i})_H + \dfrac{1}{h}\left(\td{I}{\b s} (\b{c}(t+h)) - \td{I}{\b s}  (\b{c}(t)), \dq{\b{c}}\right) _H \\
			+ \dfrac{1}{2h} \sum_{i \in \{A,B\}} \left(\pd{\epsilon}{c_i}(t+h)|\b{E}_0 - \nabla \Phi(t+h)|^2 -\pd{\epsilon}{c_{i}}(t)|\b{E}_0 - \nabla \Phi(t)|^2, \dq{c_i} \right)_H \\
			+ \sum_{(i,j) \in \{A,B\}^2} \alpha_{ij} (\mathcal{N}(\dq{c_i}), \dq{c_j})_H = -\left( \mathcal{N}_{\b{c}(t+h)} (\dq{\b{c}})_t, \dq{\b{c}}\right)_H -\dfrac{1}{h}\left( [\mathcal{N}_{\b{c}(t+h)} - \mathcal{N}_{\b{c}(t)}] \b{c}_t(t), \dq{\b{c}} \right)_H.	
		\end{multline}
		It is then clear that the structure of \eqref{eq:reg1} is very similar to the one of \eqref{eq:uniqueG9}. Therefore, we can follow the exact same computations of Subsection \ref{ssec:general}, with only notational modifications. In particular, applying the uniform Gronwall lemma in place of the standard one (see \cite[Chapter III, Lemma 1.1]{RTDDS}) leads immediately that
		\[
		\|\b{c}_t\|_{L^\infty(\sigma,t;V^*)} + \|\b{c}_t\|_{L^2(t,t+1;V)} \leq C_1,
		\]
		for every $t \geq \sigma$. Of course, once again, the additional regularity assumptions are not needed if the mobility coefficients are constant. Choosing $\b{v} = \b{w}-\overline{\b{w}}$ in the weak formulation given in Proposition \ref{prop:limitd2}, we have, by the Poincaré inequality
		\[
		\lambda_0 \|\nabla\b{w}\|^2_H \leq C\|\b{c}_t\|_{V^*}\|\nabla\b{w}\|_H,
		\]
		and thus $\nabla\b{w} \in L^\infty(\sigma,t;H)$. In order to bound its full $V$-norm, recall that the computations in the proof of Lemma \ref{lem:uniform2} show that
		\[
		\|\b{w}\|^2_V \leq C(1 + \|\nabla \b{w}\|^2_H),
		\]
		so we infer $\b{w} \in L^\infty(\sigma,t;V)$. If $\b{E}_0 \in W^{1,\alpha}(\Omega;\mathbb{R}^d)$, with $\alpha > \frac 65$, then, as a byproduct, owing to \eqref{eq:mu-ccomp}, we also have $\b{c} \in L^\infty(\sigma,t;H^2(\Omega;\mathbb{R}^3))$, hence, by the same reasoning showed in the proof of Proposition \ref{prop:higher4}, we further deduce $\Phi \in L^\infty(\sigma,t;W^{2,\beta}(\Omega))$ for the same $\beta > 1$ defined therein. Hence, for any $\sigma > 0$ 
		\[
		\|\Phi\|_{L^\infty(\sigma,t;W^{2,\beta}(\Omega))} + \|\b{w}\|_{L^\infty(\sigma,t;V)} + \|\b{c}\|_{L^\infty(\sigma,t;H^2(\Omega;\mathbb{R}^3))} + \sum_{i \in \mathcal{I}} \|\Psi'(c_i)\|_{L^\infty(\sigma,t;H)} \leq C
		\]
		provided that $t \geq \sigma$. 
		\section{Proof of Theorem \ref{thm:longtime}} \label{sec:longtime} 
		The present Section is devoted to proving Theorem \ref{thm:longtime}. The investigation of the asymptotic behavior of the weak solution of a Cahn--Hilliard system usually relies significantly on the validity of the strict separation property in some form. That means, loosely speaking, that the order parameters are uniformly bounded away from the pure states $0$ and $1$. Here, we provide a proof that makes use of the strict separation property only for large times, inspired by \cite{AW07}. The first result shows that the family of solution maps defined in Section \ref{sec:notation} possesses a semigroup structure.
		\begin{lem}
			The family of operators $\mathcal{S} := \{S(t)\}_{t \geq 0}$ defines a dynamical system over $\mathcal{X}$.
		\end{lem}
		\begin{proof}
			It amounts to verify the properties listed, for instance, in \cite[Definition 9.1.1]{CH98}.
			\begin{enumerate}[(i)] \itemsep 0.5em
				\item First, we prove that $S(0)$ is the identity map. Of course, we have
				\[
				S(0)(\b{c}_0) = \b{c}(0) = \b{c}_0,
				\]
				by definition of solution.  This shows that $S(0)$ is the identity map.
				\item Next, we show that $S(t) \in C^0(\mathcal{X};\mathcal{X})$ for all $t \geq 0$. If $t = 0$, the claim follows from (i). Otherwise, fix any $t > 0$, consider any pair  $\b{c}_{01}$ and $\b{c}_{02}$ of elements of $\mathcal{X}$ and set $\b{c}_k(t) := S(t)(\b{c}_{0k})$ for $k \in \{1,2\}$. We have
				\[
				\begin{split}
					\operatorname{dist}_{\mathcal{X}}(S(t)(\b{c}_{01}),\,S(t)(\b{c}_{02}) ) & =  \|\b{c}_{1}(t) - \b{c}_{2}(t)\|_V \\
					& \leq C\|\b{c}_{1}(t) - \b{c}_{2}(t)\|_{V^*}^\frac{1}{4}\|\b{c}_{1}(t) - \b{c}_{2}(t)\|^\frac{1}{4}_V \\
					& \leq \dfrac{1}{2}\operatorname{dist}_\mathcal{X}(S(t)(\b{c}_{01}),\,S(t)(\b{c}_{02}) ) + \|\b{c}_{1}(t) - \b{c}_{2}(t)\|_{V^*}^\frac{1}{3}.
				\end{split}
				\]
				Using the continuous dependence estimate given by Theorem \ref{thm:uniqueness} we conclude the argument.
				\item The concatenation property immediately follows since the system is autonomous.
				\item Finally, we show that $S(\cdot)(\b{c}_0) \in C^0([0,+\infty);\mathcal{X})$. This is a trivial restatement of the fact that $\b{c} \in C^0([0,+\infty);V)$.
			\end{enumerate}
			The proof is complete. 
		\end{proof} \noindent
		Let $\vartheta \in (\frac 12, 1)$ be arbitrary but fixed, and define
		\[
		\mathcal{Z} := H^{2\vartheta}(\Omega;\mathbb{R}^3),
		\] 
		equipped with its natural metric $\operatorname{dist}_\mathcal{Z}$.
		As customary, we consider the $\omega$-limit set given by
		\[
		\omega(\b{c}_0) := \left\{ \b{c}_0 \in \mathcal{X} \cap \mathcal{Z} : \exists \: \{t_n\}_{n \in \mathbb{N}} \text{ such that } t_n \nearrow +\infty \text{ and } \b{c}(t_n) \to \b{c} \text{ in } \mathcal{Z}\right\}.
		\]
		Some properties of the $\omega$-limit set can be deduced by relative compactness of the orbits. Observe that since, by definition,
		\[
		\omega(\b{c}_0) = \bigcap_{s \geq 0} \overline{\bigcup_{t \geq s} \{\b{c}(t)\}}^\mathcal{Z} = \bigcap_{s \geq \sigma} \overline{\bigcup_{t \geq s} \{\b{c}(t)\}}^\mathcal{Z}
		\]
		for all $\sigma > 0$ and $\b{c}$ is uniformly bounded in $H^2(\Omega;\mathbb{R}^3)$ for $t \geq 0$ (see Theorem \ref{thm:regularization}), the (truncated) orbit is relatively compact in $\mathcal{Z}$. Therefore, we can apply \cite[Theorem 9.1.8]{CH98} to get that the $\omega$-limit set is a non-empty, connected and compact subset of $\mathcal{Z}$. Furthermore, the $\omega$-limit set attracts the orbits of the system, as 
		\begin{equation} \label{eq:attracting}
			\lim_{t \to +\infty} \operatorname{dist}_\mathcal{Z}(S(t)(\b{c}_0),\, \omega(\b{c}_0)) = 0.
		\end{equation}
		Next, we establish the existence of a strict Lyapunov functional for the dynamical system $\mathcal{S}$.
		\begin{lem} \label{lem:lyapunov}
			The energy functional $\mathcal{E}$ is a strict Lyapunov functional for the dynamical system $\mathcal{S}$.
		\end{lem}
		\begin{proof}
			The claim follows if the only orbits conserving the energy of the system are equilibrium points. Assume that $\b{c}(t) = [c_A(t),\, c_B(t),\, c_S(t)]^T$ is a solution of \eqref{eq:strongform2}-\eqref{eq:conditions2} such that 
			\[
			\mathcal{E}(\b{c}(t), \Phi(t)) = 	\mathcal{E}(\b{c}(t),\, \mathscr{S}(c_A(t),\,c_B(t))) \equiv \mathcal{E}(\b{c}_0, \Phi_0) \qquad \forall \: t \geq 0.	
			\]
			The energy identity then implies that
			\[
			\lambda_0\|\nabla \b{w}\|^2_H \leq \int_{0}^{t} \int_\Omega M\nabla \b{w}(s) \cdot \nabla \b{w}(s) \: \d x \: \d s = 0,
			\]
			hence $\nabla \b{w}$ vanishes almost everywhere in $\Omega \times \mathbb{R}^+$. In turn, this implies that 
			the time derivative $\tpd{\b{c}}{t}$ also vanishes almost everywhere in $\Omega \times \mathbb{R}^+$, and therefore $\b{c} \equiv \b{c}_0$, and the claim is proved.
		\end{proof} \noindent
		In all of the remainder of the present Section, the symbol \reviewF{$\mathscr{S}:L^\infty(\Omega;\mathbb{R}^2) \to V_0$} denotes a realization of the solution map defined in Appendix \ref{app:A}, so that 
		\[
		\begin{cases}
			\div(\epsilon(c_A(t),c_B(t))(\b E_0 - \nabla \mathscr{S}(c_A(t),c_B(t)))) = 0 & \qquad \text{in } \Omega,\\
			\Phi = 0  & \qquad \text{on } \partial\Omega,
		\end{cases}
		\]
		for all $t \geq 0$.
		\begin{remark} \label{rem:phiconstant}
			Since $\Phi$ depends on time only through $\b{c}$, it is immediate to see that if $\b{c} \equiv \b{c}_0 = [c_{0A},\,c_{0B},\,c_{0S}]$, then $\Phi \equiv \mathscr{S}(c_{0A},\,c_{0B}) = \Phi_0$.
		\end{remark} \noindent
		Recalling the stationary problem \eqref{eq:stationary}-\eqref{eq:stationaryBCs}, we have the following result.
		\begin{lem} \label{lem:solutionsstat}
			The function $\b{c}^\infty = [c^\infty_A,\,c^\infty_B,\,c^\infty_S]^T \in H^2(\Omega;\mathbb{R}^3) \cap \mathcal{X}$ is an equilibrium point for $\mathcal{S}$ if and only if it solves \eqref{eq:stationary}-\eqref{eq:stationaryBCs} almost everywhere in $\Omega$, with $\Phi^\infty = \mathscr{S}(c^\infty_A,\,c^\infty_B)$. Furthermore, for any $\b{c}_0 \in \mathcal{X}$, we have that $\omega(\b{c}_0)$ only contains equilibrium points.
		\end{lem}
		\begin{proof}
			Assume that $\b{c}^\infty \in H^2(\Omega;\mathbb{R}^3) \cap \mathcal{X}$ is an equilibrium point for $\mathcal{S}$. Then, the claim follows immediately by the proof of Lemma \ref{lem:lyapunov} and Remark \ref{rem:phiconstant}. Viceversa, assume that $\b{c}^\infty \in H^2(\Omega;\mathbb{R}^3) \cap \mathcal{X}$ is a strong solution of \eqref{eq:stationary}-\eqref{eq:stationaryBCs}, with $\Phi^\infty = \mathscr{S}(c^\infty_A,\,c^\infty_B)$. Then, it is also immediate to show that if fulfills the non-stationary version \eqref{eq:strongform2}-\eqref{eq:conditions2}, with the first equation holding as a trivial identity. Uniqueness of solutions then implies the claim. The second part of the lemma follows instead from Lemma \ref{lem:lyapunov} and applying \cite[Theorem 9.2.7]{CH98}.
		\end{proof} \noindent
		We are now in a position to exhibit the strict separation property for large times.
		\begin{lem} \label{lem:largetimessp}
			Let $\vartheta \in (\frac d4, 1)$. For any $\b{c}_0 \in \mathcal{X}$, there exists a time instant $t^* \geq 0$ and a quantity $\delta \in (0,1)$ such that 
			\[
			\delta \leq c_i(\b{x}, t), \quad \forall \: (\b{x},t) \in \overline{\Omega} \times [t^*,+\infty), \quad \forall \: i \in \mathcal{I},
			\]
			where $\b{c}(t) := [c_A(t), c_B(t), c_S(t)]^T = S(t)(\b{c}_0)$.
		\end{lem}
		\begin{proof}
			Fix any $\b{c}_0 \in \mathcal{X}$ and consider $\b{c}^\infty\in \omega(\b{c}_0)$. By virtue of Lemma \ref{lem:solutionsstat}, $\b{c}^\infty$ is a strong solution of the stationary problem \eqref{eq:stationary}-\eqref{eq:stationaryBCs}, setting therein $\Phi^\infty = \mathscr{S}(c_A^\infty,\, c_B^\infty)$. First, we show that $\b{c}^\infty$ fulfills the strict separation proprerty. This can be done by using \cite[Theorem 8.1]{GP23}, once we show that 
			\[
			\b{f} := \b{P}\begin{bmatrix}
				\pd{I}{c_A^\infty} + \alpha_{AA}\mathcal{N}(c_A^\infty - \overline{c_A^\infty}) + \alpha_{AB}\mathcal{N}(c_B^\infty - \overline{c_B^\infty}) + \dfrac{1}{2}\pd{\epsilon}{c_A^\infty}|\b{E}_0 - \nabla \Phi^\infty|^2  \\[0.3cm] 
				\pd{I}{c_B^\infty} + \alpha_{BA}\mathcal{N}(c_A^\infty - \overline{c_A^\infty}) + \alpha_{BB}\mathcal{N}(c_B^\infty - \overline{c_B^\infty}) + \dfrac{1}{2}\pd{\epsilon}{c_B^\infty}|\b{E}_0 - \nabla \Phi^\infty|^2  \\[0.3cm]
				\pd{I}{c_S^\infty}
			\end{bmatrix} \in L^\infty(\Omega;\mathbb{R}^3).
			\]
			By the action of the projector $\b{P}$ (see \eqref{eq:projector}), it is not restrictive to drop it. The claim follows easily for the local and nonlocal interaction terms. For the electric term, instead, we need to prove additional regularity on $\Phi^\infty$. Observe that $\b{E}_0 \in L^\infty(\Omega;\mathbb{R}^d)$ by Assumption \ref{hyp:additionalrefined}. Moreover, since $\b{c}^\infty \in H^2(\Omega;\mathbb{R}^3)$, we can exploit the final part of the proof of Proposition \ref{prop:higher4} and the fact that $\beta > d$ to find $\Phi^\infty \in W^{2, \beta}(\Omega)$ and thus $\nabla \Phi^\infty \in L^\infty(\Omega;\mathbb{R}^d)$. This then proves, by arbitrarity of choice for $\b{c}_0$, that every element of the $\omega$-limit set is bounded away from zero, namely for every $\b{c}^\infty \in \omega(\b{c}_0)$ there exists $\varepsilon = \varepsilon(\b{c}^\infty) \in (0,1)$ such that
			\[
			\varepsilon \leq c_i^\infty(\b{x}), \qquad \forall\: \b{x} \in \overline{\Omega}, 	\quad \forall \: i \in \mathcal{I}.
			\]
			Since $\beta \in (\frac d4, 1)$, we have that $\mathcal{Z} \hookrightarrow C^0(\overline{\Omega};\mathbb{R}^3)$, thus the $\omega$-limit set is a compact set of $L^\infty(\Omega;\mathbb{R}^3)$. A standard compactness argument yields that $\varepsilon$ can be chosen uniformly with respect to $\b{c}^\infty \in \omega(\b{c}_0)$. Owing to the fact that the $\omega$-limit set attracts the orbit of the system in the $\mathcal{Z}$-metric, as given by \eqref{eq:attracting}, we have that for every $\b{c}_0 \in \mathcal{X}$ there exists some $t^* = t^*(\b{c}_0) \geq 0$ such that $S(t)(\b{c}_0)$ belongs to some open set $U \supset \omega(\b{c}_0)$ and
			\[
			0 < \delta \leq g(\b{x}), \qquad \forall\: \b{x} \in \overline{\Omega},
			\]
			for all $\b{g} \in U$ and for all $t \geq t^*$, with $\delta$ being some quantity in $(0,1)$ such that $\delta \leq \varepsilon$. The proof is complete.
		\end{proof} \noindent
		The asymptotic separation property given by Lemma \ref{lem:largetimessp} is crucial, since it essentially states that the singularities of $\Psi$ and its derivatives at 0 do not play any role in determining the asymptotic behavior of the system. Therefore, we can redefine the potential $\Psi$ on the set $(-\infty, \delta) \cup (1,+\infty)$, obtaining a regular function $\widetilde{\Psi}$ satisfying $\widetilde{\Psi} \in C^3(\mathbb{R})$ with uniformly bounded derivatives (up to order three) and satisfying
		\[
		\widetilde{\Psi}(s) = \Psi(s) \qquad \forall \: s \in [\delta, 1].
		\] 
		The function $I$ is also extended outside $[0,1]^3$, originating a function $\widetilde I$ satisfying $\widetilde I\in C^3(\mathbb{R}^3)$ and being such that its derivatives are uniformly bounded (up to order three). Moreover,
		\[
		\widetilde{I}(\b{s}) = I(\b{s}) \qquad \forall \: s \in [0,1]^3.
		\]
		Analogously, we also redefine the permittivity $\epsilon$ outside a suitably small closed neighbourhood $\mathcal{U}$ satisyfing $\mathcal{U}_0 \subset \mathcal{U} \subset \mathcal{U}_1$
		in such a way that the resulting extension $\widetilde{\epsilon} \in C^3(\mathbb{R}^2)$ is bounded, has bounded partial derivatives up to order three and is still affine-linear within $\mathcal{U}$.
		Let $\b{m} = [m_A,m_B,m_S]^T \in {\b{G}}$, and consider the modified energy functionals
		\begin{align*}
			\widetilde{\mathcal{E}}_1(\b{c}) & := \ii{\Omega}{}{\sum_{i \in \mathcal{I}} \left[  \dfrac{\gamma_i}{2}|\nabla c_i|^2 + \theta_i\widetilde\Psi(c_i + m_i) \right]}{x},\\
			\widetilde{\mathcal{E}}_2(\b{c})&  := \ii{\Omega}{}{\widetilde{I}(\b{c} + \b{m})}{x}, \\
			\widetilde{\mathcal{E}}_3(c_A,\, c_B) & := \sum_{(i,j) \in \{A,B\}^2} \dfrac{\alpha_{ij}}{2}\int_\Omega \int_\Omega G(\b{x},\b{y})c_i(\b{x})c_j(\b{y}) \:\d y \: \d x,\\
			\widetilde{\mathcal{E}}_4(c_A,\,c_B, \Phi) & := \dfrac{1}{2}\ii{\Omega}{}{\widetilde\epsilon(c_A + m_A,\,c_B + m_B)|\b{E}_0 - \nabla \Phi|^2}{x},
		\end{align*}
		as well as the collective modified version of the energy functional
		\[
		\widetilde{\mathcal{E}}(\b{c},\,\Phi) := \widetilde{\mathcal{E}}_1(\b{c}) + \widetilde{\mathcal{E}}_2(\b{c}) + \widetilde{\mathcal{E}}_3(\b{c}) + \widetilde{\mathcal{E}}_4(\b{c}, \Phi),
		\]
		defined in some domain
		\[
		\operatorname{dom}\widetilde{\mathcal{E}} \subset  (L^\infty(\Omega;\mathbb{R}^3) \cap TX^{\b{m}}) \times V_0 =: \mathcal{O}_0.
		\]
		Observe that $\widetilde{\mathcal{E}}(\b{c},\,\Phi) = \mathcal{E}(\b{c} + \b{m},\,\Phi)$, whence the domain of $\widetilde{\mathcal{E}}$ can be obtained by suitably projecting the one of $\mathcal{E}$. Next, we shall rephrase the energy functional as a function of only the order parameter $\b c$, exploiting to this end the solution map $\mathscr{S}$. Hence, we also define
		\[
		\widehat{\mathcal{E}} : \operatorname{dom} \widehat{\mathcal{E}} \subset \mathcal{O} :=  L^\infty(\Omega;\mathbb{R}^3) \cap T X^{\b  m} \to \mathbb{R}, \qquad \widehat{\mathcal{E}}(\b c) := \widetilde{\mathcal{E}}(\b c + \b m ,\,\mathscr{S}(c_A+m_A,\,c_B+m_B)),
		\]
		and the energy contributions $\widehat{\mathcal{E}}_i$ accordingly.
		For the sake of readability, in what follows, we may still denote the dependence of some quantity $a$ on the first two components of $\b c$ by writing $a(\b c)$ in place of $a(c_A,c_B)$. In order to characterize the first and second Fréchet derivatives of the modified energy functional of a single variable $\widehat{\mathcal{E}}$, we exploit the results proven in Appendix \ref{app:A}. Assume $\mathcal{O}$ to be endowed with its natural norm.
		\begin{lem} \label{lem:frechet2}
			The modified energy functional $\widehat{\mathcal{E}}$ satisfies the following properties:
			\begin{enumerate}[(i)]
				\item \label{me2:twice}the functional $\widehat{\mathcal{E}}$ is twice Fréchet-differentiable with continuous derivatives;
				\item \label{me2:hat}the first Fréchet derivative $D\widehat{\mathcal{E}} : \operatorname{dom} \widehat{\mathcal{E}} \to \mathcal{O}^*$ satisfies
				\begin{multline*}
					\langle D\widehat{\mathcal{E}}(\b{c}),\,\b{h} \rangle_{\mathcal{O}^*, \mathcal{O}} = \sum_{i \in \mathcal{I}} \gamma_i(\nabla c_i, \nabla h_i)_H + \theta_i \int_\Omega \widetilde{\Psi}'(c_i+m_i)h_i \: \d x + \int_\Omega \td{\widetilde{I}}{\b s}  (\b{c}+\b{m}) \cdot \b{h} \: \d x \\+ \sum_{(i,j) \in \{A,B\}^2} \alpha_{ij}(\nabla\mathcal{N}c_i, \nabla\mathcal{N}h_j)_H  +\dfrac{1}{2}\ii{\Omega}{}{|\b{E}_0 - \nabla \mathscr{S}(\b c + \b m)|^2\td{\widetilde{\epsilon}}{\b s} (\b c + \b m) \cdot \b{h}}{x} 
				\end{multline*}
				for all $\b{c} \in \operatorname{dom}\widehat{\mathcal{E}}$ and $\b{h} \in \mathcal{O}$;
				\item the second Fréchet derivative $D^2\mathcal{\widehat E}: \operatorname{dom} \widehat{\mathcal{E}} \to \mathcal{L}(\mathcal{O};\mathcal{O}^*)$ satisfies
				\begin{multline*}
					\langle D^2\widehat{\mathcal{E}}(\b{c})(\b h),\,\b k \rangle_{\mathcal{O}^*, \mathcal{O}} = \sum_{i \in \mathcal{I}} \gamma_i(\nabla h_i, \nabla k_i)_H + \theta_i \int_\Omega \widetilde{\Psi}''(c_i+m_i)h_ik_i \: \d x + \int_\Omega \td{^2 \widetilde{I}}{\b s^2} (\b{c}+\b{m})\b{k} \cdot \b{h} \: \d x\\ + \sum_{(i,j) \in \{A,B\}^2} \alpha_{ij}(\nabla\mathcal{N}h_i, \nabla \mathcal{N}k_j)_H + \dfrac{1}{2}\ii{\Omega}{}{|\b{E}_0 - \nabla \mathscr{S}(\b c + \b m)|^2\td{^2 \widetilde{\epsilon}}{\b s^2}(\b{c}+\b{m}) \b{k}\cdot \b{h}}{x} \\ 
					-\ii{\Omega}{}{\widetilde{\epsilon}(\b c + \b m)\nabla D\mathscr{S}(\b c + \b m )[\b h] \cdot \nabla D\mathscr{S}(\b c + \b m )[\b k]}{x}.
				\end{multline*}
				for all $\b{c} \in \operatorname{dom}\widehat{\mathcal{E}}$ and $\b{h}$, $\b k \in \mathcal{O}$.
			\end{enumerate}
		\end{lem}
		\begin{proof}
			Let us start by computing the first Fréchet derivative $\widehat{\mathcal{E}}$. For clarity, we split computations addressing each energy contribution separately.
			\begin{enumerate}[label=\fbox{$\widehat{\mathcal{E}}_{\arabic*}$}]
				\item The first contribution is a standard term appearing in the context of Cahn--Hilliard type equations. Referring, for example, to \cite[pp. 3829-3830]{scarpa21}, it is readily shown that
				\[
				\langle D\widehat{\mathcal{E}}_1(\b{c}),\,\b{h} \rangle_{\mathcal{O}^*, \mathcal{O}} = \sum_{i \in \mathcal{I}} \gamma_i(\nabla c_i, \nabla h_i)_H + \theta_i \int_\Omega \widetilde{\Psi}'(c_i+m_i)h_i \: \d x
				\]
				and
				\[
				\langle D^2\widehat{\mathcal{E}}_1(\b{c})\b{h},\,\b{k} \rangle_{\mathcal{O}^*,\mathcal{O}} = \sum_{i \in \mathcal{I}} \gamma_i(\nabla h_i, \nabla k_i)_H + \theta_i \int_\Omega \widetilde{\Psi}''(c_i+m_i)k_ih_i \: \d x
				\]
				for all $\b h$ and $\b k$ in $\mathcal O$.
				\item As a straightforward consequence of the mean value theorem, we get
				\[
				\langle D\widehat{\mathcal{E}}_2(\b{c}),\,\b{h} \rangle_{\mathcal{O}^*,\mathcal{O}} = \int_\Omega \td{\widetilde{I}}{\b s} (\b{c}+\b{m}) \cdot \b{h} \: \d x
				\]
				and
				\[
				\langle D^2\widehat{\mathcal{E}}_2(\b{c})\b{h},\,\b{k} \rangle_{\mathcal{O}^*,\mathcal{O}} = \int_\Omega \td{^2\widetilde{I}}{\b s^2} (\b{c}+\b{m})\b{k} \cdot \b{h} \: \d x
				\]
				for all $\b h$ and $\b k$ in $\mathcal O$.
				\item Exploiting the parallelogram law in $H$, we recast the contribution as
				\[
				\begin{split}
					\widehat{\mathcal{E}}_3(c_A,\, c_B) &  = \widetilde{\mathcal{E}}_3(c_A,\, c_B) \\ & = \sum_{(i,j) \in \{A,B\}^2} \dfrac{\alpha_{ij}}{2}\int_\Omega \int_\Omega G(\b{x},\b{y})c_i(\b{x})c_j(\b{y}) \:\d y \: \d x \\
					& = \sum_{(i,j) \in \{A,B\}^2} \dfrac{\alpha_{ij}}{2}(\mathcal{N}c_i, c_j)_H \\
					& = \dfrac{\alpha_{AA}}{2}\|\nabla \mathcal{N}c_A\|^2_H + \dfrac{\alpha_{AA}}{2}\|\nabla\mathcal{N} c_B\|_H^2 + \alpha_{AB}(\mathcal{N}c_A, c_B)_H \\
					& = \dfrac{\alpha_{AA}}{2}\|\nabla \mathcal{N}c_A\|^2_H + \dfrac{\alpha_{AA}}{2}\|\nabla\mathcal{N} c_B\|_H^2 + \dfrac{\alpha_{AB}}{4}\|\nabla\mathcal{N}(c_A + c_B)\|^2_H - \dfrac{\alpha_{AB}}{4}\|\nabla\mathcal{N}(c_A - c_B)\|^2_H.
				\end{split}
				\]
				Therefore, recalling that $\nabla \mathcal{N}$ is a linear operator on $H$, we have
				\[
				\langle D\widehat{\mathcal{E}}_3(\b{c}),\,\b{h} \rangle_{\mathcal{O}^*,\mathcal{O}} = \sum_{(i,j) \in \{A,B\}^2} \alpha_{ij}(\nabla \mathcal{N} c_i, \nabla \mathcal{N} h_j)_H
				\]
				and hence
				\[
				\langle D^2\widehat{\mathcal{E}}_3(\b{c})\b{h},\,\b{k} \rangle_{\mathcal{O}^*,\mathcal{O}} = \sum_{(i,j) \in \{A,B\}^2} \alpha_{ij}(\nabla \mathcal{N} h_i, \nabla \mathcal{N} k_j)_H
				\]
				for all $\b h$ and $\b k$ in $\mathcal O$.
				\item Finally, we tackle the electric term. First, we claim that, for all $ \b h \in \mathcal{O}$, 
				\begin{equation*}
					\langle D\widehat{\mathcal{E}}_4(\b{c}),\,\b{h}\rangle_{\mathcal{O}^*, \mathcal{O}} = \dfrac{1}{2}\ii{\Omega}{}{|\b{E}_0 - \nabla \mathscr{S}(\b c + \b m )|^2\td{\widetilde{\epsilon}}{\b s}(\b{c}+\b{m}) \cdot \b{h}}{x}.
				\end{equation*}
				To prove this, we proceed first by computing the variation $\widehat{\mathcal{E}}_4(\b{c} + \b h) - \widehat{\mathcal{E}}_4(\b{c})$. This yields, after elementary manipulations,
				\begin{equation} \label{eq:fr1}
					\small
					\begin{split}
						& \widehat{\mathcal{E}}_4(\b{c} + \b h) - \widehat{\mathcal{E}}_4(\b{c}) \\
						&\hspace{1cm}= \dfrac{1}{2}\ii{\Omega}{}{\left[\widetilde\epsilon(\b c + \b m + \b h) - \widetilde\epsilon(\b c + \b m) \right]|\b{E}_0 - \nabla \mathscr{S}(\b c + \b m)|^2}{x} \\
						&	\hspace{4cm}+ \dfrac{1}{2}\ii{\Omega}{}{\widetilde\epsilon(\b c + \b m + \b h)\left[ |\b{E}_0 - \nabla \mathscr{S}(\b c + \b m + \b h)|^2 - |\b{E}_0 - \nabla \mathscr{S}(\b c + \b m)|^2 \right]}{x} \\
						&\hspace{1cm}= \dfrac{1}{2}\ii{\Omega}{}{\left[\widetilde\epsilon(\b c + \b m + \b h) - \widetilde\epsilon(\b c + \b m) \right]|\b{E}_0 - \nabla \mathscr{S}(\b c + \b m)|^2}{x} \\
						&	\hspace{4cm}+ \dfrac{1}{2}\int_\Omega \widetilde\epsilon(\b c + \b m + \b h)\left[ (2\b{E}_0 - \nabla \mathscr{S}(\b c + \b m + \b h)- \nabla \mathscr{S}(\b c + \b m)) \right.\\
						& \hspace{7cm}\left. \cdot \, (\nabla \mathscr{S}(\b c + \b m) - \nabla \mathscr{S}(\b c + \b m + \b h)) \right] \, \d x. 
					\end{split}
				\end{equation}
				By the very definition of the solution map $\mathscr{S}$, it is easy to see that
				\[
				\dfrac{1}{2}\ii{\Omega}{}{\widetilde\epsilon(\b c + \b m + \b h)\left[ (\b{E}_0 - \nabla \mathscr{S}(\b c + \b m + \b h)) \cdot (\nabla \mathscr{S}(\b c + \b m) - \nabla \mathscr{S}(\b c + \b m + \b h)) \right]}{x} = 0,
				\]
				and therefore the second integral in the last term of \eqref{eq:fr1} reduces to
				\[
				\begin{split}
					&\dfrac{1}{2}\ii{\Omega}{}{\widetilde\epsilon(\b c + \b m + \b h)\left[ (\b{E}_0 - \nabla \mathscr{S}(\b c + \b m)) \cdot (\nabla \mathscr{S}(\b c + \b m) - \nabla \mathscr{S}(\b c + \b m + \b h)) \right]}{x} \\
					& \hspace{1cm}= \dfrac{1}{2}\ii{\Omega}{}{\widetilde\epsilon(\b c + \b m + \b h)\left[ (\nabla \mathscr{S}(\b c + \b m + \b h) - \nabla \mathscr{S}(\b c + \b m)) \cdot (\nabla \mathscr{S}(\b c + \b m) - \nabla \mathscr{S}(\b c + \b m + \b h)) \right]}{x},
				\end{split}
				\]
				by summing and subtracting $\nabla \mathscr{S}(\b c + \b m + \b h)$. Exploiting the definition of Fréchet differentiability for $\mathscr{S}$, i.e., by means of Proposition \ref{prop:frechet1}, we hence have that 
				\[
				\begin{split}
					&	\left| \dfrac{1}{2}\ii{\Omega}{}{\widetilde\epsilon(\b c + \b m + \b h)\left[ (\nabla \mathscr{S}(\b c + \b m + \b h) - \nabla \mathscr{S}(\b c + \b m)) \cdot (\nabla \mathscr{S}(\b c + \b m) - \nabla \mathscr{S}(\b c + \b m + \b h)) \right]}{x} \right|\\
					& \hspace{1cm} =\left| \dfrac{1}{2}\ii{\Omega}{}{\widetilde\epsilon(\b c + \b m + \b h)|\nabla\left[ D\mathscr{S}(\b c+ \b m)[\b h] + o(\b h)\right]|^2}{x} \right| = o(\|\b h\|_\mathcal{O}),				
				\end{split}
				\]
				where we also used the boundedness of $\widetilde{\epsilon}$ and the fact that $D\mathscr{S}(\b c+ \b m)$ is a continuous linear operator. Hence, owing to all of the above, we arrive to
				\begin{equation}
					\small
					\begin{split}
						& \left| \widehat{\mathcal{E}}_4(\b{c} + \b h) - \widehat{\mathcal{E}}_4(\b{c} ) - \langle D\widehat{\mathcal{E}}_4(\b{c}), \b h\rangle_{\mathcal{O}^*,\mathcal{O}}\right|\\
						& \hspace{1cm} = \dfrac{1}{2}\left| \ii{\Omega}{}{\left[\widetilde\epsilon(\b c + \b m + \b h) - \widetilde\epsilon(\b c + \b m) - \td{\widetilde\epsilon}{\b s}(\b c + \b m) \cdot \b h\right]|\b{E}_0 - \nabla \mathscr{S}(\b c + \b m)|^2}{x} \right| + o(\|\b h\|_{\mathcal{O}})\\
						& \hspace{1cm} = \dfrac{1}{2}\left| \ii{\Omega}{}{\left[o(|\b h|)\right]|\b{E}_0 - \nabla \mathscr{S}(\b c + \b m)|^2}{x} \right| + o(\|\b h\|_{\mathcal{O}})\\
						& \hspace{1cm} \leq o(\|\b h\|_{\mathcal{O}})\left( 1 + \|\b{E}_0 - \nabla \mathscr{S}(\b c + \b m)\|_{V_0}^2\right),
					\end{split}
				\end{equation} 
				hence the claim is shown. To prove it is continuous, we directly show its differentiability, arguing similarly. We claim that, for all $\b h$ and $\b k$ in $\mathcal{O}$, we have
				\begin{multline*}
					\langle D^2\widehat{\mathcal{E}}_4(\b{c})[\b k],\,\b{h}\rangle_{\mathcal{O}^*, \mathcal{O}} = \dfrac{1}{2}\ii{\Omega}{}{|\b{E}_0 - \nabla \mathscr{S}(\b c + \b m )|^2\td{{^2\widetilde\epsilon}}{\b s^2}(\b{c}+\b{m}) \cdot \b{h}}{x}\\
					-\ii{\Omega}{}{\widetilde{\epsilon}(\b c + \b m)\nabla D\mathscr{S}(\b c + \b m )[\b h] \cdot \nabla D\mathscr{S}(\b c + \b m )[\b k]}{x}.
				\end{multline*}
				Once again, we start from the variation
				\begin{equation} \label{eq:fr2}
					\small
					\begin{split}
						& \langle D\widehat{\mathcal{E}}_4(\b{c} + \b k), \b h\rangle_{\mathcal{O}^*,\mathcal{O}} - 	\langle D\widehat{\mathcal{E}}_4(\b{c}), \b h\rangle_{\mathcal{O}^*,\mathcal{O}} \\
						& \hspace{1cm} = \dfrac{1}{2}\ii{\Omega}{}{|\b{E}_0 - \nabla \mathscr{S}(\b c + \b m )|^2\left[\td{\widetilde{\epsilon}}{\b s}(\b{c}+\b{m}+ \b k) \cdot \b{h} -  \td{\widetilde{\epsilon}}{\b s}(\b{c}+\b{m}) \cdot \b{h} \right]}{x} \\
						& \hspace{2.5cm} + \dfrac{1}{2}\ii{\Omega}{}{\left[ |\b{E}_0 - \nabla \mathscr{S}(\b c + \b m + \b k)|^2-|\b{E}_0 - \nabla \mathscr{S}(\b c + \b m )|^2\right]\td{\widetilde{\epsilon}}{\b s}(\b{c}+\b{m} + \b k) \cdot \b{h}}{x} \\
						&\hspace{1cm}= \dfrac{1}{2}\ii{\Omega}{}{|\b{E}_0 - \nabla \mathscr{S}(\b c + \b m )|^2\left[\td{\widetilde{\epsilon}}{\b s}(\b{c}+\b{m}+ \b k) \cdot \b{h} -  \td{\widetilde{\epsilon}}{\b s}(\b{c}+\b{m}) \cdot \b{h} \right]}{x} \\
						&	\hspace{4cm}+ \dfrac{1}{2}\int_\Omega \left[ (2\b{E}_0 - \nabla \mathscr{S}(\b c + \b m + \b k)- \nabla \mathscr{S}(\b c + \b m)) \right.\\
						& \hspace{7cm}\left. \cdot \, (\nabla \mathscr{S}(\b c + \b m) - \nabla \mathscr{S}(\b c + \b m + \b k)) \right]\td{\widetilde{\epsilon}}{\b s}(\b{c}+\b{m} + \b k) \cdot \b{h} \, \d x. 
					\end{split}
				\end{equation}
				In order to simplify the expression above, we can use the characterization of the Fréchet derivative of $\mathscr{S}$ given in Proposition \ref{prop:frechet1} to infer that
				\[
				\dfrac{1}{2}\int_\Omega \left[ (\b{E}_0 - \nabla \mathscr{S}(\b c + \b m + \b k)) \cdot \, \nabla w \right]\td{\widetilde{\epsilon}}{\b s}(\b{c}+\b{m} + \b k) \cdot \b{h} \, \d x = \dfrac{1}{2}\int_\Omega \widetilde{\epsilon}(\b{c}+\b{m} + \b k) \nabla D \mathscr{S}(\b c + \b m + \b k)[\b h] \cdot \nabla w\, \d x
				\]
				for all $w \in V_0$, and therefore the last integral in \eqref{eq:fr2} reduces to 
				\begin{equation} \label{eq:fr3}			
					\begin{split}
						& \dfrac{1}{2} \int_\Omega \widetilde{\epsilon}(\b{c}+\b{m} + \b k) \nabla D \mathscr{S}(\b c + \b m + \b k)[\b h] \cdot (\nabla \mathscr{S}(\b c + \b m) - \nabla \mathscr{S}(\b c + \b m + \b k)) \: \d x \\
						& \hspace{1cm} + \dfrac{1}{2} \int_\Omega \left[(\b{E}_0 - \nabla \mathscr{S}(\b c + \b m)) \cdot (\nabla \mathscr{S}(\b c + \b m) - \nabla \mathscr{S}(\b c + \b m + \b k)) \right]\td{\widetilde{\epsilon}}{\b s}(\b{c}+\b{m} + \b k) \cdot \b{h} \, \d x \\
						& = \int_\Omega \widetilde{\epsilon}(\b{c}+\b{m} + \b k) \nabla D \mathscr{S}(\b c + \b m + \b k)[\b h] \cdot (\nabla \mathscr{S}(\b c + \b m) - \nabla \mathscr{S}(\b c + \b m + \b k)) \: \d x \\
						& \hspace{1cm} - \dfrac{1}{2} \int_\Omega \left| \nabla \mathscr{S}(\b c + \b m + \b k) - \nabla \mathscr{S}(\b c + \b m ) \right|^2\td{\widetilde{\epsilon}}{\b s}(\b{c}+\b{m} + \b k) \cdot \b{h} \, \d x \\
					\end{split}
				\end{equation}
				Observe that the last integral in \eqref{eq:fr3} satisfies
				\[
				\begin{split}
					& \left| \int_\Omega \left| \nabla \mathscr{S}(\b c + \b m + \b k) - \nabla \mathscr{S}(\b c + \b m ) \right|^2\td{\widetilde{\epsilon}}{\b s}(\b{c}+\b{m} + \b k) \cdot \b{h} \, \d x \right| \\
					& \hspace{1cm}=\left| \int_\Omega  \left|\nabla \left[ D\mathscr{S}(\b c + \b m)[\b k] + o(\b k)\right] \right|^2\td{\widetilde{\epsilon}}{\b s}(\b{c}+\b{m} + \b k) \cdot \b{h} \: \d x \right| \\
					& \hspace{1cm}= o(\|\b k\|_{\mathcal{O}}),			\end{split}
				\]
				still applying the boundedness of the first derivative of $\widetilde{\epsilon}$ and Proposition \ref{prop:frechet1}. On account of all of the previous relations, we obtain
				\begin{equation} \label{eq:fr4}
					\small
					\begin{split}
						& \langle D\widehat{\mathcal{E}}_4(\b{c} + \b k), \b h\rangle_{\mathcal{O}^*,\mathcal{O}} - 	\langle D\widehat{\mathcal{E}}_4(\b{c}), \b h\rangle_{\mathcal{O}^*,\mathcal{O}} - \langle D^2\widehat{\mathcal{E}}_4(\b{c})[\b k], \b h \rangle_{\mathcal{O}^*,\mathcal{O}}\\
						& \hspace{1cm} = \dfrac{1}{2}\ii{\Omega}{}{|\b{E}_0 - \nabla \mathscr{S}(\b c + \b m )|^2\left[\td{\widetilde{\epsilon}}{\b s}(\b{c}+\b{m}+ \b k) -  \td{\widetilde{\epsilon}}{\b s}(\b{c}+\b{m}) - \td{^2 \widetilde{\epsilon}}{\b s^2}(\b{c}+\b{m})\b k  \right]\cdot \b{h}}{x} \\
						& \hspace{2cm}-\int_\Omega \widetilde{\epsilon}(\b{c}+\b{m} + \b k) \nabla D \mathscr{S}(\b c + \b m + \b k)[\b h] \cdot (\nabla \mathscr{S}(\b c + \b m + \b k) - \nabla \mathscr{S}(\b c + \b m)) \: \d x \\
						& \hspace{3cm}+ \ii{\Omega}{}{\widetilde{\epsilon}(\b c + \b m)\nabla D\mathscr{S}(\b c + \b m )[\b h] \cdot \nabla D\mathscr{S}(\b c + \b m )[\b k]}{x} + o(\|\b k\|_{\mathcal{O}}) \\
						& \hspace{1cm} = \dfrac{1}{2}\ii{\Omega}{}{|\b{E}_0 - \nabla \mathscr{S}(\b c + \b m )|^2\left[o(\b k) \right]\cdot \b{h}}{x} \\
						& \hspace{2cm}-\int_\Omega \left[\widetilde{\epsilon}( \b c + \b m) + \td{\widetilde{\epsilon}}{\b s}( \b c + \b m) \cdot \b k + o(\b k)\right]\nabla \left[D\mathscr{S}(\b c + \b m)[\b h]  + D^2\mathscr{S}(\b c + \b m)[\b h, \b k] + o(\b k)\right] \\
						&\hspace{3cm}\cdot \nabla \left[ D\mathscr{S}(\b c + \b m)[\b k] + o(\b k)\right] \: \d x \\
						& \hspace{4cm}+ \ii{\Omega}{}{\widetilde{\epsilon}(\b c + \b m)\nabla D\mathscr{S}(\b c + \b m )[\b h] \cdot \nabla D\mathscr{S}(\b c + \b m )[\b k]}{x} + o(\|\b k\|_{\mathcal{O}}) 
					\end{split}
				\end{equation}
				It is then easily seen that the remainder so obtained decays superlinearly as $\b \|\b k\|_{\mathcal{O}} \to 0$ (we also exploited Proposition \ref{prop:frechet2}). Finally, we are left to show that the second Fréchet derivative is also continuous. Consider then two different elements $\b c_1$ and $\b c_2$ in $\mathcal{O}$. For all $\b h$ and $\b k$ in $\mathcal{O}$, we have
				\[ \small
				\begin{split}
					& \langle D^2\widehat{\mathcal{E}}_4(\b{c}_1)[\b k], \b h \rangle_{\mathcal{O}^*,\mathcal{O}} - \langle D^2\widehat{\mathcal{E}}_4(\b{c}_2)[\b k], \b h \rangle_{\mathcal{O}^*,\mathcal{O}} \\
					& = \dfrac{1}{2}\ii{\Omega}{}{|\b{E}_0 - \nabla \mathscr{S}(\b c_1 + \b m )|^2\td{{^2\widetilde\epsilon}}{\b s^2}(\b{c}_1+\b{m}) \cdot \b{h} -|\b{E}_0 - \nabla \mathscr{S}(\b c_2 + \b m )|^2\td{{^2\widetilde\epsilon}}{\b s^2}(\b{c}_2+\b{m}) \cdot \b{h} }{x} \\
					& \hspace{1cm}
					+\ii{\Omega}{}{\widetilde{\epsilon}(\b c_2 + \b m)\nabla D\mathscr{S}(\b c_2 + \b m )[\b h] \cdot \nabla D\mathscr{S}(\b c_2 + \b m )[\b k]-\widetilde{\epsilon}(\b c_1 + \b m)\nabla D\mathscr{S}(\b c_1 + \b m )[\b h] \cdot \nabla D\mathscr{S}(\b c_1 + \b m )[\b k]}{x}.
				\end{split}
				\]
				The first integral satisfies
				\[
				\begin{split}
					& \left|\dfrac{1}{2}\ii{\Omega}{}{|\b{E}_0 - \nabla \mathscr{S}(\b c_1 + \b m )|^2\td{{^2\widetilde\epsilon}}{\b s^2}(\b{c}_1+\b{m}) \cdot \b{h} -|\b{E}_0 - \nabla \mathscr{S}(\b c_2 + \b m )|^2\td{{^2\widetilde\epsilon}}{\b s^2}(\b{c}_2+\b{m}) \cdot \b{h} }{x}\right|\\
					& \leq \dfrac{1}{2}\left|\ii{\Omega}{}{\left[ |\b{E}_0 - \nabla \mathscr{S}(\b c_1 + \b m )|^2 - |\b{E}_0 - \nabla \mathscr{S}(\b c_2 + \b m )|^2 \right] \td{{^2\widetilde\epsilon}}{\b s^2}(\b{c}_1+\b{m}) \cdot \b{h}}{x}\right| \\
					& \hspace{1cm} + \dfrac{1}{2}\left|\ii{\Omega}{}{|\b{E}_0 - \nabla \mathscr{S}(\b c_2 + \b m )|^2 \left[ \td{{^2\widetilde\epsilon}}{\b s^2}(\b{c}_1+\b{m}) -\td{{^2\widetilde\epsilon}}{\b s^2}(\b{c}_2+\b{m}) \right] \cdot \b{h} }{x}\right| \\
					& \leq \dfrac{1}{2}\left|\ii{\Omega}{}{\left[ (2\b{E}_0 - \nabla \mathscr{S}(\b c_1 + \b m )- \nabla \mathscr{S}(\b c_2 + \b m )) \cdot (\nabla \mathscr{S}(\b c_2 + \b m )- \nabla \mathscr{S}(\b c_1 + \b m )) \right] \td{{^2\widetilde\epsilon}}{\b s^2}(\b{c}_1+\b{m}) \cdot \b{h}}{x}\right| \\
					& \hspace{1cm} + \dfrac{1}{2}\left|\ii{\Omega}{}{|\b{E}_0 - \nabla \mathscr{S}(\b c_2 + \b m )|^2 \left[ \td{{^2\widetilde\epsilon}}{\b s^2}(\b{c}_1+\b{m}) -\td{{^2\widetilde\epsilon}}{\b s^2}(\b{c}_2+\b{m}) \right] \cdot \b{h} }{x}\right| \\
					& \leq C\left[ 2\|\b E_0\|_{H} + \|\mathscr{S}(\b c_1 + \b m)\|_{V_0} + \|\mathscr{S}(\b c_2 + \b m)\|_{V_0}\right]\|\mathscr{S}(\b c_1 + \b m)-\mathscr{S}(\b c_2 + \b m)\|_{V_0}\|\b h\|_{\mathcal{O}} \\
					& \hspace{1cm} + C\left[ \|\b{E}_0\|_H^2 + \|\mathscr{S}(\b c_2 + \b m )\|_{V_0}^2\right]\|\b c_1 - \b c_2\|_{\mathcal{O}}\|\b h\|_{\mathcal{O}} \\
					& \leq C\left[ 1 + \|\b{E}_0\|_H^2 +\|\mathscr{S}(\b c_1 + \b m )\|_{V_0} + \|\mathscr{S}(\b c_2 + \b m )\|_{V_0}^2 \right]\|\b c_1 - \b c_2\|_{\mathcal{O}}\|\b h\|_{\mathcal{O}} 
				\end{split}
				\]
				where we applied Proposition \ref{prop:stability}. The proof is complete.
			\end{enumerate}		
		\end{proof} \noindent
		The second fundamental tool to show convergence to the equilibrium state is the \L ojasiewicz--Simon inequality. 
		\begin{prop} \label{prop:lojasiewicz}
			The energy functional $\widehat{\mathcal{E}}$ satisfies the \L ojasiewicz--Simon inequality locally around points belonging to some $\omega$-limit set, namely for all $\b{c}_0 \in \mathcal{X}$ and $\b{c}^\infty \in \omega(\b{c}_0)$, there exist two constants $\rho \in (0,\frac 12]$ and $R > 0$ such that
			\[
			|\widehat{\mathcal{E}}(\b{c}) - \widehat{\mathcal{E}}(\b{c}^\infty-\b{m})|^{1-\rho} \leq C\|D\widehat{\mathcal{E}}(\b{c})\|_{\reviewF{(TX^{\b m})}^*},
			\]
			for all $\b{c} \in \operatorname{dom}\widehat{\mathcal{E}}$ such that $\|\b{c}- \b{c}^\infty+\b m\|_{\mathcal{O}} \leq R$.
		\end{prop}
		\begin{proof}
		Our aim is to apply \cite[Theorem 1.2]{Rupp} (which follows from the more general result \cite[Corollary 3.11]{Chill03}). To this end, we check all the needed hypotheses:
			\begin{enumerate}[(i)]
				\item the space $\mathcal{O}$ embeds into \reviewF{$TX^{\b m}$} densely;
				\item  \reviewF{we show that, if $U$ is a sufficiently small neighborhood of $\b c^\infty - \b m$ in the space $\mathcal O$, then} the energy functional $\widehat{\mathcal{E}} : \reviewF{U} \to \mathbb{R}$ is analytic. Fixed any $\b c_0 \in \mathcal{X}$, let $t^*$ denote the corresponding asymptotic time of strict separation given by Lemma \ref{lem:largetimessp}. Since $\omega(\b c_0) = \omega(S(t^*)(\b c_0))$ and \reviewF{recalling that the $\omega$-limit set attracts the orbits of the system in $\mathcal{Z}$ and hence also in $\mathcal{O}$ (if $\theta$ is chosen to be sufficiently large), we may assume without loss of generality that}
				\[
				\delta \leq c_{i}(\b x, t) \leq 1-\delta \qquad  \forall \: i \in \mathcal{I},
				\] 
				for almost all $\b x \in \Omega$,  all $t \geq 0$ and for some $\delta \in (0,1)$\reviewF{, i.e., that the initial state is sufficiently close (with respect to the $\mathcal O$-metric) to some separated equilibrium state $\b c^\infty \in \omega(\b c_0)$ (cf. proof of Lemma \ref{lem:largetimessp}).} It is readily seen that $\widehat{\mathcal{E}}_1$ is analytic: this follows from the analyticity of $\widetilde{\Psi}$ in $(\delta,1-\delta)$ and from the fact that the first term of the derivative $D\widehat{\mathcal{E}}_1$ is linear (hence analytic) in $\b c$. Quite similarly, also $\widehat{\mathcal{E}}_2$ and $\widehat{\mathcal{E}}_3$ are analytic as well. We are then left to show that $\widehat{\mathcal{E}}_4$ is analytic. \reviewF{Recalling that
				\[
				\mathcal{U}_0 = \{ (x,y) \in \mathbb{R}^2: x \geq 0,\,y \geq 0,\, x+y\leq 1\}, 
				\]
				and that the electric permittivity $\epsilon$ is affine linear within $\mathcal{U}_0$,} the claim steadily follows observing that any strictly separated initial condition $\b c_0$ satisfies $[c_{0A}, \, c_{0B}]^T \in L^\infty(\Omega;\mathcal{U}_0)$, hence we have
				\[
				\widetilde{\mathcal{E}}_4(\b c, \Phi) = \dfrac{1}{2}\ii{\Omega}{}{\left( \b a \cdot \begin{bmatrix}
						c_A \\ c_B
					\end{bmatrix}+ \b b\right)|\b E_0 - \nabla \Phi|^2}{x}
				\]
				and its first derivative reads
				\[
				D\widetilde{\mathcal{E}}_4(\b c, \Phi)[\b h, \Lambda] = \dfrac{1}{2}\ii{\Omega}{}{\left( \b a \cdot \b h\right)|\b E_0 - \nabla \Phi|^2}{x} - \ii{\Omega}{}{\left( \b a \cdot \begin{bmatrix}
						c_A \\ c_B
					\end{bmatrix}+ \b b\right)(\b E_0 - \nabla \Phi)\cdot \nabla \Lambda}{x}
				\]
				\reviewF{for all $h,\,\Lambda \in \mathcal O$ and some fixed $\b a,\,\b b \in \mathbb R^2$.} The second term is affine-bilinear in both variables, hence analytic, while the first is a quadratic form (still analytic). Since the solution map $\mathscr{S}$ is also analytic by Proposition \ref{prop:analytic}, we get the result by composition of analytic maps;
				\reviewF{\item we show that the first derivative $D\widehat{\mathcal{E}} : U \to \mathcal{O}^*$ is an analytic map with values in $(TX^{\b m})^*$. Consider back the expression of $D\widehat{\mathcal{E}}$ given in Lemma \ref{lem:frechet2}, and assume that $\b h \in TX^{\b m}$ (and not necessarily in $\mathcal{O}$). We now show that $D\widehat{\mathcal{E}}$ is actually continuous with values in $TX^{\b m}$. Let $\b c_1$ and $\b c_2$ be different elements of $U$. Then,
							\begin{equation*}
									\begin{split}
											& \left| \langle D\widehat{\mathcal{E}}_1(\b{c}_1) - D\widehat{\mathcal{E}}_1(\b{c}_2),\,\b{h} \rangle_{V_{(0)}^*, V_{(0)}}\right| \\
											& \hspace{2cm}= \left| \sum_{i \in \mathcal{I}} \gamma_i(\nabla c_{1i} - \nabla c_{2i}, \nabla h_i)_H + \theta_i \int_\Omega \left[ \widetilde{\Psi}'(c_{1i}+m_i) - \widetilde{\Psi}'(c_{2i}+m_i)h_i \right] \: \d x\right| \\
											& \hspace{2cm} \leq C\|\b c_1 - \b c_2\|_{H}\|\b h\|_{V_{(0)}},
										\end{split}
								\end{equation*}
							by the mean value theorem and the Cauchy--Schwarz inequality. Analogously, we have
							\[
							\begin{split}
									\left| \langle D\widehat{\mathcal{E}}_2(\b{c}_1)-D\widehat{\mathcal{E}}_2(\b{c}_2),\,\b{h} \rangle_{V_{(0)}^*, V_{(0)}}\right| & = \left| \int_\Omega \left[ \td{\widetilde{I}}{\b s} (\b{c}_1+\b{m}) - \td{\widetilde{I}}{\b s} (\b{c}_2+\b{m}) \right] \cdot \b{h} \: \d x \right| \\
									& \leq C\|\b c_1 - \b c_2\|_{H}\|\b h\|_{V_{(0)}}
								\end{split}
							\]
							and
							\[
							\begin{split}
									\left| \langle D\widehat{\mathcal{E}}_3(\b{c}_1)-D\widehat{\mathcal{E}}_3(\b{c}_2),\,\b{h} \rangle_{V_{(0)}^*, V_{(0)}}\right| & = \left| \sum_{(i,j) \in \{A,B\}^2} \alpha_{ij}(\nabla \mathcal{N} c_{1i}-\nabla \mathcal{N} c_{2i}, \nabla \mathcal{N} h_j)_H \right| \\
									& \leq C\|\b c_1 - \b c_2\|_{H}\|\b h\|_{V_{(0)}}.
								\end{split}
							\]
							Concerning the last term, we have\
							\[
							\begin{split}
									& \left| \langle D\widehat{\mathcal{E}}_4(\b{c}_1)-D\widehat{\mathcal{E}}_4(\b{c}_2),\,\b{h} \rangle_{V_{(0)}^*, V_{(0)}}\right| \\ 
									& \hspace{1cm}= \dfrac{1}{2}\left|\ii{\Omega}{}{|\b{E}_0 - \nabla \mathscr{S}(\b c_1 + \b m )|^2\td{\widetilde{\epsilon}}{\b s}(\b{c}_1+\b{m}) \cdot \b{h} - |\b{E}_0 - \nabla \mathscr{S}(\b c_2 + \b m )|^2\td{\widetilde{\epsilon}}{\b s}(\b{c}_2+\b{m}) \cdot \b{h}}{x} \right| \\
									& \hspace{1cm} \leq \dfrac{1}{2}\left| \ii{\Omega}{}{\left[ |\b{E}_0 - \nabla \mathscr{S}(\b c_1 + \b m )|^2 -  |\b{E}_0 - \nabla \mathscr{S}(\b c_2 + \b m )|^2 \right]\td{\widetilde{\epsilon}}{\b s}(\b{c}_1+\b{m}) \cdot \b{h}}{x}  \right| \\
									& \hspace{4cm} + \dfrac{1}{2}\left|\ii{\Omega}{}{|\b{E}_0 - \nabla \mathscr{S}(\b c_2 + \b m )|^2\left[ \td{\widetilde{\epsilon}}{\b s}(\b{c}_1+\b{m}) - \td{\widetilde{\epsilon}}{\b s}(\b{c}_2+\b{m}) \right] \cdot \b{h}}{x} \right| \\
									& \hspace{1cm} \leq C\left[ \|\b E_0\|_{L^3(\Omega;\mathbb{R}^d)} + \|\mathscr{S}(\b c_1 + \b m )\|_{W^{1,3}(\Omega)} + \|\mathscr{S}(\b c_2 + \b m )\|_{W^{1,3}(\Omega)}\right]\|\b h\|_{L^6(\Omega;\mathbb{R}^2)}\|\mathscr{S}(\b c_1)-\mathscr{S}(\b c_2)\|_{V_0} \\
									& \hspace{4cm} + C\left[ \|\b E_0\|_{L^4(\Omega;\mathbb{R}^d)}^2 + \|\mathscr{S}(\b c_2 + \b m )\|_{W^{1,4}(\Omega)}^2 \right]\|\b h\|_{L^4(\Omega;\mathbb{R}^2)}\|\b c_1-\b c_2\|_{L^4(\Omega;\mathbb{R}^3)} \\
									& \hspace{1cm} \leq C\left[ 1+ \|\b E_0\|_{L^4(\Omega;\mathbb{R}^d)}^2 + \|\mathscr{S}(\b c_1 + \b m )\|_{W^{1,3}(\Omega)}+ \|\mathscr{S}(\b c_2 + \b m )\|_{W^{1,4}(\Omega)}^2 \right]\|\b c_1 - \b c_2\|_{L^\infty(\Omega;\mathbb{R}^3)}\|\b h\|_{V_{(0)}} \\
								\end{split}
							\]
							by means of Proposition \ref{prop:stability} and Sobolev embeddings. Observe that the additional regularity required on $\b c_1$ and $\b c_2$ is granted instantaneously after detaching from the initial conditions from the dynamics of the system (see Theorem \ref{thm:regularization} and, in turn, Proposition \ref{prop:compact}). Analyticity follows following the same considerations in item (i).
				}
				\item the Fredholm index of $D^2\widehat{\mathcal{E}}(\b c^\infty): \mathcal{O} \to \reviewF{(TX^{\b m})}^*$ is zero. To show this, we analyze the explicit expression given by Lemma \ref{lem:frechet2}. Indeed, recall that
				\begin{multline*}
					\langle D^2\widehat{\mathcal{E}}(\b{c}^\infty)[\b h],\,\b k \rangle_{V_{(0)}^*, V_{(0)}} = \sum_{i \in \mathcal{I}} \gamma_i(\nabla h_i, \nabla k_i)_H + \theta_i \int_\Omega \widetilde{\Psi}''(c_i^\infty+m_i)h_ik_i \: \d x + \int_\Omega \td{^2 \widetilde{I}}{\b s^2} (\b{c}^\infty+\b{m})\b{k} \cdot \b{h} \: \d x\\ + \sum_{(i,j) \in \{A,B\}^2} \alpha_{ij}(\nabla\mathcal{N}h_i, \nabla \mathcal{N}k_j)_H + \dfrac{1}{2}\ii{\Omega}{}{|\b{E}_0 - \nabla \mathscr{S}(\b c^\infty + \b m)|^2\td{^2 \widetilde{\epsilon}}{\b s^2}(\b{c}^\infty+\b{m}) \b{k}\cdot \b{h}}{x} \\ 
					-\ii{\Omega}{}{\widetilde{\epsilon}(\b c^\infty + \b m)\nabla D\mathscr{S}(\b c^\infty + \b m )[\b h] \cdot \nabla D\mathscr{S}(\b c^\infty + \b m )[\b k]}{x}.
				\end{multline*}
				\reviewF{Consider the operator $B: \mathcal{O} \to \reviewF{(TX^{\b m})}^*$ such that
					\[
					\langle B\b h, \b k\rangle_{(TX^{\b m})^*, TX^{\b m}} = \sum_{i \in \mathcal{I}} \gamma_i(\nabla h_i, \nabla k_i)_H,
					\] 
					for all $\b h \in \mathcal O$ and $\b k \in TX^{\b m}$. The operator $B$ can be regarded as the restriction of the Laplace operator to the space $\mathcal O$. It is then straightforward to see that $\ker B$ is trivial. Its adjoint operator $B^*: ( TX^{\b m})^{**} \cong TX^{\b m} \to \mathcal O^*$ satifies by definition
					\[
					\langle B^*\b h^*, \b k\rangle_{\mathcal O^*, \mathcal O} =  \langle\b h^*,  B\b k\rangle_{(TX^{\b m})^{**}, (TX^{\b m})^*} = \langle B\b k, \b h^*  \rangle_{(TX^{\b m})^{*}, TX^{\b m}},
					\] 
					for all $\b h^* \in (TX^{\b m})^{**}$ and $\b k \in \mathcal O$, where the last equality follows identifying the $TX^{\b m}$ with its bidual through the canonical map. This implies that $\ker B^*$ is trivial as well, and, in turn, that $B$ is a Fredholm operator of index $0$. Indeed, if $\b h^* \in \ker B^*$, then 
					\[
					\langle B\b k, \b h^*  \rangle_{(TX^{\b m})^{*}, TX^{\b m}} = 0
					\]
					for all $\b k \in \mathcal O$. In particular, along some sequence $\{\b h^*_n\}_{n \in \mathbb N} \subset \mathcal O$ such that $\b h_n^* \to \b h^*$ in $TX^{\b m}$, we have
					\[
					0 = \langle B\b h_n^*, \b h^*  \rangle_{(TX^{\b m})^{*}, TX^{\b m}} = \sum_{i \in \mathcal{I}} \gamma_i(\nabla h_{n,i}^*, \nabla h_i^*)_H.
					\]
					for all $n \in \mathbb N$. Since we have $\b h_n^* \rightharpoonup \b h^*$ in $TX^{\b m}$ as well, passing to the limit as $n \to +\infty$ yields that $\b h^* = 0$. Finally, we argue that the remaining terms are compact.}  Indeed,
				\begin{enumerate}
					\item the second summand, depending on $\widetilde{\Psi}''$, clearly gives a compact operator with values in $(TX^{\b m})^*$, as we can identify it with the multiplication operator
					\[
					\b h \mapsto \begin{bmatrix}
						\theta_A\widetilde{\Psi}''(c_A^\infty+m_A) & 0 & 0 \\
						0 & \theta_B\widetilde{\Psi}''(c_B^\infty+m_B) & 0 \\
						0 &0 & \theta_S\widetilde{\Psi}''(c_S^\infty+m_S)
					\end{bmatrix}\b h
					\]
					which has even values in $H$, given that $\b c^\infty$ is strictly separated,
					\item the same holds true for the third summand, where the multiplication matrix is the Hessian matrix of $\widetilde{I}$;
					\item the fourth summand is a realization of the \reviewF{restriction of the} inverse Laplace operator $\mathcal{N}: V_{(0)}^* \to V_{(0)}$ \reviewF{to the space $\mathcal O$}, and hence is compact;
					\item the fifth summand, \reviewF{owing also to Proposition \ref{prop:compact}}, also defines a compact multiplication operator;
					\item the sixth summand, lastly, can be written in the following form. Let 
					\[
					\langle J(\b{c}^\infty)[\b h],\,\b k \rangle_{\mathcal{O}^*, \mathcal{O}} := \ii{\Omega}{}{\widetilde{\epsilon}(\b c^\infty + \b m)\nabla D\mathscr{S}(\b c^\infty + \b m )[\b h] \cdot \nabla D\mathscr{S}(\b c^\infty + \b m )[\b k]}{x}.
					\]
					If we consider the continuous bilinear form
					\[
					A_0: V_0 \times V_0 \to \mathbb{R}, \qquad 
					A_0(u, v) := \ii{\Omega}{}{\widetilde{\epsilon}(\b c^\infty + \b m)\nabla u \cdot \nabla v}{x}
					\]
					we have
					\[
					\langle J(\b{c}^\infty)[\b h],\,\b k \rangle_{\mathcal{O}^*, \mathcal{O}} = A_0(D\mathscr{S}(\b c^\infty + \b m )[\b h],\, D\mathscr{S}(\b c^\infty + \b m )[\b k]]).
					\]
					Since $\b c^\infty \in H^2(\Omega;\mathbb{R}^3)$ (cf. Lemma \ref{lem:solutionsstat}), and owing to Assumption \ref{hyp:additionalrefined} and Proposition \ref{prop:compact}, we infer the compactness of $J(\b c^\infty)$ as well.
				\end{enumerate}
			\end{enumerate}
			The claim is therefore a direct application of \cite[Theorem 1.2]{Rupp}. 
		\end{proof} \noindent
		We are now in a position to show Theorem \ref{thm:longtime}.
		\begin{proof}[Proof of Theorem \ref{thm:longtime}]
			As a consequence of the LaSalle's invariance principle, we have that the energy functional is constant over elements of the same $\omega$-limit set, namely
			\[
			\restr{\widehat{\mathcal{E}}(\b c)}{\omega(\b c_0)-\b m} = \restr{\widehat{\mathcal{E}}(\b c,\, \mathscr{S}(c_A,\,c_B))}{\omega(\b c_0)-\b m} \equiv \mathcal{E}_\infty,
			\]
			for some limit energy $\mathcal{E}_\infty \in \mathbb{R}$. Next, since the $\omega$-limit set $\omega(\b c_0)$ is compact in $\mathcal{Z}$, the translated set $\omega(\b c_0) - \b m$ is compact in $\mathcal{Z}$ as well. As $\mathcal{Z} \hookrightarrow \mathcal{O}$ compactly (granted that $\theta > \frac d4$), we can consider an open cover $U_0$ formed by a finite number of sufficiently small $\mathcal{O}$-balls, namely
			\[
			\omega(\b c_0) - \b m \subset U_0 := \bigcup_{i = 1}^N  B_{\varpi_i}(\b c^i),
			\]
			for some $N \in \mathbb{N}$, and for a finite family of centers $\{\b c^i\}_{i=1}^N \subset \omega(\b c_0) - \b m$ and radii $\varpi_i > 0$ satisfying
			\[ \max_{1 \leq i \leq N} \varpi_i \leq \min_{_{1 \leq i \leq N}}R(\b c^i),\]
			where $R(\b c^i)$ is the radius appearing in Proposition \ref{prop:lojasiewicz} corresponding to the element $\b c^i + \b m \in \omega(\b c_0)$. Therefore, in each of the balls $B_{\varpi_i}(\b c^i)$ we can find positive constants $C_i$, $\rho_i$ ensuring local validity of the \L ojasiewicz--Simon inequality, and since they are a finite number, we can extract uniform constants $C > 0$ and $\rho \in (0,\frac{1}{2}]$ such that
			\begin{equation}
				\label{eq:LSin}
				|\widehat{\mathcal{E}}(\b{c}) - \widehat{\mathcal{E}}(\b{c}^\infty-\b{m})|^{1-\rho} = |\widehat{\mathcal{E}}(\b{c}) - \mathcal{E}_\infty|^{1-\rho} \leq C\|D\widehat{\mathcal{E}}(\b{c})\|_{\mathcal{O}^*}, 
			\end{equation}
			holding for all $\b c \in U_0$. In the proof of the asymptotic separation property (Lemma \ref{lem:largetimessp}), we proved the existence of some open set $U$ covering $\omega(\b c_0)$ within which every trajectory is eventually attracted, becoming strictly separated from the pure phases. Thus, there exists a time instant $t^\sharp > 0$ such that $\b c(t) \in (U - \b m) \cap U_0$ for every $t \geq t^\sharp$. Accordingly, let us consider the functional $H : [t^\sharp, +\infty) \to \mathbb{R}$ defined by
			\[
			H(t) = (\widehat{\mathcal{E}}(\b c(t)) - \mathcal{E}_\infty)^\rho.
			\]
			Then, a straightforward computation, jointly with the differential form of the energy identity and Proposition \ref{prop:lojasiewicz}, entail that
			\begin{equation} \label{eq:conv1}
				\begin{split}
					-\td{}{t}H(t) & = -\rho (\widehat{\mathcal{E}}(\b c(t)) - \mathcal{E}_\infty)^{\rho-1} \td{}{t}\widehat{\mathcal{E}}(\b c(t)) \\ 
					& \geq C\rho\|D\widehat{\mathcal{E}}(\b{c})\|_{\mathcal{O}^*}^{-1} \int_\Omega M(\b{c}(t))\nabla \b{w}(t) \cdot \nabla \b{w}(t) \: \d \b x \\
					& \geq C\rho\lambda_0\|D\widehat{\mathcal{E}}(\b{c})\|_{\mathcal{O}^*}^{-1} \|\nabla \b w(t)\|_{H}^2.
				\end{split}
			\end{equation}
			In order to develop further \eqref{eq:conv1}, we study the quantity $\|D\widehat{\mathcal{E}}(\b{c})\|_{\mathcal{O}^*}$. By definition,
			\[
			\begin{split}
				\|D\widehat{\mathcal{E}}(\b{c})\|_{\mathcal{O}^*}  & = \sup_{\|\b h\|_{\mathcal{O}} = 1} |\langle D\widehat{\mathcal{E}}(\b{c}), \b h \rangle_{\mathcal{O}^*, \mathcal{O}}| \\
				& = \sup_{\|\b h\|_{\mathcal{O}} = 1} \left| \sum_{i \in \mathcal{I}} \gamma_i(\nabla c_i, \nabla h_i)_H + \theta_i \int_\Omega \widetilde{\Psi}'(c_i+m_i)h_i \: \d x + \int_\Omega \td{\widetilde{I}}{\b s} (\b{c}+\b{m}) \cdot \b{h} \: \d x \right.\\
				& \hspace{1cm} \left. + \sum_{(i,j) \in \{A,B\}^2} \alpha_{ij}(\nabla\mathcal{N}c_i, \nabla\mathcal{N}h_j)_H  +\dfrac{1}{2}\ii{\Omega}{}{|\b{E}_0 - \nabla \mathscr{S}(\b c + \b m)|^2\td{\widetilde \epsilon}{\b s}(\b c + \b m) \cdot \b{h}}{x} \right| \\
				& = \sup_{\|\b h\|_{\mathcal{O}} = 1} (\b w, \b h)_H \\
				& \leq \reviewF{C}\|\nabla \b w\|_H
			\end{split}
			\]
			and hence, from \eqref{eq:conv1}
			\[
			\begin{split}
				- \td{}{t}H(t) & \geq C\rho\lambda_0\|D\widehat{\mathcal{E}}(\b{c})\|_{\mathcal{O}^*}^{-1} \|\nabla \b w(t)\|_{H}^2 \\
				& \geq C\rho\lambda_0\|\nabla \b w(t)\|_{H}
			\end{split}
			\]
			and an integration on $[t^\sharp, +\infty)$ yields that $\nabla \b w \in L^1([t^\sharp,+\infty), H)$. By comparison in the evolution equations, we also have $\partial_t \b c \in L^1([t^\sharp,+\infty), V^*_{(0)})$. Therefore, the limit
			\[ \lim_{t \to +\infty} \b c(t) = \b c^\infty \]
			holds in the $V_{(0)}^*$-sense and, by compactness of the $\omega$-limit set, also in the $H^{2\vartheta}$-sense, for all $\vartheta \in (0,1)$. The convergence of $\Phi(t)$ to $\Phi^\infty = \mathscr{S}(\b c^\infty)$ follows from the continuity of $\mathscr{S}$, given by Proposition \ref{prop:stability}.
		\end{proof} \noindent
		\appendix
		\section{On the regularity of a certain solution map}
		\label{app:A}
		In this appendix, we investigate the Fréchet regularity of the solution map of an elliptic problem with nonconstant coefficients resembling the equation for the electric potential $\Phi$ appearing in \eqref{eq:strongform2}. To this end, let $\Omega \subset \mathbb{R}^d$, with $d \in \{2,3\}$, be a bounded domain with boundary $\partial \Omega$ of class $C^{1,1}$ and consider the partial differential equation
		\begin{equation} \label{eq:elliptic}
			\begin{cases}
				-\div(a(\b \eta)\nabla u) = -\div(a(\b \eta)\b v) & \quad \text{in } \Omega, \\
				u = 0 & \quad \text{on }\partial\Omega.
			\end{cases}
		\end{equation}
		where the given data comply with the following assumption:
		\begin{enumerate}[label = \textbf{(H)}]
			\item The function $a: \mathbb{R}^2 \to \mathbb{R}$ is of class $C^3(\mathbb{R}^2)$, affine-linear in an open neighbourhood $\mathcal{U}$ of the set
			\[
			\mathcal{U}_0 := \{ (x,y) \in \mathbb{R}^2: x \geq 0,\,y \geq 0,\, x+y\leq 1\}, 
			\]
			and such that
			\[
			0 < a_* \leq a(\b s) \leq a^* \qquad \forall \: \b s \in \mathbb{R}^2.
			\]
			Furthermore, we assume the partial derivatives of $a$ up to order three to be bounded over $\mathbb{R}^2$.
			Finally, the vector field $\b v \in L^\infty(\Omega;\mathbb{R}^d)$.
			\label{hyp:appendix} 
		\end{enumerate} 
		In the following, we consider $\b \eta := (\eta_1, \eta_2) \in L^\infty(\Omega;\mathbb{R}^2)$ as a control parameter, while $\b v$ is to be fixed once and for all. 
		\begin{prop} \label{prop:existence1}
			Assume Assumption \ref{hyp:appendix}. For any $\b \eta \in L^\infty(\Omega;\mathbb{R}^2)$, problem \eqref{eq:elliptic} admits a unique weak solution $u \in V_0$. Furthermore, there exists a constant $C$ independent of $u$ and $\b \eta$ such that
			\[
			\|u\|_{V_0} \leq C.
			\]
		\end{prop} 
		\begin{proof}
			It amounts to apply the Lax--Milgram theorem. Indeed, we have that, for any $u \in V_0$, we have
			\[
			(a(\b \eta)\nabla u, \nabla u)_H \geq a_*\|u\|^2_{V_0},
			\]
			and on the other hand,
			\[
			\left| (a(\b \eta)\nabla u, \nabla v)_H \right| \leq a^*\|u\|_{V_0}\|v\|_{V_0}
			\]
			\reviewF{for all $u,\,v \in V_0$}.
			Furthermore, the map
			\[
			v \mapsto (a(\b \eta)\b v, \nabla v)_H
			\]
			defines an element of $V_0^*$. This is enough to conclude the existence and uniqueness of a solution $u$ to problem \eqref{eq:elliptic}. In particular, since
			\[
			\left| (a(\b \eta)\b v, \nabla v)_H \right| \leq a^*\|\b v\|_H \|v\|_{V_0},
			\]
			the thesis follows with $C = a^*a_*^{-1}\|\b v\|_H$.
		\end{proof} \noindent
		Proposition \ref{prop:existence1} shows that the solution map $\mathscr{S}: L^\infty(\Omega;\mathbb{R}^2) \to V_0$, associating to each $\b \eta \in L^\infty(\Omega;\mathbb{R}^2)$ the corresponding $u \in V_0$ solving \eqref{eq:elliptic}, is well defined. The next result proves it is Lipschitz-continuous.
		\begin{prop} \label{prop:stability}
			There exist a constant $C > 0$ such that
			\[
			\|\mathscr{S}(\b \eta^1) - \mathscr{S}(\b \eta^2)\|_{V_0} \leq C \|\b \eta^1 - \b \eta^2\|_{L^\infty(\Omega;\mathbb{R}^2)}
			\]
			for all $\b \eta^1$ and $\b \eta^2 \in L^\infty(\Omega;\mathbb{R}^2)$.
		\end{prop}
		\begin{proof}
			Let $u := \mathscr{S}(\b \eta^1)$ and $v := \mathscr{S}(\b \eta^2)$. By definition of the solution map $\mathscr{S}$ the identity
			\[
			(a(\b \eta^1)\nabla u - a(\b \eta^2)\nabla v, \nabla w) = ((a(\b \eta^1) - a(\b \eta^2))\b v, \nabla w)_H
			\]
			holds for all $w \in V_0$. Moreover, we have, by summing and subtracting $a(\b \eta^1)\nabla v$ at left hand side, 
			\[
			(a(\b \eta^1)\nabla (u-v) + (a(\b \eta^1)-a(\b \eta^2))\nabla v, \nabla w)_H = ((a(\b \eta^1) - a(\b \eta^2))\b v, \nabla w)_H,
			\]
			for all $w \in V_0$. Choose now $w = u-v$. On one hand,
			\[
			(a(\b \eta^1)\nabla (u-v), \nabla(u-v))_H \geq a_*\|u-v\|^2_{V_0}
			\]
			and on the other, by the mean value theorem,
			\[
			|((a(\b \eta^1)-a(\b \eta^2))(\b v -\nabla v), \nabla (u-v))_H| \leq \|\b \eta^1- \b \eta^2\|_{L^\infty(\Omega;\mathbb{R}^2)}\|\b v - \nabla v\|_H\|\nabla(u-v)\|_H\sup_{\b x \in \mathbb{R}^2}  \left\|\td{a}{\b s} (\b x)\right\|_{\mathbb{R}^2}.
			\]
			By Proposition \ref{prop:existence1}, we get
			\[
			a_*\|u-v\|_{V_0} \leq 2a^*a_*^{-1}\left( 1 + a^*a_*^{-1}\right)\| \b v\|_H^2\sup_{\b x \in \mathbb{R}^2}\left\|\td{a}{\b s} (\b x)\right\|_{\mathbb{R}^2}\|\b \eta^1- \b \eta^2\|_{L^\infty(\Omega;\mathbb{R}^2)},
			\]
			and the claim is therefore shown with 
			\[
			C := 2a^*a_*^{-2}\left( 1 + a^*a_*^{-1}\right)\| \b v\|_H^2\sup_{\b x \in \mathbb{R}^2}\left\|\td{a}{\b s} (\b x)\right\|_{\mathbb{R}^2}.
			\]
			The proof is complete.
		\end{proof} \noindent
		The next result concerns the Fréchet differentiability of the solution map $\mathscr{S}$. The following proofs are inspired, for instance, by \cite{Blank2014}. Indeed, problems of this kind naturally appear in optimal control settings.
		\begin{prop} \label{prop:frechet1}
			The map \reviewF{$\mathscr{S}:  L^\infty(\Omega;\mathbb{R}^2) \to V_0$} is continuously Fréchet differentiable. For any $\b \eta$ and $\b h$ in $L^\infty(\Omega;\mathbb{R}^2)$, the derivative $D\mathscr{S}: L^\infty(\Omega;\mathbb{R}^2) \to \mathcal{L}(L^\infty(\Omega;\mathbb{R}^2), V_0)$ is \reviewF{given} by
			\[
			D\mathscr{S}(\b \eta)[\b h] = u^*,
			\]
			where $u^* \in V_0$ is the unique weak solution of the problem
			\[
			\begin{cases}
				-\div\left(a(\b \eta) \nabla u^*+\left[\td{a}{\b s}(\b \eta) \cdot \b h\right]\nabla u\right) = -\div\left( \left[\td{a}{\b s}(\b \eta) \cdot \b h\right]\b v\right) & \quad \text{in } \Omega, \\
				u^* = 0 & \quad \text{on }\partial\Omega,
			\end{cases}
			\]
			provided that $u := \mathscr{S}(\b \eta)$. Moreover, for all $\b \eta_1$ and $\b \eta_2$ in $L^\infty(\Omega; \mathbb{R}^2)$ it holds 
			\[
			\|D\mathscr{S}(\b \eta^1) - D\mathscr{S}(\b \eta^2)\|_{\mathcal{L}(L^\infty(\Omega;\mathbb{R}^2);V_0)} \leq C\|\b \eta^1 - \b \eta^2\|_{ L^\infty(\Omega;\mathbb{R}^2)},
			\]
			where the positive constant $C$ only depends on $a$ and $\b v$. 
		\end{prop}
		\begin{proof}
			The existence of a unique weak solution $u^* \in V_0$ to the problem follows again from the Lax--Milgram theorem. Indeed, we have that for all $u^*$ and $v \in V_0$
			\[
			(a(\b \eta) \nabla u^*, \nabla u^*)_H \geq a_* \| u^*\|_{V_0}^2, \qquad |(a(\b \eta) \nabla u^*, \nabla v)_H| \leq a^* \| u^*\|_{V_0}\| v\|_{V_0}
			\]
			and
			\[
			\begin{split}
				\left|\left(\left[\td{a}{\b s}(\b \eta) \cdot \b h\right]\nabla u - \left[\td{a}{\b s}(\b \eta) \cdot \b h\right]\b v, \nabla w\right)_H\right| &\leq \sup_{\b x \in \mathbb{R}^2}\left\|\td{a}{\b s}(\b x) \right\|_{\mathbb{R}^2}\|\b h\|_{L^\infty(\Omega;\mathbb{R}^2)}\left( \|u\|_{V_0} + \|\b v\|_{H} \right)\|w\|_{V_0}.
			\end{split}
			\]
			Hence, the map 
			\[
			w \mapsto\left(\left[\td{a}{\b s}(\b \eta) \cdot \b h\right]\nabla u - \left[\td{a}{\b s}(\b \eta) \cdot \b h\right]\b v, \nabla w\right)_H
			\]
			defines an element of $V_0^*$ and the problem is therefore well-posed. Next, we characterize the Fréchet derivative by exploiting its very definition. In particular, we prove that for any $\b \eta$ and $\b h \in L^\infty(\Omega;\mathbb{R}^2)$, with $\b h$ in a sufficiently small neighbourhood of zero, we have
			\[
			\|\mathscr{S}(\b \eta + \b h) - \mathscr{S}(\b \eta) - u^*\|_{V_0} = o(\|\b h\|_{L^\infty(\Omega;\mathbb{R}^2)}).
			\]
			For simplicity, let $u^{\b h} := \mathscr{S}(\b \eta + \b h)$ and $u := \mathscr{S}(\b \eta)$. Further define $v := u^{\b h} - u - u^*$. By definition of the quantities involved, for all $w \in V_0$, we have
			\begin{multline*}
				(a(\b \eta + \b h) \nabla u^{\b h} - a(\b \eta) \nabla u - a(\b \eta) \nabla u^*, \nabla w)_H = \left(\left[\td{a}{\b s}(\b \eta) \cdot \b h\right]\nabla u, \nabla w\right)_H \\ + \left(\left(a(\b \eta + \b h) - a(\b \eta) - \td{a}{\b s}(\b \eta) \cdot \b h\right)\b v, \nabla w\right)_H,
			\end{multline*}
			and by summing and subtracting $a(\b \eta)\nabla u^{\b h}$ at left hand side we arrive at
			\begin{multline} \label{eq:frechet1}
				(a(\b \eta) \nabla v, \nabla w)_H = - \left(\left[ a(\b \eta + \b h) - a(\b \eta) - \td{a}{\b s}(\b \eta) \cdot \b h\right] \nabla u^{\b h}, \nabla w\right)_H \\ - \left(\left[\td{a}{\b s}(\b \eta) \cdot \b h\right]\nabla (u ^{\b h} - u), \nabla w\right)_H + \left(\left(a(\b \eta + \b h) - a(\b \eta) - \td{a}{\b s}(\b \eta) \cdot \b h\right)\b v, \nabla w\right)_H.
			\end{multline}
			Choose now $w = v$. \reviewF{Since $a \in C^2(\mathbb R^2)$, we have that, }for $\b h$ in a sufficiently small neighbourhood of 0,
			\[
			\left\|a(\b \eta + \b h) - a(\b \eta) - \td{a}{\b s}(\b \eta) \cdot \b h \right\|_{L^\infty(\Omega)}= o(\|\b h\|_{L^\infty(\Omega;\mathbb{R}^2)}),
			\]
			and we can estimate the right hand side of \eqref{eq:frechet1} as follows:
			\[
			\begin{split}
				& \left| \left(\left[ a(\b \eta + \b h) - a(\b \eta) - \td{a}{\b s}(\b \eta) \cdot \b h\right] \nabla u^{\b h}, \nabla v\right)_H + \left(\left[\td{a}{\b s}(\b \eta) \cdot \b h\right]\nabla (u ^{\b h} - u), \nabla v\right)_H \right|\\
				& \hspace{4cm} + \left| \left(\left(a(\b \eta + \b h) - a(\b \eta) - \td{a}{\b s}(\b \eta) \cdot \b h\right)\b v, \nabla v\right)_H\right| \\
				& \hspace{1cm} \leq o(\|\b h\|_{L^\infty(\Omega;\mathbb{R}^2)})\|u^{\b h}\|_{V_0}\|v\|_{V_0} + C\|\b h\|_{L^\infty(\Omega;\mathbb{R}^2)}\|u^{\b h}-u\|_{V_0}\|v\|_{V_0} \\
				& \hspace{1cm} = o(\|\b h\|_{L^\infty(\Omega;\mathbb{R}^2)})\reviewF{\|v\|_{V_0}}
			\end{split}
			\]
			by virtue of Propositions \ref{prop:existence1} and \ref{prop:stability}. We are left to show that $D\mathscr{S}$ is also continuous. This follows essentially the same steps of Proposition \ref{prop:stability} on the solution map of the problem whose solution is $u^*$. Let $\b \eta^1$ and $\b \eta^2$ be elements of $L^\infty(\Omega;\mathbb{R}^2)$. For convenience, let us define $u^1 := D\mathscr{S}(\b \eta^1)[\b h]$ and $u^2 := D\mathscr{S}(\b \eta^2)[\b h]$ for some fixed $ \b h \in L^\infty(\Omega;\mathbb{R}^2)$. Then, it holds
			\begin{multline*}
				(a(\b \eta^1) \nabla u^1 - a(\b \eta^2) \nabla u^2, \nabla w)_H + \left(\left[\td{a}{\b s}(\b \eta^1) \cdot \b h \right] \nabla \mathscr{S}(\b \eta^1) - \left[\td{a}{\b s}(\b \eta^2)  \cdot \b h \right] \nabla \mathscr{S}(\b \eta^2) , \nabla w\right)_H \\ = \left(\left( \left[\td{a}{\b s}(\b \eta^1) - \td{a}{\b s}(\b \eta^2) \right]\cdot \b h \right)\b v, \nabla w\right)_H
			\end{multline*}
			for all $w \in V_0$. By suitably adding and subtracting mixed terms, and employing the mean value theorem, we end up with
			\begin{multline*}
				(a(\b \eta^2)\nabla(u^1-u^2), \nabla w)_H \\ \leq |([a(\b \eta^1)-a(\b \eta^2)]\nabla u^1, \nabla w)_H| + \left|\left( \left( \left[\td{a}{\b s}(\b \eta^1) - \td{a}{\b s}(\b \eta^2) \right]\cdot \b h \right)\nabla \mathscr{S}(\b \eta^1), \nabla w\right)_H \right|\\ + \left|\left(\left[\td{a}{\b s}(\b \eta^2)  \cdot \b h \right] \nabla(\mathscr{S}(\b \eta^1) - \mathscr{S}(\b \eta^2)) , \nabla w\right)_H \right|\\+ \sup_{\b x \in \mathbb{R}^2} \left\|\td{^2 a}{\b s^2}(\b x)\right \|_{\mathbb{R}^{2\times 2}}\|\b h\|_{ L^\infty(\Omega;\mathbb{R}^2)}\|\b v\|_{H}\|w\|_{V_0}\|\b \eta^1 - \b \eta^2\|_{ L^\infty(\Omega;\mathbb{R}^2)}.
			\end{multline*}
			Choose now $w = u^1-u^2$. Developing the estimates we have
			\begin{multline*}
				a_*\|u^1-u^2\|^2_{V_0} \leq \sup_{\b x \in \mathbb{R}^2}\left\| \td{a}{\b s}(\b x)\right\|_{\mathbb{R}^2}\|u^1\|_{V_0}\|u^1-u^2\|_{V_0}\|\b \eta^1 - \b \eta^2\|_{ L^\infty(\Omega;\mathbb{R}^2)} \\
				+ \sup_{\b x \in \mathbb{R}^2}\left\| \td{a}{\b s}(\b x)\right\|_{\mathbb{R}^2}\|\b h\|_{ L^\infty(\Omega;\mathbb{R}^2)}\|\mathscr{S}(\b \eta^1) - \mathscr{S}(\b \eta^2)\|_{L^\infty(\Omega;\mathbb{R}^2)}\|u^1-u^2\|_{V_0} \\ + \sup_{\b x \in \mathbb{R}^2} \left\|\td{^2 a}{\b s^2}(\b x)\right \|_{\mathbb{R}^{2\times 2}}\|\b h\|_{ L^\infty(\Omega;\mathbb{R}^2)}\left( \|\b v\|_{H} + \|\mathscr{S}(\b \eta^1)\|_{V_0} \right)\|u^1-u^2\|_{V_0}\|\b \eta^1 - \b \eta^2\|_{ L^\infty(\Omega;\mathbb{R}^2)}
			\end{multline*}
			and owing to all the previous results and the Lax--Milgram estimate
			\[
			\|u^i\|_{V_0} \leq a_*^{-1}\sup_{\b x \in \mathbb{R}^2}\left\|\td{a}{\b s}(\b x) \right\|_{\mathbb{R}^2}\|\b h\|_{L^\infty(\Omega;\mathbb{R}^2)}\left( \|u\|_{V_0} + \|\b v\|_{H} \right)
			\] 
			for $i = 1,\,2$, we end up with
			\[
			\dfrac{\|u^1-u^2\|_{V_0}}{\|\b h\|_{ L^\infty(\Omega;\mathbb{R}^2)}} \leq C\|\b \eta^1 - \b \eta^2\|_{ L^\infty(\Omega;\mathbb{R}^2)},
			\]
			where $C$ only depends on $a$ and $\b v$. This also means that
			\[
			\|D\mathscr{S}(\b \eta^1) - D\mathscr{S}(\b \eta^2)\|_{\mathcal{L}(L^\infty(\Omega;\mathbb{R}^2);V_0)} \leq C\|\b \eta^1 - \b \eta^2\|_{ L^\infty(\Omega;\mathbb{R}^2)},
			\]
			proving the claim. 
		\end{proof} \noindent
		The following technical result will be useful in the analysis of the longtime behavior of the system of interest.
		\begin{prop} \label{prop:compact}
			\reviewF{Assume there exists $\alpha > d$ such that} $\b \eta \in W^{1,\reviewF{\alpha}}(\Omega;\mathbb{R}^2)$ and $\b v \in W^{1,\alpha}(\Omega;\mathbb{R}^d)$. Then, the following hold:
			\begin{enumerate}[(i)]
				\item $\mathscr{S}(\b \eta) \in W^{2,q}(\Omega)$ for all $q \in (\reviewF{1},d)$;
				\item \reviewF{the linear map $D\mathscr{S}(\b \eta): L^\infty(\Omega;\mathbb R^2) \to V_0$ is compact for all $\b \eta \in W^{1,\alpha}(\Omega;\mathbb{R}^2)$}.
			\end{enumerate}
		\end{prop}
		\begin{proof}
			By the assumptions on $a$, it is immediate to show that $a \circ \b \eta \in \reviewF{W^{1,\alpha}(\Omega)} \hookrightarrow W^{1,d}(\Omega)$. Therefore, we also have that $a \circ \b \eta$ lies in the space of vanishing mean oscillation functions $\operatorname{VMO}(\Omega)$ (see, for the definition of the space and for the proof of this result, \cite[Section II.1]{BrezNir96}). Now, problem \eqref{eq:elliptic} is an elliptic problem with VMO coefficients. Moreover, as $a(\b \eta)\b v \in W^{1,\alpha}(\Omega;\mathbb{R}^d)$ and
			\[
			\div \begin{bmatrix}
				a(\b \eta) & 0 & 0 \\
				0 & a(\b \eta) & 0 \\
				0 & 0 & a(\b \eta) 
			\end{bmatrix} = \nabla (a(\b \eta)) \in L^\reviewF{\alpha}(\Omega),
			\]
			the first claim follows from applying \cite[Theorem 2.2]{Kang19}. Next, we address the compactness claim. In particular, we shall perform a comparison argument. Let $\b h \in L^\infty(\Omega;\mathbb{R}^2) \cap V$. From Proposition \ref{prop:frechet1}, we have \reviewF{that $D \mathscr{S}(\b \eta)[\b h] \in V_0$ and that}
			\[
			-\div\left(a(\b \eta) \nabla D \mathscr{S}(\b \eta)[\b h]+\left[\td{a}{\b s}(\b \eta) \cdot \b h\right]\left( \nabla \mathscr{S}(\b \eta) - \b v\right) \right)= 0
			\]
			\reviewF{in the weak sense. Therefore,}
			\begin{multline} \label{eq:comparison}
				\Delta D \mathscr{S}(\b \eta)[\b h] = \dfrac{1}{a(\b \eta)}   \left[\td{a}{\b s}(\b \eta) \cdot \b h\right]\left( \Delta \mathscr{S}(\b \eta) - \div \b v\right) + \dfrac{1}{a(\b \eta)} \left( \nabla \mathscr{S}(\b \eta) - \b v\right) \cdot \nabla\left[\td{a}{\b s}(\b \eta) \cdot \b h\right] \\ + \dfrac{1}{a(\b \eta)} \nabla\left[ a(\b \eta) \right] \cdot \nabla D \mathscr{S}(\b \eta)[\b h] .
			\end{multline}
			In the right hand side of \eqref{eq:comparison}, \reviewF{owing also to the first claim,} the first term clearly defines an element of $L^q(\Omega)$ for all $q \in (\review{1},d)$. \reviewF{As for the second term, we observe that by Sobolev embeddings we readily have $\mathscr{S}(\b \eta) \in W^{1,p}(\Omega)$ for any finite $p \geq 1$. Hence, for any fixed $s \in (\frac 65, 2)$,} by interpolation we have
			\begin{multline*}
				\left\| \dfrac{1}{a(\b \eta)} \left( \nabla \mathscr{S}(\b \eta) - \b v\right) \cdot \nabla\left[\td{a}{\b s}(\b \eta) \cdot \b h\right]  \right\|_{L^\reviewF{s}(\Omega)} \\ \leq C\|\nabla \mathscr{S}(\b \eta) - \b v\|_{L^\reviewF{\frac{2s}{2-s}}(\Omega;\mathbb{R}^d)}\left( \|\b h\|_{L^\infty(\Omega);\mathbb{R}^2)}\|\nabla \b \eta\|_{L^\reviewF{\alpha}(\Omega;\mathbb{R}^{2\times d})} + \|\b h\|_V\right).
			\end{multline*}
			\reviewF{
			Here, we are exploiting that, still by Sobolev embeddings,
			\[
			\nabla \mathscr{S}(\b \eta) - \b v \in L^{p}(\Omega;\mathbb{R}^d)
			\]
			for any finite $p \geq 1$. Finally, for the third term,}
			\begin{equation*}
				\left\| \dfrac{1}{a(\b \eta)} \nabla\left[ a(\b \eta) \right] \cdot \nabla D \mathscr{S}(\b \eta)[\b h] \right\|_{L^\reviewF{\frac{2\alpha}{\alpha+2}}(\Omega)}\leq C\|\nabla \b \eta\|_{L^\reviewF{\alpha}(\Omega;\mathbb{R}^d)}\|D \mathscr{S}(\b \eta)[\b h]\|_{V_0},
			\end{equation*}
			\reviewF{and, since $\alpha > d$, it holds
			\[
			\frac{2\alpha}{\alpha + 2} \in \left(\frac{2d}{d+2},\, 2\right) = \begin{cases}
				(1,2) & \quad \text{if }d=2,\\
				(\frac 65,2) \subset (1,3) & \quad \text{if }d=3.
			\end{cases}
			\]
			On account of the previous observations, we can choose $q = s = \frac{2\alpha}{\alpha + 2}$ and get}
			\[
			\Delta D \mathscr{S}(\b \eta)[\b h] \in L^\frac{2\alpha}{\alpha + 2}(\Omega)
			\]
			and, by elliptic regularity and the \reviewF{compact} embedding $W^{2, \frac{2\alpha}{\alpha + 2}}(\Omega) \hookrightarrow V$ in two and three dimensions, we \reviewF{conclude} that $D \mathscr{S}(\b \eta) : L^\infty(\Omega;\mathbb R^2) \to V_0$ is a linear and compact operator.
		\end{proof} \noindent
		The next result concerns the second Fréchet derivative of the solution map $\mathscr{S}$.
		\begin{prop} \label{prop:frechet2}
			The map \reviewF{$\mathscr{S}:L^\infty(\Omega;\mathbb{R}^2) \to V_0$} is twice Fréchet differentiable. For any $\b \eta$, $\b h$ and $\b k$ in $L^\infty(\Omega;\mathbb{R}^2)$, the second derivative $D^2\mathscr{S}: L^\infty(\Omega;\mathbb{R}^2) \to \mathcal{L}(L^\infty(\Omega;\mathbb{R}^2); \mathcal{L}(L^\infty(\Omega;\mathbb{R}^2);V_0))$ is \reviewF{given} by
			\[
			D^2\mathscr{S}(\b \eta)[\b k, \b h] = u^\sharp,
			\]
			where $u^\sharp \in V_0$ is the unique weak solution of the problem
			\[
			\footnotesize
			\begin{cases}
				-\div\left(a(\b \eta) \nabla u^\sharp + \left[\td{a}{\b s}(\b \eta) \cdot \b h \right] \nabla u^*_1 + \left[\td{a}{\b s}(\b \eta) \cdot \b k\right]\nabla u^*_2 + \left[\td{^2a}{\b s^2}(\b \eta) \b k \cdot \b h\right]\nabla u\right) = -\div\left( \left[\td{^2a}{\b s^2}(\b \eta) \b k \cdot \b h\right]\b v\right) & \quad \text{in } \Omega, \\
				u^\sharp = 0 & \quad \text{on }\partial\Omega,
			\end{cases}
			\]
			provided that $u := \mathscr{S}(\b \eta)$, $u^*_1 := D\mathscr{S}(\b \eta)[\b k]$ and $u^*_2 := D\mathscr{S}(\b \eta)[\b h]$.
		\end{prop}
		\begin{proof}
			The existence of a unique weak solution $u^\sharp \in V_0$ follows once again by means of the Lax--Milgram theorem. Indeed, we have that  for all $u^\sharp$ and $v \in V_0$
			\[
			(a(\b \eta) \nabla u^\sharp, \nabla u^\sharp)_H \geq a_* \| u^\sharp\|_{V_0}^2, \qquad |(a(\b \eta) \nabla u^\sharp, \nabla v)_H| \leq a^* \| u^\sharp\|_{V_0}\|v\|_{V_0}
			\]
			and,
			\begin{multline*}
				\left|\left(\left[\td{a}{\b s}(\b \eta) \cdot \b h \right] \nabla u^*_1 + \left[\td{a}{\b s}(\b \eta) \cdot \b k\right]\nabla u^*_2 + \left[\td{^2a}{\b s^2}(\b \eta) \b k \cdot \b h\right]\nabla u - \left[\td{^2a}{\b s^2}(\b \eta) \b k \cdot \b h\right]\b v, \nabla w\right)_H\right| \\ \leq \left[ \left(\| \b h\|_{L^\infty(\Omega;\mathbb{R}^2)} + \| \b k\|_{L^\infty(\Omega;\mathbb{R}^2)} \right)\sup_{\b x \in \mathbb{R}^2} \left\| \td{a}{\b s}(\b x) \right\|_{\mathbb{R}^2} + \| \b h\|_{L^\infty(\Omega;\mathbb{R}^2)}\| \b k\|_{L^\infty(\Omega;\mathbb{R}^2)}\sup_{\b x \in \mathbb{R}^2} \left\| \td{^2a}{\b s^2}(\b x) \right\|_{\mathbb{R}^{2\times2}} \right]\\ \times\left(\|u^*_1\|_{V_0} + \|u^*_2\|_{V_0} + \|u\|_{V_0} + \|\b v\|_H \right)\|w\|_{V_0}.
			\end{multline*}
			Therefore, the map 
			\[
			w \mapsto \left(\left[\td{a}{\b s}(\b \eta) \cdot \b h \right] \nabla u^*_1 + \left[\td{a}{\b s}(\b \eta) \cdot \b k\right]\nabla u^*_2 + \left[\td{^2a}{\b s^2}(\b \eta) \b k \cdot \b h\right]\nabla u - \left[\td{^2a}{\b s^2}(\b \eta) \b k \cdot \b h\right]\b v, \nabla w\right)_H
			\]
			defines an element of $V_0^*$ and hence the problem is well posed. Next, we invoke the definition of Fréchet derivative and aim to show that for any $\b \eta$, $\b h$ and $ \b k \in L^\infty(\Omega;\mathbb{R}^2)$, with $\b k$ in a sufficiently small neighbourhood of zero, we have
			\[
			\|D\mathscr{S}(\b \eta + \b k)[\b h] - D\mathscr{S}(\b \eta)[\b h] - u^\sharp\|_{V_0} = o(\|\b k\|_{L^\infty(\Omega;\mathbb{R}^2)}).
			\]
			Let us introduce convenient notation by setting $u^{\b k} := D\mathscr{S}(\b \eta + \b k)[\b h]$ and $u^0 := D\mathscr{S}(\b \eta)[\b h]$, while define the remainder $v := u^{\b k} - u^0 - u^\sharp$. By Proposition \ref{prop:frechet1}, for all $w \in V_0$, we have
			\begin{multline*}
				(a(\b \eta + \b k) \nabla u^{\b k} - a(\b \eta) \nabla u^0 - a(\b \eta) \nabla u^\sharp, \nabla w)_H \\
				+ \left( \left[ \td{a}{\b s}(\b \eta + \b k) \cdot \b h\right]\nabla \mathscr{S}(\b \eta + \b k) - \left[ \td{a}{\b s}(\b \eta) \cdot \b h\right]\nabla \mathscr{S}(\b \eta) - \left[ \td{^2a}{\b s^2}(\b \eta)\b k \cdot \b h\right]\nabla \mathscr{S}(\b \eta), \nabla w\right)_H \\ 
				- \left(  \left[\td{a}{\b s}(\b \eta) \cdot \b h \right] \nabla D\mathscr{S}(\b \eta)[\b k] + \left[\td{a}{\b s}(\b \eta) \cdot \b k\right]\nabla u^0, \nabla w\right)_H\\= \left( \left( \left[  \td{a}{\b s}(\b \eta + \b k) - \td{a}{\b s}(\b \eta) - \td{^2a}{\b s^2}(\b \eta)\b k\right] \cdot \b h\right)\b v, \nabla w\right)_H.
			\end{multline*}
			Let us tackle the various terms separately. First, by summing and subtracting $a(\b \eta)\nabla u^{\b k}$, we have
			\begin{equation} \label{eq:frechet21}
				(a(\b \eta + \b k) \nabla u^{\b k} - a(\b \eta) \nabla u^0 - a(\b \eta) \nabla u^\sharp, \nabla w)_H = (a(\b \eta)\nabla v, \nabla w)_H + ((a(\b \eta + \b k)  - a(\b \eta))\nabla u^{\b k}, \nabla w)_H.
			\end{equation} 
			Next, we have
			\begin{equation} \label{eq:frechet22}
				\begin{split}
					& \left( \left[ \td{a}{\b s}(\b \eta + \b k) \cdot \b h\right]\nabla \mathscr{S}(\b \eta + \b k) - \left[ \td{a}{\b s}(\b \eta) \cdot \b h\right]\nabla \mathscr{S}(\b \eta) - \left[ \td{^2a}{\b s^2}(\b \eta)\b k \cdot \b h\right]\nabla \mathscr{S}(\b \eta), \nabla w\right)_H \\
					& \hspace{1cm} = \left( \left( \left[ \td{a}{\b s}(\b \eta + \b k) - \td{a}{\b s}(\b \eta) - \td{^2a}{\b s^2}(\b \eta)\b k \right]\cdot \b h\right)\nabla \mathscr{S}(\b \eta), \nabla w\right)_H \\
					& \hspace{3cm} + \left(\left[\td{a}{\b s}(\b \eta + \b k) \cdot \b h\right](\nabla \mathscr{S}(\b \eta + \b k)-\nabla \mathscr{S}(\b \eta)) , \nabla w\right)_H.
				\end{split}
			\end{equation}
			On account of \eqref{eq:frechet21} and \eqref{eq:frechet22}, and choosing $w = v$, we have
			\begin{multline*}
				(a(\b \eta)\nabla v, \nabla v)_H + \left((a(\b \eta + \b k)  - a(\b \eta))\nabla u^{\b k} - \left[\td{a}{\b s}(\b \eta) \cdot \b k\right]\nabla u^0, \nabla v\right)_H \\
				\left( \left( \left[ \td{a}{\b s}(\b \eta + \b k) - \td{a}{\b s}(\b \eta) - \td{^2a}{\b s^2}(\b \eta)\b k \right]\cdot \b h\right)\nabla \mathscr{S}(\b \eta), \nabla v\right)_H \\+ \left(\left[\td{a}{\b s}(\b \eta + \b k) \cdot \b h\right](\nabla \mathscr{S}(\b \eta + \b k)-\nabla \mathscr{S}(\b \eta)) -  \left[ \td{a}{\b s}(\b \eta) \cdot \b h\right]\nabla D\mathscr{S}(\b \eta)[\b k], \nabla v\right)_H
				\\= \left( \left( \left[  \td{a}{\b s}(\b \eta + \b k) - \td{a}{\b s}(\b \eta) - \td{^2a}{\b s^2}(\b \eta)\b k\right] \cdot \b h\right)\b v, \nabla v\right)_H.
			\end{multline*}
			In order to tackle the second and fourth terms in the above, we proceed as follows. For the former, we get
			\begin{equation} \label{eq:frechet23}
				\begin{split}
					&\left((a(\b \eta + \b k)  - a(\b \eta))\nabla u^{\b k} - \left[\td{a}{\b s}(\b \eta) \cdot \b k\right]\nabla u^0, \nabla v\right)_H \\
					&\hspace{2cm}=\left(\left(a(\b \eta + \b k)  - a(\b \eta) - \td{a}{\b s}(\b \eta) \cdot \b k\right)\nabla u^{\b k} + \left[\td{a}{\b s}(\b \eta) \cdot \b k\right]\nabla (u^{\b k} - u^0), \nabla v\right)_H 
				\end{split}
			\end{equation}
			while
			\begin{equation} \label{eq:frechet24} \small
				\begin{split}
					&\left(\left[\td{a}{\b s}(\b \eta + \b k) \cdot \b h\right](\nabla \mathscr{S}(\b \eta + \b k)-\nabla \mathscr{S}(\b \eta)) -  \left[ \td{a}{\b s}(\b \eta) \cdot \b h\right]\nabla  D\mathscr{S}(\b \eta)[\b k], \nabla v\right)_H \\
					&\hspace{1cm}=\left(\left[\td{a}{\b s}(\b \eta + \b k) \cdot \b h\right](\nabla \mathscr{S}(\b \eta + \b k)-\nabla \mathscr{S}(\b \eta)- \nabla  D\mathscr{S}(\b \eta)[\b k]), \nabla v \right)_H \\
					& \hspace{7cm}+  \left( \left[ \left( \td{a}{\b s}(\b \eta + \b k) - \td{a}{\b s}(\b \eta) \right) \cdot \b h\right]\nabla  D\mathscr{S}(\b \eta)[\b k], \nabla v\right)_H.
				\end{split}
			\end{equation}
			We are now in a position to perform the needed estimates. Observing that
			\begin{equation} \label{eq:frechet25}
				(a(\b \eta)\nabla v, \nabla v)_H \geq a_*\|v\|^2_{V_0}
			\end{equation}
			and in light of the regularity of $a$, we have, from \eqref{eq:frechet23} and Proposition \ref{prop:frechet1}, that
			\begin{equation*}
				\begin{split}
					& \left|\left(\left(a(\b \eta + \b k)  - a(\b \eta) - \td{a}{\b s}(\b \eta) \cdot \b k\right)\nabla u^{\b k} + \left[\td{a}{\b s}(\b \eta) \cdot \b k\right]\nabla (u^{\b k} - u^0), \nabla v\right)_H\right| \\
					& \hspace{2cm}\leq o(\|\b k\|_{L^\infty(\Omega;\mathbb{R}^2)})\|u^{\b k}\|_{V_0}\|v\|_{V_0} + C\| \b k\|_{L^\infty(\Omega;\mathbb{R}^2)}\|u^{\b k} - u^0\|_{V_0}\|v\|_{V_0} \\
					& \hspace{2cm} \leq C\left[o(\|\b k\|_{L^\infty(\Omega;\mathbb{R}^2)})\left(\|\b \eta\|_{L^\infty(\Omega;\mathbb{R}^2)} + \|\b k\|_{L^\infty(\Omega;\mathbb{R}^2)} \right) + \|\b k\|_{L^\infty(\Omega;\mathbb{R}^2)}^2 \right]\|v\|_{V_0}.
				\end{split}
			\end{equation*}
			Moreover, owing also to Proposition \ref{prop:existence1}, we have
			\[
			\begin{split}
				& \left|\left( \left( \left[ \td{a}{\b s}(\b \eta + \b k) - \td{a}{\b s}(\b \eta) - \td{^2a}{\b s^2}(\b \eta)\b k \right]\cdot \b h\right)\nabla \mathscr{S}(\b \eta), \nabla v\right)_H \right| \leq Co(\|\b k\|_{L^\infty(\Omega;\mathbb{R}^2)})\|\b h\|_{L^\infty(\Omega;\mathbb{R}^2)}\|v\|_{V_0}
			\end{split}
			\]
			and analogously
			\[
			\begin{split}
				& \left| \left( \left( \left[  \td{a}{\b s}(\b \eta + \b k) - \td{a}{\b s}(\b \eta) - \td{^2a}{\b s^2}(\b \eta)\b k\right] \cdot \b h\right)\b v, \nabla v\right)_H  \right| \leq o(\|\b k\|_{L^\infty(\Omega;\mathbb{R}^2)})\|\b h\|_{L^\infty(\Omega;\mathbb{R}^2)}\| \b v\|_{\b H}\|v\|_{V_0}.
			\end{split}
			\]
			Finally, we have
			\[
			\begin{split}
				& \left| \left(\left[\td{a}{\b s}(\b \eta + \b k) \cdot \b h\right](\nabla \mathscr{S}(\b \eta + \b k)-\nabla \mathscr{S}(\b \eta)) -  \left[ \td{a}{\b s}(\b \eta) \cdot \b h\right]\nabla D\mathscr{S}(\b \eta)[\b k], \nabla v\right)_H \right|\\
				& \hspace{1cm} \leq  \left| \left(\left[\td{a}{\b s}(\b \eta + \b k) \cdot \b h\right](\nabla \mathscr{S}(\b \eta + \b k)-\nabla \mathscr{S}(\b \eta) -  \nabla D\mathscr{S}(\b \eta)[\b k]), \nabla w \right)_H\right| \\
				& \hspace{4cm}+ \left| \left( \left[ \td{a}{\b s}(\b \eta + \b k) \cdot \b h-  \td{a}{\b s}(\b \eta) \cdot \b h\right]\nabla D\mathscr{S}(\b \eta)[\b k], \nabla v\right)_H \right| \\
				& \hspace{1cm} \leq C\left[ \|\b h\|_{L^\infty(\Omega;\mathbb{R}^2)}\|D\mathscr{S}(\b \eta)\|_{\mathcal{L}(L^\infty(\Omega;\mathbb{R}^2);V_0)}\|\b k\|_{L^\infty(\Omega;\mathbb{R}^2)}^2\right]\|v\|_{V_0},
			\end{split}
			\]
			and collecting all of the above yields
			\[
			\|v\|_{V_0} \leq C(\b h, a, \b \eta)o(\|\b k\|_{L^\infty(\Omega;\mathbb{R}^2)})
			\]
			where the positive constant $C$ only depends on $\b h$, $a$ and $\b \eta$. The claim is proved.
		\end{proof} 
		\begin{prop} \label{prop:analytic}
			The solution map $\mathscr{S}$ is locally analytic around points of $L^\infty(\Omega;\mathcal{U}_0)$.
		\end{prop}
		\begin{proof}
			The argument exploits the analytic version of the implicit function theorem. Define the operator $\mathscr{F}: L^\infty(\Omega;\mathbb{R}^2) \times V_0 \to V_0^*$ so that
			\[
			\mathscr{F}(\b \eta, \xi) := - \div(a(\b \eta)(\nabla \xi - \b v)) \qquad \forall \: \b \eta \in L^\infty(\Omega;\mathbb{R}^2),\, \xi \in V_0,
			\]
			to be interpreted in the weak sense. Of course, we have $\mathscr{F}(\b \eta, \mathscr{S}(\b \eta)) = 0$ for all $\b \eta \in L^\infty(\Omega;\mathbb{R}^2)$ by the definition of the solution map. Observe that the operator $\mathscr{F}$ is affine-linear in the variable $\xi$, and the partial Fréchet derivative $D_\xi \mathscr{F}(\b \eta, \xi)$ defines a linear isomorphism between $V_0$ and $V_0^*$ by the Lax--Milgram theorem (here we also exploit Assumption \ref{hyp:appendix}). If $\b \eta \in L^\infty(\Omega; \mathcal{U}_0)$ \reviewF{(see Assumption \ref{hyp:appendix})}, as $\mathcal{U}_0$ is compactly contained in $\mathcal{U}$, we conclude that $\mathcal{F}$ is affine-bilinear, and hence analytic. This is enough to conclude the argument by means of the analytic version of the implicit function theorem (see, for instance, \cite[Chapter 4, Theorem 4.B]{Zeidler}).
		\end{proof} \noindent
		\textbf{Acknowledgments.} Part of this work was done when ADP was visiting HA and HG at the Department of Mathematics at the University of Regensburg. ADP wishes to thank HA and HG for their kind hospitality. The second author is a member of Gruppo Nazionale per l'Analisi Matematica, la Probabilità e le loro Applicazioni
		(GNAMPA), Istituto Nazionale di Alta Matematica  (INdAM). The second author has been partially funded by MIUR-PRIN
		Grant 2020F3NCPX  ``Mathematics  for industry 4.0 (Math4I4)'', and is supported by the MUR Grant ``Dipartimento di Eccellenza 2023-2027''. \review{The work of all three authors is supported by the DFG Graduiertenkolleg 2339 \textit{IntComSin} -- Project-ID 321821685.}
		\printbibliography
\end{document}